\documentclass{article}
\usepackage{amssymb}
\usepackage{amsfonts}
\usepackage{amsmath}
\usepackage{geometry}

\newtheorem{theorem}{Theorem}
\newtheorem{acknowledgement}[theorem]{Acknowledgement}

\newtheorem{corollary}[theorem]{Corollary}

\newtheorem{lemma}[theorem]{Lemma}
\newtheorem{notation}[theorem]{Notation}

\newtheorem{proposition}[theorem]{Proposition}
\newtheorem{remark}[theorem]{Remark}

\newenvironment{proof}[1][Proof]{\noindent\textbf{#1.} }{\ \rule{0.5em}{0.5em}}
\input{tcilatex}

\begin{document}

\title{On Uniformly Subelliptic Operators and Stochastic Area}
\author{Peter Friz\thanks{%
Corresponding author. Department of Pure Mathematics and Mathematical
Statistics, University of Cambridge. Email: P.K.Friz@statslab.cam.ac.uk. } \
\ \thanks{%
Leverhulme Fellow.} \and Nicolas Victoir}
\maketitle

\begin{abstract}
Let $X^{a}$ be a Markov process with generator $\sum_{i,j}\partial
_{i}\left( a^{ij}\partial _{j}\cdot \right) $ where $a$ is a uniformly
elliptic symmetric matrix. Thanks to the fundamental works of T. Lyons,
stochastic differential equations driven by $X^{a}$ can be solved in the
"rough path sense"; that is, pathwise by using a suitable stochastic area
process.

Our construction of the area, which generalizes previous works of
Lyons-Stoica and then Lejay, is based on Dirichlet forms associated to
subellitpic operators. This enables us in particular to discuss large
deviations and support descriptions in suitable rough path topologies. As
typical rough path corollary, Freidlin-Wentzell theory and the
Stroock-Varadhan support theorem remain valid for stochastic differential
equations driven by $X^{a}.$
\end{abstract}

\section{Introduction}

Let $V=\left( V_{1},...,V_{d}\right) $ be a collection of sufficiently nice
vector fields on $\mathbb{R}^{e}$ and consider the stochastic differential
equation in the Stratonovich sense $dY=V\left( Y\right) dB,\,\,\,Y\left(
0\right) =y_{0}\in \mathbb{R}^{e}$, driven by a $d$-dimensional Brownian
motion, a diffusion with generator $\frac{1}{2}\sum_{i=1}^{d}\partial
_{i}^{2}$. \ We try to understand what happens when $B$ is replaced by a $d$%
-dimensional diffusion process $X=X^{a}$ with uniformly elliptic generator
in divergence form $\sum_{i,j=1}^{d}\partial _{i}\left( a^{ij}\partial
_{j}\cdot \right) $. Of course, $dY=V\left( Y\right) dX$ still makes sense
as Stratonovich equation if $a$ is smooth but this breaks down when $a$ is
only assumed to be measurable. Such an assumption is not only standard in
the theory of partial differential equations but also a basic example in the
theory of Dirichlet forms \cite{Fu94} and the construction of the
corresponding diffusion process $X^{a}$ is well-known, e.g. \cite{St88, Fu94}%
.

We recall that one can construct $X^{a}$ as weak limit of semi-martingales $%
X^{a\left( \varepsilon \right) }$ along a sequence of mollifier
approximations $\left\{ a\left( \varepsilon \right) :\varepsilon >0\right\} $%
. It is a natural question \cite{LejII06} if the sequence of SDE\ solutions
driven by $X^{a\left( \varepsilon \right) }$ converges. One can also replace 
$X^{a}$ by piecewise linear approximations $X^{a}\left( n\right) $ and ask
if the resulting ODE solutions converge. It turns out they all converge to
the same limiting object which can be constructed intrinsically as solution
to the \textit{rough differential equation} \cite{lyons-98, lyons-qian-02}
of form $dY=V\left( Y\right) d\mathbf{X}.$A stochastic area process $A^{a}$
is now considered part of the driving signal $\mathbf{X}=\left(
X^{a},A^{a}\right) $. The construction of $A^{a}$ was carried out by subtle
forward-backward martingale arguments in \cite{LySt99}, together with a
convergence statement for piecewise linear approximations. It is verified in 
\cite{LejI06} that convergence takes place in suitable rough path metrics .
By the fundamental continuity result of rough path theory this implies the
convergence of ODE solutions driven by $X^{a}\left( n\right) $, i.e. a
Wong-Zakai theorem.

In contrast to \cite{LySt99, LejI06, LejII06} we emphasize and exploit the
Markovian nature of $\left( X^{a},A^{a}\right) $\textit{. }The basic
observation is that for smooth $a$ we are dealing with semi-martingales $%
X^{a}$ so that the stochastic area process should be given in terms of It%
\^{o} stochastic integrals,%
\begin{equation*}
t\mapsto A_{t}^{a}\equiv \frac{1}{2}\int_{0}^{t}\left( X^{a}\otimes
dX^{a}-dX^{a}\otimes X^{a}\right) \in so\left( d\right) .
\end{equation*}%
It is a simple exercise in It\^{o} calculus\footnote{%
Once can proceed as follows. First write $X=X^{a}$ as solution to a
Stratonovich SDE involving a smooth square-root of $a$. In combination with
the fact the the lift of $X$, denoted by $Y$ say, is obtained by solving the
Stratonovich equation $dY=\sum_{i=1}^{d}U_{i}\left( Y\right) \circ dX^{i}$
along the left-invariant vectorfields $U_{1},...,U_{d}$ on $g_{2}\left( 
\mathbb{R}^{d}\right) $ as defined in (\ref{DefOf_U_on_g2}), a few lines of
It\^{o} calculus identify the generator of the lift.} to see that the\
process $\left( X^{a},A^{a}\right) $ is Markov with (uniformly subelliptic)
generator of form%
\begin{equation}
L^{a}=\sum_{i,j=1}^{d}U_{i}\left( a^{ij}U_{j}\cdot \right) .
\label{Lsubelliptic}
\end{equation}

The vector fields $U_{1},...,U_{d}$ are defined in (\ref{DefOf_U_on_g2}) and
play the r\^{o}le of coordinate vector fields $\partial _{1},...,\partial
_{d}$ on $g^{2}\left( \mathbb{R}^{d}\right) \equiv \mathbb{R}^{d}\oplus
so\left( d\right) $, which is given the structure of a Lie group $G$. Of
course, $L^{a}$ is understood in a weak sense and the correct mathematical
object is the Dirichlet form\footnote{%
Lebesgue measure on $g^{2}\left( \mathbb{R}^{d}\right) $ coincides with Haar
measure $m$ on $G$. Then $U_{i}^{\ast }=-U_{i}$ where $^{\ast }$ denotes the
formal adjoint with respect to $m.$} 
\begin{equation}
\mathcal{E}^{a}\left( f,g\right)
=\sum_{i,j=1}^{d}\int_{G}dm\,a^{ij}U_{i}fU_{j}g.  \label{EaIntro}
\end{equation}

We can thus use the highly developed analytic machinery of Dirichlet forms 
\cite{Dav89, Fu94}; the collections of results in \cite{SturmIII}, in
conjunction with \cite{SaSt91}, applies directly to (\ref{EaIntro}). Leaving
precise references to those papers, the relevant results in \cite{SturmIII}
are based on the seminal works of De Giorgi, Nash, Moser for the elliptic
case and the various extensions to subelliptic/H\"{o}rmander type operators
as studied in papers by Rothschild, Stein, Jerison, S\'{a}nchez-Calle,
Nagel, Waigner and many others.

This paper is organized as follows. In Sections 2 and 3 we spezialise the
toolbox of Dirichlet forms to our situation and settle the notation. In
Section 4 we show that the $g^{2}\left( \mathbb{R}^{d}\right) $-valued
Markov process $\mathbf{X}^{a}$ has, just as Brownian motion and L\'{e}vy
area, $\left( 1/2-\varepsilon \right) $-H\"{o}lder regularity with respect
to Carnot-Caratheodory distance on $g^{2}\left( \mathbb{R}^{d}\right) $. It
follows that a.e. sample path $\mathbf{X}^{a}\left( \omega \right) $ is a
geometric H\"{o}lder rough path in the sense of Lyons, \cite{lyons-98,
friz-victoir-04-Note}. In fact, the H\"{o}lder norm of $\mathbf{X}^{a}$ is
seen to have Gaussian tail which answers a question raised in Lyons' St.
Flour lecture \cite{lyons-04}. In\ Section 5 we study both weak
approximations, $a_{n}\rightarrow a$ a.e. is seen to imply $\mathbf{X}%
^{a_{n}}\rightarrow \mathbf{X}^{a}$ in distribution, and a strong Wong-Zakai
type theorem. The latter shows that our stochastic area associated to $X^{a}$
coincides with the area constructed by Lyons and Stoica \cite{LySt99} and we
improve on results in \cite{LejI06, LejII06}. In Section 6 we note that an
RDE solution jointly with its driving signal $\mathbf{X}^{a}$ is Markov and
describe its generator, using stochastic Taylor expansions for random RDEs
obtained in \cite{friz-victoir-06-Euler}. In Section 7 we prove a sample
path large deviation principle for $\mathbf{X}^{a}$ making crucial use of Ram%
\'{\i}rez's result \cite{Ra01}. As a typical rough paths corollary, we
obtain Freidlin-Wentzell type large deviations for stochastic differential
equations driven by $\mathbf{X}^{a}$ in the rough path sense. Finally, in
Section 8 \ we revert to the case where $\mathbf{X}^{a}$ is the lift of $%
X^{a}$ (that is, $a$ is defined on $\mathbb{R}^{d}$ rather than $g^{2}\left( 
\mathbb{R}^{d}\right) $) and prove that $\mathbf{X}^{a}$ has full support in
suitable H\"{o}lder topologies. As a typical rough paths corollary, we
obtain a Stroock-Varadhan type support theorem for stochastic differential
equations driven by $\mathbf{X}^{a}$ in the rough path sense. Such a support
description was conjectured by T. Lyons in \cite{lyons-04}.

\begin{notation}
Although the key notations are introduced in the main text as appropriate we
feel the reader will be helped by this brief summary. The space of real
antisymmetric $d\times d$ matrices is denoted by $so\left( d\right) $ and is
given the standard Euclidean structure with $\cdot $ denoting the scalar
product. The corresponding norm is denoted by $\left\vert \cdot \right\vert $%
. It will cause no confusion to use $\cdot $ and $\left\vert \cdot
\right\vert $ also for standard scalar product and Euclidean norm on $%
\mathbb{R}^{d}$. The vector space $g^{2}\left( \mathbb{R}^{d}\right) =%
\mathbb{R}^{d}\oplus so\left( d\right) $ will be given a nilpotent Lie
algebra structure so that the corresponding Lie group can and will be
realized on the same space, $\left( g^{2}\left( \mathbb{R}^{d}\right) ,\ast
,0\right) $. Points in $g^{2}\left( \mathbb{R}^{d}\right) $ are denoted by $%
x,y,z$, ... and may be written out in coordinates as $\left( \left(
x^{1;i}\right) ,\left( x^{2;jk}\right) :i,j,k=1,...,d\text{ with }j<k\right) 
$. We also write $x^{1}=\pi _{1}\left( x\right) $, $x^{2}=\pi _{2}\left(
x\right) $ for the projections to $\mathbb{R}^{d},\,so\left( d\right) $
respectively. Haar measure on $g^{2}\left( \mathbb{R}^{d}\right) $ coincides
with Lebesgue measure on $\mathbb{R}^{d}\oplus so\left( d\right) $ and is
denoted by $m$, in integrations we write $dm,dm\left( x\right) $ or simply $%
dx$. We use $\left\langle \cdot ,\cdot \right\rangle $ for the scalar
product in $L^{2}\left( g^{2}\left( \mathbb{R}^{d}\right) ,m\right) $ and
the corresponding $L^{2}$-norm is written as $\left\vert \cdot \right\vert
_{L^{2}}$ or $\left\vert \cdot \right\vert _{L^{2}\left( D\right) }$ for $%
D\subset g^{2}\left( \mathbb{R}^{d}\right) $.The Lie group $g^{2}\left( 
\mathbb{R}^{d}\right) $ has a dilation structure $\delta _{\lambda }\left(
x\right) \mapsto \left( \lambda \pi _{1}\left( x\right) ,\lambda ^{2}\pi
_{2}\left( x\right) \right) $, carries a Carnot-Carathedory continuous norm $%
x\mapsto $ $\left\Vert x\right\Vert $, homogenuous in the sense that $%
\left\Vert \delta _{\lambda }\left( x\right) \right\Vert =\left\vert \lambda
\right\vert \left\Vert x\right\Vert ,$ and equivalent to $\left\vert \pi
_{1}\left( x\right) \right\vert +\left\vert \pi _{2}\left( x\right)
\right\vert ^{1/2}$. It induces the left invariant Carnot-Caratheodory
distance $d\left( x,y\right) =\left\Vert x^{-1}\ast y\right\Vert $under
which $g^{2}\left( \mathbb{R}^{d}\right) $ is a metric (in fact: geodesic)
space. This distance coincides with the intrinsic metric from a reference
Dirichlet form $\mathcal{E}$. A family of Dirichlet forms $\left\{ \mathcal{E%
}^{a}:a\in \Xi \left( \Lambda \right) \right\} $, where $\Xi \left( \Lambda
\right) $ denotes a class of certain diffusion matrices with ellipticity
constant $\Lambda $, gives rise to a family of intrinsic metrics on $%
g^{2}\left( \mathbb{R}^{d}\right) $, denoted by $d^{a}$, all Lipschitz
equivalent to $d$. Stochastic process with values in $g^{2}\left( \mathbb{R}%
^{d}\right) $ are denoted by capital bold letter such as $\mathbf{X}$ or $%
\mathbf{X}^{a,x}$ to indicate dependence on $a\in \Xi \left( \Lambda \right) 
$ and starting point. The $so\left( d\right) $-valued area process $A:=\pi
_{2}\left( \mathbf{X}\right) $ will be of interest. A fixed path in $C\left( %
\left[ 0,1\right] ,g^{2}\left( \mathbb{R}^{d}\right) \right) $ may be
written as $\mathbf{x=x}\left( \cdot \right) $ or $\omega $, the latter is
used when $C\left( \left[ 0,1\right] ,g^{2}\left( \mathbb{R}^{d}\right)
\right) $ is equipped with a Borel measure such as the law of $\mathbf{X}%
^{a,x}$ for which we write $\mathbb{P}^{a,x}$. $L^{p}$-norms with respect to 
$\mathbb{P}^{a,x}$ are denoted by $\left\Vert \cdot \right\Vert
_{L^{p}\left( \mathbb{P}^{a,x}\right) }$. A path $\mathbf{x}\in C\left( %
\left[ 0,1\right] ,g^{2}\left( \mathbb{R}^{d}\right) \right) $ has
increments $\mathbf{x}_{s,t}=\mathbf{x}_{s}^{-1}\ast \mathbf{x}_{t}=:\left( 
\mathbf{x}_{s,t}^{1},\mathbf{x}_{s,t}^{2}\right) $. Note $\mathbf{x}_{t}^{1}-%
\mathbf{x}_{s}^{1}=\mathbf{x}_{s,t}^{1}$ but $\mathbf{x}_{t}^{2}-\mathbf{x}%
_{s}^{2}\neq \mathbf{x}_{s,t}^{2}=\mathbf{x}_{t}^{2}-\mathbf{x}_{s}^{2}-%
\left[ \mathbf{x}_{s}^{1},\mathbf{x}_{s,t}^{1}\right] /2$. (Semi-)norms and
distances are defined naturally on this path space over $g^{2}\left( \mathbb{%
R}^{d}\right) $. In particular,%
\begin{equation*}
\left\Vert \mathbf{x}\right\Vert _{\alpha \text{-H\"{o}l}}=\sup_{0\leq
s<t\leq 1}\frac{d\left( \mathbf{x}_{s}\mathbf{,x}_{t}\right) }{\left\vert
t-s\right\vert ^{\alpha }}=\sup_{0\leq s<t\leq 1}\frac{\left\Vert \mathbf{x}%
_{s,t}\right\Vert }{\left\vert t-s\right\vert ^{\alpha }}\sim \sup_{0\leq
s<t\leq 1}\frac{\left\vert \mathbf{x}_{s,t}^{1}\right\vert +\left\vert 
\mathbf{x}_{s,t}^{2}\right\vert ^{1/2}}{\left\vert t-s\right\vert ^{\alpha }}%
.
\end{equation*}%
and%
\begin{equation*}
d_{\alpha \text{-H\"{o}l}}\left( \mathbf{x},\mathbf{y}\right) =\sup_{0\leq
s<t\leq 1}\frac{d\left( \mathbf{x}_{s,t}\mathbf{,y}_{s,t}\right) }{%
\left\vert t-s\right\vert ^{\alpha }}.
\end{equation*}%
We write $d_{0}\equiv d_{0\text{-H\"{o}l}}$ and $d_{\infty }\left( \mathbf{x}%
,\mathbf{y}\right) =\sup_{0\leq t\leq 1}d\left( \mathbf{x}_{t}\mathbf{,y}%
_{t}\right) $. Care must be taken since $d_{0}$ and $d_{\infty }$ are not
Lipschitz equivalent. We avoid the double bar notation, i.e. $\left\Vert
\cdot \right\Vert _{\left( ...\right) }$, for semi-norms resp. distances on
the path space over some Euclidean space $\mathbb{R}^{e},e\in \mathbb{N}$.
For instance, when $y\in C\left( \left[ 0,1\right] ,\mathbb{R}^{e}\right) $
we write%
\begin{equation*}
\left\vert y\right\vert _{\alpha \text{-H\"{o}l}}=\sup_{0\leq s<t\leq 1}%
\frac{\left\vert y_{t}-y_{s}\right\vert }{\left\vert t-s\right\vert ^{\alpha
}}=\sup_{0\leq s<t\leq 1}\frac{\left\vert y_{s,t}\right\vert }{\left\vert
t-s\right\vert ^{\alpha }}.
\end{equation*}%
Vector fields (usually on some Euclidean space $\mathbb{R}^{e},e\in \mathbb{N%
}$) are denoted by $V$ and usually assumed to be in some regularity class $%
\mathrm{Lip}^{p}$ which means bounded derivatives up to order $\lfloor
p\rfloor $, and H\"{o}lder regularity of the $\lfloor p\rfloor $th
derivative with exponent $p-\lfloor p\rfloor $. In particular, such vector
fields are bounded. The (smooth but unbounded) invariant vector fields on $%
g^{2}\left( \mathbb{R}^{d}\right) $ are denoted by $U_{i}$. A dissection $D$
of $[0,1]$ is a collection $\left\{
0=t_{0}<t_{1}<...<t_{\#D-1}<t_{\#D}=1\right\} $. Its mesh is defined as $%
\left\vert D\right\vert =\sup_{i=1,...,\#D}t_{i}-t_{i-1}$. Given $t\in \left[
0,1\right] $ we write $t_{D}$ for its lower neighbour in $D$ that is $%
t_{D}=\max \left\{ t_{i}\in D:t_{i}\leq t\right\} $. Similarly, $t^{D}$
denotes the upper neighbour in $D$. Constants which appears in statement are
typically indexed by the statement number. To indicate changing constant in
proofs we sometimes number them with upper indices. (This will cause no
confusion with powers.) We try to be explicit about the dependence of all
constant with the exception of $d=\dim \left( \mathbb{R}^{d}\right) .$
\end{notation}

\section{Analysis on the Group}

Let $g^{2}\left( \mathbb{R}^{d}\right) $ be the free step-$2$ nilpotent Lie
algebra over $\mathbb{R}^{d}$, that is $g^{2}\left( \mathbb{R}^{d}\right) =%
\mathbb{R}^{d}\oplus so\left( d\right) $ ($so\left( d\right) $ being the
space of antisymmetric $d\times d$ matrices) with Lie bracket%
\begin{equation*}
\left[ x,y\right] \equiv \left[ \left( x^{1},x^{2}\right) ,\left(
y^{1},y^{2}\right) \right] =x^{1}\otimes y^{1}-y^{1}\otimes x^{1}.
\end{equation*}%
Due to nilpotency and the Campbell-Baker-Hausdorff formula, we can and will
realize the associated Lie group on the same space $g^{2}\left( \mathbb{R}%
^{d}\right) =\mathbb{R}^{d}\oplus so\left( d\right) $ with product%
\begin{equation*}
x\ast y=x+y+\frac{1}{2}\left[ x,y\right]
\end{equation*}%
and unit element $0$. Lebesgue-measure $dx$ on $\mathbb{R}^{d}\oplus
so\left( d\right) $ is the (left- and right-invariant) Haar measure $m;$ in
symbols $dx=dm$, see \cite{CoSaVa92} for instance. For $i=1,..,d$ we define
left-invariant vector fields by%
\begin{equation}
U_{i}\left( x\right) =\partial _{i}+\frac{1}{2}\left( \sum_{1\leq j<i\leq
d}x^{1;j}\partial _{j,i}-\sum_{1\leq i<j\leq d}x^{1;j}\partial _{i,j}\right)
\label{DefOf_U_on_g2}
\end{equation}%
where $\partial _{i}$ denotes the coordinate vector field on $\mathbb{R}^{d}$
and $\partial _{i,j}$ with $i<j$ the coordinate vector field on $so\left(
d\right) $, identified with its upper diagonal elements. A simple
computation shows that $\left[ U_{i},U_{j}\right] =\partial _{i,j}$ and all
higher brackets are zero. Since H\"{o}rmander's condition is satisfied, we
call $\nabla =\left( U_{1},\cdots ,U_{d}\right) $ the \textit{hypoelliptic
gradient}. A (symmetric, regular, strongly local) Dirichlet form on $%
L^{2}\left( g^{2}\left( \mathbb{R}^{d}\right) ,dx\right) $ is defined by 
\begin{equation*}
\mathcal{E}\left( f,g\right) =\int_{g^{2}\left( \mathbb{R}^{d}\right)
}\,\nabla f\cdot \nabla g\,dm
\end{equation*}%
with domain $\mathcal{F}:=D\left( \mathcal{E}\right) :=\left\{ f\in L^{2}:%
\mathcal{E}\left( f,f\right) <\infty \right\} ,$ closure of smooth compactly
support functions with respect to%
\begin{equation*}
\left\Vert f\right\Vert _{\mathcal{F}}=\left( \mathcal{E}\left( f,f\right)
+\left\langle f,f\right\rangle _{L^{2}\left( g^{2}\left( \mathbb{R}%
^{d}\right) \right) }\right) ^{1/2}.
\end{equation*}%
This is a very standard setting, see \cite{Fu94} and \cite{CoSaVa92}, and as
pointed out in the introduction, $\mathcal{E}$ is the Dirichlet form
associated to the Markov process Brownian Motion plus its canonical Levy
area. The Dirichlet form $\mathcal{E}$ is based on the \textit{carr\'{e} du
champ operator}%
\begin{equation*}
\Gamma \left( f,g\right) =\nabla f\cdot \nabla g=\sum_{i=1}^{d}U_{i}f\left(
.\right) U_{i}g\left( .\right) ,
\end{equation*}%
which can be defined for all $f,g\in \mathcal{F}_{loc}=\left\{ f\in
L^{2}:\Gamma \left( f,f\right) \in L_{loc}^{1}\left( dm\right) \right\} $.
The associated \textit{energy measure} is simply $d\Gamma \left( f,g\right)
:=\Gamma \left( f,g\right) dm.$ Given $x,y\in g^{2}\left( \mathbb{R}%
^{d}\right) $ the (left-invariant) \textit{Carnot-Caratheodory }or \textit{%
control} distance $d\left( x,y\right) $ is defined as the length of the
shortest path from \thinspace $x$ to $y$ which remains tangent to $%
span\left\{ U_{1},...,U_{d}\right\} $, and the induced topology coincides
with the original topology of $g^{2}\left( \mathbb{R}^{d}\right) $; the 
\textit{Carnot-Caratheodory norm} is defined as $\left\Vert x\right\Vert
=d\left( 0,x\right) $. See \cite{CoSaVa92}, \cite{montgomery-02} or \cite%
{friz-victoir-04-Note}. From \cite[Lemma 5.29]{CKS87}, this distance
coincides with the \textit{intrinsic metric} of $\mathcal{E}$,\textit{\ } 
\begin{equation*}
d\left( x,y\right) =\sup \left\{ f\left( x\right) -f\left( y\right) :f\in 
\mathcal{F}_{loc}\text{ and }f\text{ continuous, }\Gamma \left( f,f\right)
\leq 1\right\} .
\end{equation*}

\begin{proposition}
\label{Sturm} \textrm{(I)} Completeness Property: In the metric space $%
\left( g^{2}\left( \mathbb{R}^{d}\right) ,d\right) $, every closed ball $%
\bar{B}$%
\begin{equation*}
\bar{B}\left( x,r\right) =\left\{ y\in g^{2}\left( \mathbb{R}^{d}\right)
:d\left( x,y\right) \leq r\right\}
\end{equation*}%
is complete and compact.\newline
\textrm{(II)} Doubling Property: The volume-doubling property 
\begin{equation*}
\forall r\geq 0:m\left( B\left( x,2r\right) \right) \leq 2^{N}m\left(
B\left( x,r\right) \right) .
\end{equation*}%
holds with $N=d^{2}.$\newline
\textrm{(III)} Poincar\'{e} Inequality: For all $r\geq 0$ and $f\in D\left( 
\mathcal{E}\right) $ 
\begin{equation*}
\int_{B\left( x,r\right) }\left\vert f-\bar{f}_{r}\right\vert ^{2}dm\leq C_{%
\ref{Sturm}}r^{2}\int_{B\left( x,r\right) }\Gamma \left( f,f\right) \,dm
\end{equation*}%
where%
\begin{equation*}
\bar{f}_{r}=m\left( B\left( x,r\right) \right) ^{-1}\int_{B\left( x,r\right)
}fdm.
\end{equation*}%
\textrm{(IV) }Nash Inequality: For all $f\in D\left( \mathcal{E}\right) \cap
L^{1}$ we have%
\begin{equation*}
\left\Vert f\right\Vert _{L^{2}}^{2+4/d^{2}}\leq C_{\ref{Sturm}}^{\prime }%
\mathcal{E}\left( f,f\right) \left\Vert f\right\Vert _{L^{1}}^{4/d^{2}}.
\end{equation*}
\end{proposition}

\begin{proof}
Property \textrm{(I)} is a simple consequence of $\left( g^{2}\left( \mathbb{%
R}^{d}\right) ,d\right) $ being complete, property \textrm{(II)} follows
from left then, every closed subset is complete. \textrm{(II)} follows
readily from invariance of $m$ under translation, $B\left( 0,r\right)
=\delta _{r}B\left( 0,1\right) $ and the Jacobian of $\delta _{\lambda }$
(as map from $g^{2}\left( \mathbb{R}^{d}\right) =\mathbb{R}^{d}\oplus
so\left( d\right) $ into itself) being equal to $\lambda ^{d}.\left( \lambda
^{2}\right) ^{\frac{d\left( d-1\right) }{2}}=\lambda ^{d^{2}}$. Property 
\textrm{(III)} appears explicitly in an appropriate Lie group setting in 
\cite{Je86}. At last, Property \textrm{(IV)} follows from \cite{CKS87},\cite%
{Rob91} or \cite{CoSaVa92}.
\end{proof}

\section{Uniformly Subelliptic Dirichlet Forms}

For $\Lambda \geq 1$ we call $\Xi \left( \Lambda \right) $ the set of all
measurable maps $a$ from $g^{2}\left( \mathbb{R}^{d}\right) $ into the space
of symmetric matrics such that 
\begin{equation*}
\forall \xi \in \mathbb{R}^{d}:\frac{1}{\Lambda }\left\vert \xi \right\vert
^{2}\leq \mathcal{\xi }\cdot a\xi \leq \Lambda \left\vert \xi \right\vert
^{2}.
\end{equation*}%
A symmetric Dirichlet form on $L^{2}\left( g^{2}\left( \mathbb{R}^{d}\right)
,dx\right) $ is defined by%
\begin{eqnarray*}
\mathcal{E}^{a}\left( f,g\right) &=&\int_{g^{2}\left( \mathbb{R}^{d}\right)
}\,\nabla f\left( x\right) \cdot a\left( x\right) \nabla g\left( x\right)
\,dm \\
&=&\sum_{i,j=1}^{d}\int_{g^{2}\left( \mathbb{R}^{d}\right) }\,a^{ij}\left(
x\right) U_{i}f\left( x\right) U_{j}g\left( x\right) \,dx.
\end{eqnarray*}%
The associated carr\'{e} du champ operator and energy measure are given by 
\begin{equation*}
\Gamma ^{a}\left( f,g\right) =\nabla f\left( x\right) \cdot a\left( x\right)
\nabla g\left( x\right) ,\,\text{\ }d\Gamma ^{a}\left( f,g\right) =\Gamma
^{a}\left( f,g\right) dm,
\end{equation*}%
respectively. The forms $\mathcal{E}^{a}$ and $\mathcal{E}$ are \textit{%
quasi-isometric} in the sense that $D\left( \mathcal{E}\right) =D\left( 
\mathcal{E}^{a}\right) $ and for all $f$ in the common domain,%
\begin{equation}
\frac{1}{\Lambda }\mathcal{E}\left( f,f\right) \leq \mathcal{E}^{a}\left(
f,f\right) \leq \Lambda \mathcal{E}\left( f,f\right) .  \label{quasiIsometry}
\end{equation}

The intrinsic metric associated to $\mathcal{E}^{a}\left( f,f\right) ,$ 
\begin{equation*}
d^{a}\left( x,y\right) =\sup \left\{ f\left( x\right) -f\left( y\right)
:f\in \mathcal{F}_{loc}\text{ and }f\text{ continuous, }\Gamma ^{a}\left(
f,f\right) \leq 1\right\} ,
\end{equation*}%
is obviously Lipschitz equivalent to $d\left( x,y\right) $ and hence a
metric on $g^{2}\left( \left( \mathbb{R}^{d}\right) \right) $ which induces
the original topology so that, in particular, $d^{a}\left( \cdot ,\cdot
\right) $ is continuous. Moreover, $\left( g^{2}\left( \left( \mathbb{R}%
^{d}\right) \right) ,d^{a}\right) $ is complete since $\left( g^{2}\left(
\left( \mathbb{R}^{d}\right) \right) ,d\right) $ is and closed balls are
easily seen to be compact, see property \textrm{(I)} above and in
Propositions \ref{Sturm} and \ref{I,II,IIIda}. The following proposition is
a special case of a result in \cite{Sturm93}.

\begin{proposition}
\label{da_is_geodesic}For all $a\in \Xi \left( \Lambda \right) $, the space $%
\left( g^{2}\left( \left( \mathbb{R}^{d}\right) \right) ,d^{a}\right) $ is a
geodesic space in the sense that for all $x,y$ there exists a continuous map 
$\gamma :\left[ 0,1\right] \rightarrow $ $g^{2}\left( \left( \mathbb{R}%
^{d}\right) \right) $ with $\gamma _{0}=x,\,\gamma _{1}=y$ and%
\begin{equation*}
d^{a}\left( \gamma _{r},\gamma _{t}\right) =d^{a}\left( \gamma _{r},\gamma
_{s}\right) +d^{a}\left( \gamma _{s},\gamma _{t}\right) \text{ \ for all }%
0\leq r<s<t\leq 1.
\end{equation*}
\end{proposition}

\begin{proposition}
\label{I,II,IIIda}Let $a\in \Xi \left( \Lambda \right) $. Properties \textrm{%
(I),(II),(III),(IV) }in proposition \ref{Sturm} remain valid when we replace 
$\mathcal{E}$ by $\mathcal{E}^{a}$ and $d$ by $d^{a}.$
\end{proposition}

\begin{proof}
Such properties are invariant under quasi-isometry, i.e. whenever we have (%
\ref{quasiIsometry}). This is easy to see for properties \textrm{(I), (II),
(IV)}. Invariance of the Poincar\'{e} inequality \textrm{(III)}, discussed
in detail in \cite{SturmIII}, is seen by first proving that the Poincar\'{e}
inequality is equivalent to a \textit{weak Poincar\'{e} inequality} for
which quasi-isometry is obvious.)
\end{proof}

Standard semigroup theory \cite{Fu94, Dav89} allows us to associate a
non-positive self-ajoint operator $L^{a}$ to $\mathcal{E}^{a}$. We then have%
\footnote{%
In view of De Giorgi-Moser-Nash regularity, see below, we may indeed write 
\textrm{inf, sup} rather than \textrm{ess-inf, ess-sup}$.$}

\begin{proposition}
\label{HarnackP}\textrm{(V) }Parabolic Harnack Inequality: Let $a\in \Xi
\left( \Lambda \right) $. There exists a constant $C_{\ref{HarnackP}}=C_{\ref%
{HarnackP}}\left( \Lambda \right) $ such that%
\begin{equation*}
\sup_{\left( s,y\right) \in Q^{-}}u\left( s,y\right) \leq C_{\ref{HarnackP}%
}\inf_{\left( s,y\right) \in Q^{+}}u\left( s,y\right) ,
\end{equation*}%
whenever $u$ is a nonnegative weak solution of the parabolic partial
differential equation $\partial _{t}u=L^{a}u$ on some cylinder $Q=\left(
t-4r^{2},t\right) \times B\left( x,2r\right) $ for some reals $t,r>0$. Here, 
$Q^{-}=\left( t-3r^{2},t-2r^{2}\right) \times B\left( x,r\right) $ and $%
Q^{+}=\left( t-r^{2},t\right) \times B\left( x,r\right) $ are lower and
upper sub-cylinders of $Q$ separated by a lapse of time. The statement
remains valid for balls with respect to $d^{a}$.
\end{proposition}

\begin{proof}
Based on the classical ideas by Moser \cite{Mo64, Mo71}, Grigor'yan,
Saloff-Coste, it is shown in \cite{SturmIII} that if $\mathrm{(I)}$ holds
then $\mathrm{(II)+}$ $\mathrm{(III)\Leftrightarrow (V)}$. For a more direct
proof along ideas of Nash, see \cite{St88, SaSt91}.
\end{proof}

Following \cite{FaSt86, St88, SturmIII} (these paper building on the seminal
works of De Giorgi-Moser-Nash) we have also H\"{o}lder regularity of such
weak solution (and in particular of the heat kernels discussed below). We
will refer to this simply as \textit{De Giorgi-Moser-Nash regularity}:

\begin{proposition}
\label{MoserNash}Let $a\in \Xi \left( \Lambda \right) $. Then there exist
constants $\eta \in \left( 0,1\right) $ and $C_{\ref{MoserNash}}$, only
depending on $\Lambda $, such that 
\begin{equation*}
\sup_{\left( s,y\right) ,\left( s^{\prime },y^{\prime }\right) \in
Q_{1}}\left\vert u\left( s,y\right) -u\left( s^{\prime },y^{\prime }\right)
\right\vert \leq C_{\ref{MoserNash}}\sup_{u\in Q_{2}}\left\vert u\right\vert
.\left( \frac{\left\vert s-s^{\prime }\right\vert ^{1/2}+d\left( y,y^{\prime
}\right) }{r}\right) ^{\eta }.
\end{equation*}%
whenever $u$ is a nonnegative weak solution of the parabolic partial
differential equation $\partial _{s}u=L^{a}u$ on some cylinder $Q_{2}\equiv
\left( t-4r^{2},t\right) \times B\left( x,2r\right) $ for some reals $t,r>0$%
. Here $Q_{1}\equiv \left( t-r^{2},t-2r^{2}\right) \times B\left( x,r\right) 
$ is a subcylinder of $Q_{2}$.
\end{proposition}

\subsection{Upper and Lower Heat Kernel Bounds}

Heat kernel existence is not an issue here. (For instance, \cite{FaSt86,
St88, Dav89, SaSt91}, Nash's inequality $\mathrm{(IV)}$ implies an estimate
on $\left\Vert P_{t}^{a}\right\Vert _{L^{1}\rightarrow L^{2}}$ and then via
duality on $\left\Vert P_{t}^{a}\right\Vert _{L^{1}\rightarrow L^{\infty }}$
which implies existence of the heat kernel $p^{a}=p^{a}\left( t,x,y\right) $%
.) We now turn to Aronson-type \cite{Ar67} heat-kernel estimates. As a
well-known consequence of our proposition \ref{I,II,IIIda} (see \cite[%
Corollary 4.2]{SturmIII}, also \cite{FaSt86, St88, Dav89}) we get

\begin{theorem}
\label{upper}Let $a\in \Xi \left( \Lambda \right) $. The heat kernel $p^{a}$
satisfies, for $\varepsilon >0$ fixed,%
\begin{equation*}
p^{a}\left( t,x,y\right) \leq \frac{C_{\ref{upper}}}{t^{d^{2}/2}}\exp \left(
-\frac{d^{a}\left( x,y\right) ^{2}}{\left( 4+\varepsilon \right) t}\right)
\end{equation*}%
for some constant $C_{\ref{upper}}=C_{\ref{upper}}\left( \varepsilon
,\Lambda \right) .$
\end{theorem}

\begin{theorem}
\label{lower}Let $a\in \Xi \left( \Lambda \right) $. The heat kernel $p^{a}$
satisfies%
\begin{equation*}
p^{a}\left( t,x,y\right) \geq \frac{1}{C_{\ref{lower}}}\frac{1}{t^{d^{2}/2}}%
\exp \left( -\frac{C_{\ref{lower}}d^{a}\left( x,y\right) ^{2}}{t}\right)
\end{equation*}%
for some constant $C_{\ref{lower}}=C_{\ref{lower}}\left( \Lambda \right) $.
\end{theorem}

Let $a\in \Xi \left( \Lambda \right) $. Let $C_{\ref{upper}},$ $C_{\ref%
{lower}}$ denote the constants of the previous two theorems. Then 
\begin{equation*}
\frac{1}{C_{\ref{lower}}}\frac{1}{t^{d^{2}/2}}\exp \left( -\frac{C_{\ref%
{lower}}\Lambda d\left( x,y\right) ^{2}}{t}\right) \leq p^{a}\left(
t,x,y\right) \leq \frac{C_{\ref{upper}}}{t^{d^{2}/2}}\exp \left( -\frac{%
d\left( x,y\right) ^{2}}{\Lambda \left( 4+\varepsilon \right) t}\right) .
\end{equation*}

\begin{proof}
Lipschitz-equivalence of $d\left( x,y\right) $ and $d^{a}\left( x,y\right) $.
\end{proof}

\subsection{The Associated Markov Process}

Following a standard construction, the heat kernel $p^{a}$ gives rise to a
consistent family of finite-dimensional distributions and determines a $%
g^{2}\left( \mathbb{R}^{d}\right) $-valued (strong) Markov process $\left( 
\mathbf{X}_{t}^{a,x}:t\geq 0\right) $ where $a\in \Xi \left( \Lambda \right) 
$ and $\mathbf{X}_{0}^{a,x}=x\in g^{2}\left( \mathbb{R}^{d}\right) $. The
natural time horizon is $[0,\infty )$ but our focus will be on finite time
horizon and by scaling (cf. next section) there is no loss of generality to
work on $\left[ 0,1\right] $. The heat kernel estimates are more than
enough, via Kolmogorov's criterion, to guarantee that any such process can
be taken with continuous sample paths; the law of $\mathbf{X}^{a,x}$ is then
denoted by $\mathbb{P}^{a,x}$, a Borel measure on $C\left( \left[ 0,1\right]
,g^{2}\left( \mathbb{R}^{d}\right) \right) $, under which we can think of $%
\mathbf{X=X}^{a,x}$ simply as coordinate process $\mathbf{X}_{t}\left(
\omega \right) =\omega _{t}$. By construction, the density of $\mathbf{X}%
_{t} $ under $\mathbb{P}^{a,x}$, or equivalently, the density of $\mathbf{X}%
_{t}^{a,x}$, with respect to $m$ is given by $p^{a}\left( t,x,\cdot \right)
. $

\subsection{Scaling\label{scalingDivForm}}

We will refer to the following simple proposition as \textit{scaling}.\
Recall that the dilation operator $\delta $ extends scalar multiplication to 
$g^{2}\left( \mathbb{R}^{d}\right) $.

\begin{proposition}
For any $a\in \Xi \left( \Lambda \right) ,r\neq 0$ set $a^{r}\left( x\right)
:=a\left( \delta _{1/r}\,x\right) \in \Xi \left( \Lambda \right) $. Then%
\begin{equation*}
\left( \mathbf{X}_{t}^{a^{r},x}:t\geq 0\right) \overset{\mathcal{D}}{=}%
\left( \delta _{r}\mathbf{X}_{t/r^{2}}^{a,\delta _{1/r}\left( x\right)
}:t\geq 0\right) .
\end{equation*}
\end{proposition}

\subsection{Short Time Asymptotics}

When $a=I$, the identity matrix, an essentially sharp lower bound with $1/C_{%
\ref{lower}}=4\left( 1-\varepsilon \right) $ is known, see \cite{Varo90}.
This implies Varadhan's formula%
\begin{equation*}
4t\log p^{I}\left( t,x,y\right) \rightarrow -d^{I}\left( x,y\right) ^{2}%
\text{ as }t\rightarrow 0\text{.}
\end{equation*}%
The generalization to arbitrary $a\in \Xi \left( \Lambda \right) $ follows
from the recent work of Ram\'{\i}rez \cite{Ra01} and will be central to our
discussion of large deviations.

\begin{theorem}
The heat kernel associated to $L^{a}$ satisfies, for all $x,y\in g^{2}\left( 
\mathbb{R}^{d}\right) $ 
\begin{equation*}
4t\log p^{a}\left( t,x,y\right) \rightarrow -d^{a}\left( x,y\right) ^{2}%
\text{ as }t\rightarrow 0.
\end{equation*}
\end{theorem}

\subsection{A Lower Bound for the Killed Process}

\begin{theorem}
\label{lowerKilled}Let $a\in \Xi \left( \Lambda \right) $. For $x_{0}\in
g^{2}\left( \mathbb{R}^{d}\right) $ and $r>0,$ define 
\begin{eqnarray*}
\xi _{B\left( x_{0},r\right) }^{a;x} &=&\inf \left\{ t\geq 0:\mathbf{X}%
_{t}^{a;x}\notin B\left( x_{0},r\right) \right\} , \\
\mathbb{P}_{B\left( x_{0},r\right) }^{a;x}\left( t,\cdot \right) &=&\mathbb{P%
}\left( \mathbf{X}_{t}^{a,x}\in \cdot \,\,,\xi _{B\left( x_{0},r\right)
}^{a;x}>t\right) .
\end{eqnarray*}%
Then $\mathbb{P}_{B\left( x_{0},r\right) }^{a;x}\left( t,dy\right)
=p_{B\left( x_{0},r\right) }^{a}\left( t,x,y\right) dy.$ \TEXTsymbol{>}%
\TEXTsymbol{>}\TEXTsymbol{>}\ CHECK $B$ vs $B^{a}$.\TEXTsymbol{<}\TEXTsymbol{%
<}\TEXTsymbol{<}Moreover, if $x,y$ are two elements of $B^{a}\left(
x_{0},r\right) $ joined by a curve $\gamma $ which is at a $d^{a}$-distance $%
R>0$ of $g^{2}\left( \mathbb{R}^{d}\right) /B^{a}\left( x_{0},r\right) $
there exists constant $C_{\ref{lowerKilled}}$ depending only on $\Lambda ,$%
\begin{equation*}
p_{B\left( x_{0},r\right) }^{a}\left( t,x,y\right) \geq \frac{1}{C_{\ref%
{lowerKilled}}\delta ^{d^{2}/2}}\exp \left( -C_{\ref{lowerKilled}}\frac{%
d^{a}\left( x,y\right) ^{2}}{t}\right) \exp \left( -\frac{C_{\ref%
{lowerKilled}}t}{R^{2}}\right)
\end{equation*}%
where $\delta =\min \left\{ t,R^{2}\right\} $.
\end{theorem}

\begin{proof}
See \cite{SturmIII} or \cite{SaSt91}, the ideas are adapted from \cite%
{FaSt86, St88}.
\end{proof}

One should observe that $d^{a}$ can be replaced by $d,$ at the price of
changing the constants.

\section{Construction of Associated Rough Paths}

In conjunction with the ever useful Garsia-Rodemich-Rumsey's lemma, the
upper heat bounds leads to H\"{o}lder regularity of the sample paths $%
t\mapsto \mathbf{X}_{t}^{a;x}\left( \omega \right) $. Moreover, a \textit{%
Fernique estimate} holds by which we mean that the homogenous H\"{o}lder
norm of the $g^{2}\left( \mathbb{R}^{d}\right) $-valued process $\mathbf{X}%
^{a;x}$ has a Gauss tail.

\begin{lemma}
\label{GaussianIntegrabilityMarkov}For all $\eta <\frac{1}{4\Lambda }$ we
have \ 
\begin{equation*}
\sup_{a\in \Xi \left( \Lambda \right) }\sup_{x\in g^{2}\left( \mathbb{R}%
^{d}\right) }\sup_{0\leq s<t\leq 1}\mathbb{E}^{a,x}\left( \exp \left( \eta 
\frac{d\left( \mathbf{X}_{t},\mathbf{X}_{s}\right) ^{2}}{t-s}\right) \right)
<\infty .
\end{equation*}
\end{lemma}

\begin{proof}
By scaling and the Markov property, for any $a\in \Xi \left( \Lambda \right)
,$%
\begin{equation*}
\sup_{x\in g^{2}\left( \mathbb{R}^{d}\right) }\sup_{0\leq s<t\leq 1}\mathbb{E%
}^{a,x}\left( \exp \left( \eta \frac{d\left( \mathbf{X}_{t},\mathbf{X}%
_{s}\right) ^{2}}{t-s}\right) \right) \leq \sup_{x\in g^{2}\left( \mathbb{R}%
^{d}\right) }\sup_{a\in \Xi \left( \Lambda \right) }\mathbb{E}^{a,x}\left(
\exp \left( \eta \left\Vert \mathbf{X}_{0,1}\right\Vert ^{2}\right) \right) .
\end{equation*}%
(Recall that $d\left( \mathbf{X}_{t},\mathbf{X}_{s}\right) =d\left( 0,%
\mathbf{X}_{s}^{-1}\ast \mathbf{X}_{t}\right) =\left\Vert \mathbf{X}%
_{s,t}\right\Vert $where $\left\Vert \cdot \right\Vert =d\left( 0,\cdot
\right) $ denotes the Carnot-Caratheodory norm.) Fix $\eta <\frac{1}{%
4\Lambda }$, and $\varepsilon >0$ such that $\eta <\frac{1}{4(1+\varepsilon
)\Lambda }$. Then, from the heat kernel upper-bound, we obtain%
\begin{eqnarray*}
\mathbb{E}^{a,x}\left( \exp \left( \eta \left\Vert \mathbf{X}%
_{0,1}\right\Vert ^{2}\right) \right) &=&\int \exp \left( \eta d\left(
x,y\right) ^{2}\right) p^{a}\left( 1,x,y\right) dy \\
&\leq &C_{\ref{upper}}\int \exp \left( -\left( \frac{1}{4(1+\varepsilon
)\Lambda }-\eta \right) d\left( x,y\right) ^{2}\right) dy
\end{eqnarray*}%
From $m\left( B\left( x,r\right) \right) =m\left( B\left( 0,1\right) \right)
r^{d^{2}}$ we have $dm\left( B\left( x,r\right) \right) /dr=m\left( B\left(
0,1\right) \right) d^{2}r^{d^{2}-1}$ so that%
\begin{equation*}
\mathbb{E}^{a,x}\left( \exp \left( \eta \left\Vert \mathbf{X}%
_{0,1}\right\Vert ^{2}\right) \right) \leq C_{\ref{upper}}m\left( B\left(
0,1\right) \right) d^{2}\int_{r=0}^{\infty }\exp \left( -\left( \frac{1}{%
4(1+\varepsilon )\Lambda }-\eta \right) r^{2}\right) r^{d^{2}-1}dr
\end{equation*}%
and by our choice of $\eta ,\varepsilon $ the right hand side is finite,
uniformly in $x$ and $a\in \Xi \left( \Lambda \right) $ as required.
\end{proof}

The previous lemma combined with a standard application of the
Garsia-Rodemich-Rumsey lemma leads immediately to Fernique estimate for
homogenous $\alpha $-H\"{o}lder norm 
\begin{equation*}
\left\Vert \mathbf{X}\right\Vert _{\alpha -H\ddot{o}l;\left[ 0,1\right]
}=\sup_{0\leq s<t\leq 1}\frac{d\left( \mathbf{X}_{t},\mathbf{X}_{s}\right) }{%
\left\vert t-s\right\vert ^{\alpha }}.
\end{equation*}%
More precisely, we have

\begin{theorem}
\label{FerniqueEstimates}Let $0\leq \alpha <1/2.$There exists a constant $C_{%
\ref{FerniqueEstimates}}=C_{\ref{FerniqueEstimates}}\left( \Lambda ,\alpha
\right) >0$ such that%
\begin{equation*}
\sup_{a\in \Xi \left( \Lambda \right) }\sup_{x\in g^{2}\left( \mathbb{R}%
^{d}\right) }\mathbb{E}^{a,x}\left[ \exp \left( C_{\ref{FerniqueEstimates}%
}\left\Vert \mathbf{X}\right\Vert _{\alpha \text{-H\"{o}l;}\left[ 0,1\right]
}^{2}\right) \right] <\infty .
\end{equation*}%
In particular, for $\alpha \in \left( 1/3,1/2\right) $ almost every sample
path $t\mapsto \mathbf{X}_{t}^{a;x}\left( \omega \right) $ is an $\alpha $-H%
\"{o}lder geometric rough path.\newline
\end{theorem}

For later use - namely our discussion of Wong-Zakai approximations - we
record the following estimate.

\begin{corollary}
\label{CorForUniBoundsOnPWLinearApprox}Let 
\begin{equation}
M_{\eta }:=\sup_{a\in \Xi \left( \Lambda \right) }\sup_{x\in g^{2}\left( 
\mathbb{R}^{d}\right) }\sup_{0\leq s<t\leq 1}\mathbb{E}^{a,x}\left( \exp
\left( \eta \frac{d\left( \mathbf{X}_{t},\mathbf{X}_{s}\right) ^{2}}{t-s}%
\right) \right) .  \label{DefMeta}
\end{equation}%
Then there exists $C_{\ref{CorForUniBoundsOnPWLinearApprox}}=C_{\ref%
{CorForUniBoundsOnPWLinearApprox}}\left( \Lambda \right) $ such that $%
M_{\eta }\leq \exp \left( C_{\ref{CorForUniBoundsOnPWLinearApprox}}\eta
\right) $ for all $\eta \in \left[ 0,\frac{1}{16\Lambda }\right) $.
\end{corollary}

\begin{proof}
It suffices to show $M_{\eta }\leq 1+C_{\ref{CorForUniBoundsOnPWLinearApprox}%
}\eta $. From the inequality $\exp \left( x\right) \leq 1+x\exp \left(
x\right) $ for $x>0$ we obtain%
\begin{equation*}
M_{\eta }\leq 1+\eta \sup_{x\in \mathbb{R}^{d}}\sup_{s<t\in \left[ 0,1\right]
}\mathbb{E}^{a,x}\left( \frac{d\left( \mathbf{X}_{t},\mathbf{X}_{s}\right)
^{2}}{t-s}\exp \left( \eta \frac{d\left( \mathbf{X}_{t},\mathbf{X}%
_{s}\right) ^{2}}{t-s}\right) \right) .
\end{equation*}%
Define%
\begin{equation*}
Q_{4}:=\sup_{a\in \Xi \left( \Lambda \right) }\sup_{x\in g^{2}\left( \mathbb{%
R}^{d}\right) }\sup_{s<t\in \left[ 0,1\right] }\mathbb{E}^{a,x}\left( \frac{%
d\left( \mathbf{X}_{t},\mathbf{X}_{s}\right) ^{4}}{\left\vert t-s\right\vert
^{2}}\right) <\infty .
\end{equation*}%
The proof is now finished by Cauchy-Schwarz,%
\begin{eqnarray*}
M_{\eta } &\leq &1+\eta Q_{4}^{1/2}\sqrt{\sup_{x\in \mathbb{R}%
^{d}}\sup_{s<t\in \left[ 0,1\right] }\mathbb{E}^{a,x}\left( \exp \left(
2\eta \frac{d\left( \mathbf{X}_{t},\mathbf{X}_{s}\right) ^{2}}{t-s}\right)
\right) } \\
&\leq &1+\eta Q_{4}^{1/2}\sqrt{\sup_{x\in \mathbb{R}^{d}}\sup_{s<t\in \left[
0,1\right] }\mathbb{E}^{a,x}\left( \exp \left( \frac{1}{8\Lambda }\frac{%
d\left( \mathbf{X}_{t},\mathbf{X}_{s}\right) ^{2}}{t-s}\right) \right) }
\end{eqnarray*}%
and Lemma \ref{GaussianIntegrabilityMarkov}.
\end{proof}

\section{Approximations}

\subsection{Weak Convergence}

\begin{theorem}
Let $\left( a_{n}\right) $ be a sequence of (smooth) functions in $\Xi
\left( \Lambda \right) $ such that $a_{n}$ converges almost everywhere to $%
a\in \Xi \left( \Lambda \right) $. Then we have\newline
(i)\ uniformly on compacts in $\left( 0,\infty \right) \times g^{2}\left( 
\mathbb{R}^{d}\right) \times g^{2}\left( \mathbb{R}^{d}\right) $, 
\begin{equation*}
p^{a_{n}}\left( t,x,y\right) \rightarrow p^{a}\left( t,x,y\right) \text{ as }%
n\rightarrow \infty ;
\end{equation*}%
(ii) convergence in distribution $\mathbf{X}^{a_{n},x}\overset{\mathcal{D}}{%
\rightarrow }\mathbf{X}^{a,x}$ with respect to uniform topology on $\left\{
\omega :C\left( \left[ 0,1\right] ,g^{2}\left( \mathbb{R}^{d}\right) \right)
:\omega \left( 0\right) =x\right\} $, with fixed $x\in g^{2}\left( \mathbb{R}%
^{d}\right) ;$\newline
(iii) the convergence in distribution remains valid with respect to
homogenous $\alpha $-H\"{o}lder topology of exponent for $\alpha \in \lbrack
0,1/2)$.
\end{theorem}

\begin{proof}
The proof of (i) is identical to the proof of \cite[Theorem II.3.1]{St88}
and implies convergence of the finite-dimensional distributions. A standard
tightness argument leads to (ii) and (iii).
\end{proof}

\begin{remark}
\cite{LejII06} discusses the case when $a\left( x\right) $ depends only on
the projection $\pi _{1}\left( x\right) \in \mathbb{R}^{d}$.
\end{remark}

\subsection{Strong Convergence}

\subsubsection{Geodesic Approximations}

Recall that $g^{2}\left( \mathbb{R}^{d}\right) $ equipped with
Carnot-Caratheodory distance is a geodesic space. Given a dissection $D$ of $%
\left[ 0,1\right] $ and a deterministic path $\mathbf{x}\in C^{\alpha \text{%
-H\"{o}lder}}\left( \left[ 0,1\right] ,g^{2}\left( \mathbb{R}^{d}\right)
\right) $ we can approximate $\mathbf{x}$ by a path $\mathbf{x}^{D}\in C^{%
\text{Lip}}\left( \left[ 0,1\right] ,g^{2}\left( \mathbb{R}^{d}\right)
\right) $ obtained by connecting the points $\left( \mathbf{x}%
_{t_{i}}:t_{i}\in D\right) $ with geodesics run at unit speed. If there are
several geodesics between two points $\mathbf{x}_{t_{i}}$ and $\mathbf{x}%
_{t_{i+1}}$ it is immaterial which one is chosen. It is not hard to show
that 
\begin{equation}
\left\Vert \mathbf{x}^{D}\right\Vert _{\alpha \text{-H\"{o}lder}}\leq
3\left\Vert \mathbf{x}\right\Vert _{\alpha \text{-H\"{o}lder}}\text{.}
\label{UniformityInGeoApprox}
\end{equation}%
Clearly, $\mathbf{x}^{D}\rightarrow \mathbf{x}$ pointwise as $\left\vert
D\right\vert \rightarrow 0$ and, in fact, this convergence is uniform in
view of the uniform bound (\ref{UniformityInGeoApprox}). A simple
interpolation argument then gives $\alpha ^{\prime }$-H\"{o}lder
convergence, $\alpha ^{\prime }\in \left( 0,\alpha \right) $. All this
results are purely deterministic and discussed in detail in \cite%
{friz-victoir-04-Note}. By Theorem \textbf{\ref{FerniqueEstimates}} these
approximation results apply to a.e. sample path of $\mathbf{X}^{a,x}$. We
emphasize that these approximations required apriori knowledge of the area $%
\pi _{2}\left( \mathbf{X}^{a,x}\right) $. In fact, $\pi _{1}\left( \mathbf{x}%
^{D}\right) $ is simply the concatenation of path segments designed to wipe
out prescribed areas.

\subsubsection{Piecewise Linear Approximations:\ Wong-Zakai\label%
{SectionWongZakai}}

In contrast to geodesic approximation, convergence of piecewise linear
approximations, based on the $\mathbb{R}^{d}$-valued path $\pi _{1}\left( 
\mathbf{X}^{a,x}\right) $ alone and without apriori knowledge of the area $%
\pi _{2}\left( \mathbf{X}^{a,x}\right) $, is a genuine probabilistic
statement and relies on subtle cancellations. (An example by McShane, see 
\cite{ikeda-watanabe-89}, shows what can go wrong if one replaces linear
cords by general interpolation functions.)

\paragraph{The Idea}

Fix a dissection $D=\left\{ t_{i}:i\right\} $ of $[0,1]$ and $a\in \Xi
\left( \Lambda \right) $. Let us project $\mathbf{X=X}^{a}$ to the $\mathbb{R%
}^{d}$-valued process $X=X^{a}$ and consider piecewise-linear approximations
to $X$ based on $D$, denoted by $X^{D}$. Of course, $X^{D}$ has a
canonically defined area given by the usual iterated integrals and thus
gives rise to an $g^{2}\left( \mathbb{R}^{d}\right) $-valued path which we
denote by $S\left( X^{D}\right) $. For $0\leq \alpha <1/2$ as usual, the
convergence 
\begin{equation}
d_{\alpha \text{-H\"{o}lder}}\left( S\left( X^{D}\right) ,\mathbf{X}\right)
\rightarrow 0\text{ in probability }  \label{WongZakaiDivForm}
\end{equation}%
as $\left\vert D\right\vert \rightarrow 0$ is a subtle problem and the
difficulty is already present in the pointwise convergence statement $%
S\left( X^{D}\right) _{0,t}\rightarrow \mathbf{X}_{0,t}$ as $\left\vert
D\right\vert \rightarrow 0.$ Our idea is simple. Noting that straight line
segments do not produce area, it is an elementary application of the
Campbell-Baker-Hausdorff formula to see that for $t\in D=\left\{
t_{i}\right\} $ 
\begin{equation}
\left( S\left( X^{D}\right) _{0,t}\right) ^{-1}\ast \mathbf{X}%
_{0,t}=\sum_{i}A_{t_{i},t_{i+1}},  \label{SumOfLittleAreas}
\end{equation}%
where $A$ is the area of $\mathbf{X}$ and $\cup _{i}\left[ t_{i},t_{i+1}%
\right] =[0,t]$. On the other hand, it is relatively straight-forward to
show that the $L^{p}$ norm of $\left\Vert S\left( X^{D}\right) \right\Vert
_{\alpha \text{-H\"{o}l;}\left[ 0,1\right] }$ is finite uniformly over all $%
D $. In essence, this reduces (\ref{WongZakaiDivForm}) to the pointwise
convergence statement which we can rephrase as $\sum_{i}A_{t_{i},t_{i+1}}%
\rightarrow 0.$ It is natural to show this in $L^{2}$ since this allows to
write\footnote{%
We equip $so\left( d\right) \subset \mathbb{R}^{d}\otimes \mathbb{R}^{d}$
with the Euclidean structure $A\cdot \tilde{A}=\sum_{k,l=1}^{d}A^{k,l}\tilde{%
A}^{k,l}$ and $\left\vert A\right\vert ^{2}=A\cdot A$. It may be instructive
to consider $d=2$ in which case $A$ can be viewed as scalar.}%
\begin{equation*}
\mathbb{E}\left[ \left\vert \sum_{i}A_{t_{i},t_{i+1}}\right\vert ^{2}\right]
=\sum_{i}\mathbb{E}\left( \left\vert A_{t_{i},t_{i+1}}\right\vert
^{2}\right) +2\sum_{i<j}\mathbb{E}\left( A_{t_{i},t_{i+1}}\cdot
A_{t_{j},t_{j+1}}\right) .
\end{equation*}%
For simplicity only, assume $t_{i+1}-t_{i}\equiv \delta $ for all $i$. As a
sanity check, if $X$ were a Brownian motion and $A$ the usual L\'{e}vy area,
all off-diagonal terms are zero and%
\begin{equation*}
\sum_{i}\mathbb{E}\left( \left\vert A_{t_{i},t_{i+1}}\right\vert ^{2}\right)
\sim \sum_{i}\delta ^{2}\sim \frac{1}{\delta }\delta ^{2}\rightarrow 0\text{
with }\left\vert D\right\vert =\delta \rightarrow 0
\end{equation*}%
which is what we want. Back to the general case of $\mathbf{X=X}^{a}$, the
plan must be to cope with the off-diagonal sum. Since there are $\sim \delta
^{2}/2$ terms what we need is $\mathbb{E}\left( A_{t_{i},t_{i+1}}\cdot
A_{t_{j},t_{j+1}}\right) =o\left( \delta ^{2}\right) .$To this end, let us
momentarily assume that%
\begin{equation}
\sup_{x}\mathbb{E}^{a,x}\left( A_{0,\delta }\right) =o\left( \delta \right) .
\label{TooStrongAssumptionWongZakai}
\end{equation}%
holds. Then, using the Markov property,%
\begin{equation*}
\left\vert \mathbb{E}\left( A_{t_{i},t_{i+1}}\cdot A_{t_{j},t_{j+1}}\right)
\right\vert \leq \mathbb{E}\left( \left\vert A_{t_{i},t_{i+1}}\right\vert
\times \left\vert \mathbb{E}^{\mathbf{X}_{t_{j}}}A_{0,\delta }\right\vert
\right) =\mathbb{E}\left( \left\vert A_{t_{i},t_{i+1}}\right\vert \right)
\times o\left( \delta \right)
\end{equation*}%
and since $\mathbb{E}\left( \left\vert A_{t_{i},t_{i+1}}\right\vert \right)
\sim \delta $, by a soft scaling argument, we are done. Unfortunately, (\ref%
{TooStrongAssumptionWongZakai}) seems to be too strong to be true but we are
able to establish a weak version of (\ref{TooStrongAssumptionWongZakai})
which is good enough to successfully implement what we just outlined. The
key to all this (cf. the proof of the forthcoming Proposition \ref%
{weakConvegenceArea0}) is a semi-group argument which leads to the desired
cancellations.

\paragraph{Uniform H\"{o}lder Bound}

Let $X^{D}$ denote the piecewise linear approximation to $X=X\left( \omega
\right) $. We now show $L^{q}\left( \mathbb{P}^{a,x}\right) $-bounds,
uniformly over all dissections $D$, of the homogenous $\alpha $-H\"{o}lder
norm of the path $X^{D}$ and its area.

\begin{theorem}
\bigskip \label{UniformGaussTailOfApproximations}There exists $\eta =\eta $ $%
\left( \Lambda \right) >0$ such that 
\begin{equation*}
\sup_{a\in \Xi \left( \Lambda \right) ,x\in g^{2}\left( \mathbb{R}%
^{d}\right) }\sup_{D}\sup_{0\leq s<t\leq 1}\mathbb{E}^{a,x}\left( \exp
\left( \eta \frac{\left\Vert S\left( X^{D}\right) _{s,t}\right\Vert ^{2}}{t-s%
}\right) \right) <\infty .
\end{equation*}%
As a consequence, for any $\alpha \in \lbrack 0,1/2)$ there exists $C_{\ref%
{UniformGaussTailOfApproximations}}=C_{\ref{UniformGaussTailOfApproximations}%
}\left( \alpha ,\Lambda \right) >0$ so that%
\begin{equation*}
\sup_{a\in \Xi \left( \Lambda \right) ,x\in g^{2}\left( \mathbb{R}%
^{d}\right) }\sup_{D}\mathbb{E}^{a,x}\left( \exp \left( C_{\ref%
{UniformGaussTailOfApproximations}}\left\Vert S\left( X^{D}\right)
\right\Vert _{\alpha \text{-H\"{o}l;}\left[ 0,1\right] }^{2}\right) \right)
<\infty .
\end{equation*}
\end{theorem}

\begin{proof}
The consequence is an immediate application of the Garsia-Rodemich-Rumsey
lemma and we only have to discuss the first estimate. We remind the reader
that from Lemma \ref{GaussianIntegrabilityMarkov} for $\eta \in \lbrack 0,%
\frac{1}{4\Lambda }),$ 
\begin{equation*}
M_{\eta }\equiv \sup_{a\in \Xi \left( \Lambda \right) ,x\in g^{2}\left( 
\mathbb{R}^{d}\right) }\sup_{0\leq s<t\leq 1}\mathbb{E}^{a,x}\left( \exp
\left( \eta \frac{\left\Vert \mathbf{X}_{s,t}\right\Vert ^{2}}{t-s}\right)
\right) <\infty .
\end{equation*}%
By the triangle inequality (recall $t_{D},t^{D}$ were defined at the end of
the introduction) 
\begin{eqnarray*}
\frac{\left\Vert S\left( X^{D}\right) _{s,t}\right\Vert }{\sqrt{t-s}} &\leq &%
\frac{\left\Vert S\left( X^{D}\right) _{s,s^{D}}\right\Vert }{\sqrt{s^{D}-s}}%
+\frac{\left\Vert S\left( X^{D}\right) _{s^{D},t_{D}}\right\Vert }{\sqrt{%
t_{D}-s^{D}}}+\frac{\left\Vert S\left( X^{D}\right) _{t_{D},t}\right\Vert }{%
\sqrt{t-t_{D}}} \\
&\leq &\frac{\left\vert X_{s,s^{D}}^{D}\right\vert }{\sqrt{s^{D}-s}}+\frac{%
\left\Vert S\left( X^{D}\right) _{s^{D},t_{D}}\right\Vert }{\sqrt{t_{D}-s^{D}%
}}+\frac{\left\vert X_{t_{D},t}^{D}\right\vert }{\sqrt{t-t_{D}}} \\
&\leq &\frac{\left\Vert \mathbf{X}_{s,s^{D}}\right\Vert }{\sqrt{s^{D}-s}}+%
\frac{\left\Vert S\left( X^{D}\right) _{s^{D},t_{D}}\right\Vert }{\sqrt{%
t_{D}-s^{D}}}+\frac{\left\Vert \mathbf{X}_{t_{D},t}\right\Vert }{\sqrt{%
t-t_{D}}} \\
&\leq &\left( \frac{3\left\Vert \mathbf{X}_{s,s^{D}}\right\Vert ^{2}}{s^{D}-s%
}+\frac{3\left\Vert S\left( X^{D}\right) _{s^{D},t_{D}}\right\Vert ^{2}}{%
t_{D}-s^{D}}+\frac{3\left\Vert \mathbf{X}_{t_{D},t}^{x}\right\Vert ^{2}}{%
t-t_{D}}\right) ^{1/2}.
\end{eqnarray*}%
Hence%
\begin{multline*}
\mathbb{E}^{a,x}\left( \exp \left( \eta \frac{\left\Vert S\left(
X^{D}\right) _{s,t}\right\Vert ^{2}}{t-s}\right) \right) \\
\left. 
\begin{array}{l}
\leq \mathbb{E}^{a,x}\left\{ \exp \left[ \eta \left( \frac{3\left\Vert 
\mathbf{X}_{s,s^{D}}\right\Vert ^{2}}{s^{D}-s}+\frac{3\left\Vert S\left(
X^{D}\right) _{s^{D},t_{D}}\right\Vert ^{2}}{t_{D}-s^{D}}+\frac{3\left\Vert 
\mathbf{X}_{t_{D},t}^{x}\right\Vert ^{2}}{t-t_{D}}\right) \right] \right\}
\\ 
\leq M_{6\eta }^{2}\mathbb{E}^{a,x}\left( \exp \left( 6\eta \frac{\left\Vert
S\left( X^{D}\right) _{s^{D},t_{D}}\right\Vert ^{2}}{t_{D}-s^{D}}\right)
\right)%
\end{array}%
\right.
\end{multline*}%
and the proof is reduced to show that for some $\eta >0$ small enough 
\begin{equation*}
\sup_{a\in \Xi \left( \Lambda \right) ,x\in g^{2}\left( \mathbb{R}%
^{d}\right) }\sup_{D}\sup_{s<t\in D}\mathbb{E}^{a,x}\left( \exp \left( 6\eta 
\frac{\left\Vert S\left( X^{D}\right) _{s,t}\right\Vert ^{2}}{t-s}\right)
\right) <\infty .
\end{equation*}%
By the triangle inequality for the Carnot-Caratheodory distance, for $%
t_{i},t_{j}\in D,$%
\begin{equation*}
\left\Vert S\left( X^{D}\right) _{t_{i},t_{j}}\right\Vert \leq \left\Vert 
\mathbf{X}_{t_{i},t_{j}}\right\Vert +d\left( \mathbf{X}_{t_{i},t_{j}},S%
\left( X^{D}\right) _{t_{i},t_{j}}\right) .
\end{equation*}%
To proceed we note that, similar to equation (\ref{SumOfLittleAreas}), 
\begin{equation*}
\left( S\left( X^{D}\right) _{t_{i},t_{j}}\right) ^{-1}\ast \mathbf{X}%
_{t_{i},t_{j}}=\sum_{k=i}^{j-1}A_{t_{k},t_{k+1}}.
\end{equation*}%
By left-invariance of the Carnot-Caratheodory distance $d$ and equivalence
of continuous homogenous norms (so that, in particular, $\left\Vert
(x,A)\right\Vert \sim \left\vert x\right\vert +\left\vert A\right\vert
^{1/2} $ where $\left\vert \cdot \right\vert $ denotes Euclidean norm on $%
\mathbb{R}^{d}$ resp. $\mathbb{R}^{d}\otimes \mathbb{R}^{d}$) there exists $%
C $ such that%
\begin{eqnarray*}
d\left( \mathbf{X}_{t_{i},t_{j}},S\left( X^{D}\right) _{t_{i},t_{j}}\right)
&=&\left\Vert \left( 0,\sum_{k=i}^{j-1}A_{t_{k},t_{k+1}}\right) \right\Vert
\\
&\leq &C\left\vert \sum_{k=i}^{j-1}A_{t_{k},t_{k+1}}\right\vert ^{1/2}\leq C%
\sqrt{\sum_{k=i}^{j-1}\left\vert A_{t_{k},t_{k+1}}\right\vert } \\
&\leq &C\sqrt{\sum_{k=i}^{j-1}\left\Vert \mathbf{X}_{t_{k},t_{k+1}}\right%
\Vert ^{2}}.
\end{eqnarray*}%
By Cauchy-Schwartz, 
\begin{eqnarray*}
&&\mathbb{E}^{a,x}\left( \exp \left( 6\eta \frac{\left\Vert S\left(
X^{D}\right) _{t_{i},t_{j}}\right\Vert ^{2}}{t_{j}-t_{i}}\right) \right) \\
&&\left. 
\begin{array}{l}
\leq \mathbb{E}^{a,x}\left( \exp \left( 12\eta \frac{\left\Vert \mathbf{X}%
_{t_{i},t_{j}}\right\Vert ^{2}}{t_{j}-t_{i}}\right) \exp \left( 12C\eta 
\frac{\sum_{k=i}^{j-1}\left\Vert \mathbf{X}_{t_{k},t_{k+1}}\right\Vert ^{2}}{%
t_{j}-t_{i}}\right) \right) \\ 
\leq M_{24\eta }\mathbb{E}^{a,x}\left( \prod_{k=i}^{j-1}\exp \left( 24C\eta 
\frac{\left\Vert \mathbf{X}_{t_{k},t_{k+1}}\right\Vert ^{2}}{t_{j}-t_{i}}%
\right) \right) .%
\end{array}%
\right.
\end{eqnarray*}%
and the $\mathbb{E}^{a,x}\left( ...\right) $ term in the last line is
estimated using the Markov property as follows.%
\begin{multline*}
\mathbb{E}^{a,x}\left( \prod_{k=i}^{j-1}\exp \left( 24C\eta \frac{\left\Vert 
\mathbf{X}_{t_{k},t_{k+1}}\right\Vert ^{2}}{t_{j}-t_{i}}\right) \right) \\
\left. 
\begin{array}{l}
\leq \prod_{k=i}^{j-1}\sup_{x\in \mathbb{R}^{d}}\mathbb{E}\left( \exp \left(
24C\eta \frac{t_{k+1}-t_{k}}{t_{j}-t_{i}}\frac{\left\Vert \mathbf{X}%
_{0,t_{k+1}-t_{k}}^{x}\right\Vert ^{2}}{t_{k+1}-t_{k}}\right) \right) \\ 
\leq \prod_{k=i}^{j-1}M_{24C\eta \frac{t_{k+1}-t_{k}}{t_{j}-t_{i}}} \\ 
\leq \prod_{k=i}^{j-1}\exp \left( C_{\ref{CorForUniBoundsOnPWLinearApprox}%
}\times 24C\eta \frac{t_{k+1}-t_{k}}{t_{j}-t_{i}}\right) \text{ \ \ for }%
\eta \text{ small enough } \\ 
=\exp \left( 24C_{\ref{CorForUniBoundsOnPWLinearApprox}}C\eta \right)
<\infty .%
\end{array}%
\right.
\end{multline*}%
where we used Corollary \ref{CorForUniBoundsOnPWLinearApprox}, valid for $%
\eta $ small enough. The proof is finished.
\end{proof}

\paragraph{The Subtle Cancellation}

Let us define 
\begin{equation*}
r_{\delta }\left( t,x\right) =\frac{1}{\delta }\mathbb{E}^{a,x}\left(
A_{t,t+\delta }\right) \in so\left( d\right) \text{ \ \ \ \ and \ \ \ \ \ }%
r_{\delta }\left( x\right) =r_{\delta }\left( 0,x\right) .
\end{equation*}%
For instance, (\ref{TooStrongAssumptionWongZakai}) is now expressed as $%
\lim_{\delta \rightarrow 0}r_{\delta }\left( x\right) \rightarrow 0$
uniformly in $x$. Our goal here is to establish a weak version of this. We
also recall that%
\begin{equation*}
A_{t,t+\delta }=\pi _{2}\left( \mathbf{X}_{t,t+\delta }\right) =\pi
_{2}\left( \mathbf{X}_{t}^{-1}\ast \mathbf{X}_{t+\delta }\right) .
\end{equation*}

\begin{proposition}
\label{weakConvegenceArea0}(i) We have uniform boundedness of $r_{\delta
;t}\left( x\right) ,$ 
\begin{equation*}
\sup_{x\in g^{2}\left( \mathbb{R}^{d}\right) }\sup_{\delta \in \left[ 0,1%
\right] }\sup_{t\in \lbrack 0,1-\delta ]}r_{\delta }\left( t,x\right)
<\infty .
\end{equation*}%
(ii) For all $h\in L^{1}\left( g^{2}\left( \mathbb{R}^{d}\right) ,dx\right)
, $%
\begin{equation*}
\lim_{\delta \rightarrow 0}\int_{g^{2}\left( \mathbb{R}^{d}\right)
}dxh\left( x\right) r_{\delta }\left( x\right) \equiv 0.\,
\end{equation*}
\end{proposition}

\begin{proof}
(i) follows from Lemma \ref{GaussianIntegrabilityMarkov}. For (ii) we may
consider $h$ smooth and compactly supported. Now the problem is local and we
can assume that smooth locally bounded functions such as the coordinate
projections $\pi _{1;j}$ and $\pi _{2;k,l}$ are in $D\left( \mathcal{E}%
^{a}\right) $. (More formally, we could smoothly truncate outside the
support of\ $h$ and work on a big torus). Clearly, it is enough to show the
componentwise statement%
\begin{equation*}
\lim_{\delta \rightarrow 0}\int_{g^{2}\left( \mathbb{R}^{d}\right)
}dxh\left( x\right) \pi _{2;k,l}\left( r_{\delta }\left( x\right) \right)
\equiv 0
\end{equation*}%
for $k<l$ fixed in $\left\{ 1,...,d\right\} $. To keep notation short we set 
$f\equiv $ $\pi _{2;k,l}\left( \cdot \right) $ and abuse notation by writing 
$A$ instead of $A^{k,l}$. We can then write%
\begin{equation*}
\mathbb{E}^{a,\mathbf{\cdot }}\left( A_{t}\right) \equiv \mathbb{E}^{a,%
\mathbf{\cdot }}\left( f\left( \mathbf{X}_{t}\right) \right)
=:P_{t}^{a}f\left( .\right)
\end{equation*}%
and note that $P_{0}^{a}f\left( x\right) =A$ when $x=\left( x^{1},A\right)
\in g^{2}\left( \mathbb{R}^{d}\right) $. Writing $\left\langle \cdot ,\cdot
\right\rangle $ for the usual inner product on $L^{2}\left( g^{2}\left( 
\mathbb{R}^{d}\right) ,dx\right) $ we have%
\begin{eqnarray*}
\left\langle h,\mathbb{E}^{a,.}A_{0,t}\right\rangle &=&\left\langle h,%
\mathbb{E}^{a,\mathbf{\cdot }}f\left( \mathbf{X}_{t}\right) -A-\frac{1}{2}%
\mathbb{E}^{a,\mathbf{\cdot }}\left( [\cdot ,\mathbf{X}_{t}^{1}]\right)
\right\rangle \\
&=&\left\langle h,P_{t}^{a}f-P_{0}^{a}f\right\rangle -\left\langle h,\frac{1%
}{2}\mathbb{E}^{a,\mathbf{\cdot }}\left( [\cdot ,\mathbf{X}_{t}^{1}]\right)
\right\rangle \\
&=&\int_{0}^{t}\mathcal{E}^{a}\left( h,P_{s}^{a}f\right) -\left\langle h,%
\frac{1}{2}\mathbb{E}^{a,\mathbf{\cdot }}\left( [\cdot ,\mathbf{X}%
_{t}^{1}]\right) \right\rangle \\
&=&\mathcal{E}^{a}\left( h,f\right) \times t-\left\langle h,\frac{1}{2}%
\mathbb{E}^{a,\mathbf{\cdot }}\left( [\cdot ,\mathbf{X}_{t}^{1}]\right)
\right\rangle +o\left( t\right) .
\end{eqnarray*}%
Here, again, we abused notation by writing $\left[ \cdot ,\cdot \right] $
instead of picking out the $\left( k,l\right) $ component and using the
cumbersome notation $\left[ \cdot ,\cdot \right] ^{k,l}$. Note that in
general $\mathcal{E}^{a}\left( h,f\right) \times t\neq o\left( t\right) $
and our only hope is cancellation of $2\mathcal{E}^{a}\left( h,f\right) $
with the bracket term%
\begin{equation*}
\left\langle h,\mathbb{E}^{a,\mathbf{\cdot }}\left( [\cdot ,\mathbf{X}%
_{t}^{1}]\right) \right\rangle \equiv \left\langle h,\mathbb{E}^{a,\mathbf{%
\cdot }}\left( [\cdot ,\mathbf{X}_{t}^{1}]^{k,l}\right) \right\rangle .
\end{equation*}%
To see this cancellation, we compute the bracket term,%
\begin{eqnarray*}
\left\langle h,\mathbb{E}^{a,\mathbf{\cdot }}\left( [\cdot ,\mathbf{X}%
_{t}^{1}]^{k,l}\right) \right\rangle &=&\int dx\,h\left( x\right) \mathbb{E}%
^{a,x}\left( x^{1;k}\mathbf{X}_{t}^{1;l}-x^{1;l}\mathbf{X}_{t}^{1;k}\right)
\\
&=&\int dx\,h\left( x\right) \left( \left( x^{1;k}\left[ P_{t}^{a}\pi _{1;l}%
\right] \left( x\right) -x^{1;l}\left[ P_{t}^{a}\pi _{1;k}\right] \left(
x\right) \right) \right) ,
\end{eqnarray*}%
and by adding and subtracting $x^{1;k}x^{1;l}$ inside the integral this
rewrites as%
\begin{equation*}
\int dx\,h\left( x\right) x^{1;k}\left\{ \left[ P_{t}^{a}\pi _{1;l}\right]
\left( x\right) -\pi _{1;l}\left( x\right) \right\} -\int dx\,h\left(
x\right) x^{1;l}\left\{ \left[ P_{t}^{a}\pi _{1;k}\right] \left( x\right)
-\pi _{1;k}\left( x\right) \right\} .
\end{equation*}%
It now follows as earlier that%
\begin{equation*}
\left\langle h,\mathbb{E}^{a,\mathbf{\cdot }}\left( [\cdot ,\mathbf{X}%
_{t}^{1}]^{k,l}\right) \right\rangle =\left[ \mathcal{E}^{a}\left( h\pi
_{1;k},\pi _{1;l}\right) -\mathcal{E}^{a}\left( h\pi _{1;l},\pi
_{1;k}\right) \right] \times t+o\left( t\right)
\end{equation*}%
and we see that the required cancellation takes place if, for all $h$ smooth
and compactly supported,%
\begin{equation*}
\left[ \mathcal{E}^{a}\left( h\pi _{1;k},\pi _{1;l}\right) -\mathcal{E}%
^{a}\left( h\pi _{1;l},\pi _{1;k}\right) \right] \equiv 2\mathcal{E}%
^{a}\left( h,\pi _{2;k,l}\right) .
\end{equation*}%
We will check this with a direct computation. First note that%
\begin{equation*}
\mathcal{E}^{a}\left( h\pi _{1;k},\pi _{1;l}\right) -\mathcal{E}^{a}\left(
h\pi _{1;l},\pi _{1;k}\right) =\int \pi _{1,k}d\Gamma ^{a}\left( h,\pi
_{1,l}\right) -\int \pi _{1,l}d\Gamma ^{a}\left( h,\pi _{1,k}\right)
\end{equation*}%
which is immediately seen via symmetry of $d\Gamma ^{a}\left( \cdot ,\cdot
\right) $, inherited from the symmetric of $\left( a^{ij}\right) ,$ and the
Leibnitz formula%
\begin{equation*}
\mathcal{E}^{a}\left( gg^{\prime },h\right) =\int gd\Gamma ^{a}\left(
g^{\prime },h\right) +\int g^{\prime }d\Gamma ^{a}\left( g,h\right) .
\end{equation*}%
It is immediately checked from the definition of the vector fields $U_{i}$,
see equation (\ref{DefOf_U_on_g2}), that%
\begin{equation*}
U_{i}f\equiv U_{i}\pi _{2;k,l}=\left\{ 
\begin{array}{c}
\begin{array}{c}
-\left( 1/2\right) \pi _{1;l}\text{ if }i=k \\ 
\left( 1/2\right) \pi _{1;k}\text{ if }i=l%
\end{array}
\\ 
0\text{ otherwise}%
\end{array}%
\right.
\end{equation*}%
so that 
\begin{equation*}
\int \pi _{1,k}d\Gamma ^{a}\left( h,\pi _{1,l}\right) =\sum_{i,j}\int \pi
_{1,k}a^{ij}U_{i}hU_{j}\pi _{1,l}=2\sum_{i}\int \left( U_{l}f\right)
a^{il}\left( U_{i}h\right)
\end{equation*}%
and similarly%
\begin{equation*}
-\int \pi _{1,l}d\Gamma ^{a}\left( h,\pi _{1,k}\right) =\sum_{i,j}\int
\left( -\pi _{1,l}\right) a^{ij}U_{i}hU_{j}\pi _{1,k}=2\sum_{i}\int \left(
U_{k}f\right) a^{ik}\left( U_{i}h\right) .
\end{equation*}%
Therefore, using $U_{j}f=0$ for $j\neq \left\{ k,l\right\} $ in the second
equality, 
\begin{eqnarray*}
\mathcal{E}^{a}\left( h\pi _{1;k},\pi _{1;l}\right) -\mathcal{E}^{a}\left(
h\pi _{1;l},\pi _{1;k}\right) &=&2\sum_{j=k,l}\sum_{i}\int \left(
U_{j}f\right) a^{ij}\left( U_{i}h\right) \\
&=&2\sum_{i,j}\int \left( U_{j}f\right) a^{ij}\left( U_{i}h\right)
\end{eqnarray*}%
and this equals precisely $2\mathcal{E}^{a}\left( h,f\right) $ as required.
\end{proof}

\begin{corollary}
\label{weakShiftedT}For all $t\in \lbrack 0,1)$ and all $h\in L^{1}\left(
g^{2}\left( \mathbb{R}^{d}\right) ,dx\right) ,$%
\begin{equation*}
\lim_{\delta \rightarrow 0}\int_{g^{2}\left( \mathbb{R}^{d}\right)
}dxh\left( x\right) \mathbb{E}^{a,x}\left( \frac{A_{t,t+\delta }}{\delta }%
\right) \equiv 0.
\end{equation*}
\end{corollary}

\begin{proof}
We first write%
\begin{eqnarray*}
\int dxh\left( x\right) \mathbb{E}^{a,x}\left( \frac{A_{t,t+\delta }}{\delta 
}\right) &=&\int \int h\left( x\right) p^{a}\left( t,x,y\right) r_{\delta
}\left( y\right) dxdy \\
&=&\int \left( \int h\left( x\right) p^{a}\left( t,x,y\right) dx\right)
r_{\delta }\left( y\right) dy\mathbf{.}
\end{eqnarray*}%
Then, noting that $y\mapsto \int h\left( x\right) p_{t}\left( x,y\right) dx$
is in $L^{1}\left( g^{2}\left( \mathbb{R}^{d}\right) ,dx\right) $, the proof
is finished by applying the previous proposition.
\end{proof}

\begin{theorem}
\label{AreaUniformTo0}For all bounded sets $K\subset g^{2}\left( \mathbb{R}%
^{d}\right) $ and all $\sigma \in (0,1],$%
\begin{equation*}
\lim_{\delta \rightarrow 0}\sup_{t\in \left[ \sigma ,1\right] }\sup_{y\in
K}\left\vert \mathbb{E}^{a,y}\left( \frac{A_{t,t+\delta }}{\delta }\right)
\right\vert =0.
\end{equation*}
\end{theorem}

\begin{proof}
It suffices to prove this for a compact ball $K=\bar{B}\left( 0,R\right)
\subset $ $g^{2}\left( \mathbb{R}^{d}\right) $ of arbitrary radius $R>0$. We
fix $\sigma \in (0,1]$ and think of $r_{\delta }=r_{\delta }\left(
t,y\right) $ as a family of maps, indexed by $\delta >0$, defined on the
cylinder $\left[ \sigma ,1\right] \times K$, that is 
\begin{equation*}
\left( t,y\right) \in \left[ \sigma ,1\right] \times K\mapsto r_{\delta
}\left( t,y\right) \in so\left( d\right) .
\end{equation*}%
By Proposition \ref{weakConvegenceArea0}, (i) we know that $\sup_{\delta
>0}\left\vert r_{\delta }\right\vert _{\infty }<\infty $. We now show
equicontinuity of $\left\{ r_{\delta }:\delta >0\right\} $. By the Markov
property, $r_{\delta }\left( t,y\right) $ equals%
\begin{equation*}
\mathbb{E}^{a,y}\left( \frac{A_{t,t+\delta }}{\delta }\right) =\left\langle
p^{a}\left( t,y,\cdot \right) ,\frac{\mathbb{E}^{a,\cdot }\left( A_{0,\delta
}\right) }{\delta }\right\rangle =\left\langle p^{a}\left( t,y,\cdot \right)
,r_{\delta }\left( 0,\cdot \right) \right\rangle ,
\end{equation*}%
so that, for all $\left( s,x\right) ,\left( t,y\right) \in \left[ \sigma ,1%
\right] \times K$,%
\begin{eqnarray*}
\left\vert r_{\delta }\left( s,x\right) -r_{\delta }\left( t,y\right)
\right\vert &=&\left\vert \left\langle p^{a}\left( s,x,\cdot \right)
-p^{a}\left( t,y,\cdot \right) ,r_{\delta }\left( .\right) \right\rangle
\right\vert \\
&\leq &\left( \sup_{\delta \in (0,1]}\left\vert r_{\delta }\right\vert
_{\infty }\right) \left\vert p^{a}\left( s,x,\cdot \right) -p^{a}\left(
t,y,\cdot \right) \right\vert _{L^{1}}.
\end{eqnarray*}%
From Proposition \ref{MoserNash}, $\left( t,y\right) \in \left[ \sigma ,1%
\right] \times K\mapsto p^{a}\left( t,y,z\right) $ is continuous for all $z$%
; the dominated convergence theorem then gives easily continuity of $\left(
t,y\right) \mapsto p^{a}\left( t,y,\cdot \right) \in L^{1}$. In fact, this
map is uniformly continuous when restricted to the compact $\left[ \sigma ,1%
\right] \times K$ and it follows that $\left\{ r_{\delta }:\delta >0\right\} 
$ is equicontinuous as claimed. By Arzela-Ascoli, there exists a subsequence 
$\left( \delta ^{n}\right) $ such that $r_{\delta ^{n}}$ converges uniformly
on $\left[ \sigma ,1\right] \times K$ to some (continuous) function $r.$ On
the other hand, Proposition \ref{weakConvegenceArea0}, (ii), applied to $h=$ 
$p^{a}\left( t,y,\cdot \right) $, shows that $r_{\delta }\left( t,y\right)
\rightarrow 0$ as $\delta \rightarrow 0$ for all fixed $y,t>0$. This shows
that $r\equiv 0$ is the only limit point and hence%
\begin{equation*}
\lim_{\delta \rightarrow 0}\sup_{t\in \left[ \sigma ,1\right] }\sup_{y\in
K}\left\vert \mathbb{E}^{a,y}\left( \frac{A_{t,t+\delta }}{\delta }\right)
\right\vert =0.
\end{equation*}
\end{proof}

\paragraph{Convergence of the Sum of the Small Areas}

For fixed $a\in \Xi \left( \Lambda \right) $ and $x\in g^{2}\left( \mathbb{R}%
^{d}\right) $ let us define the real-valued quantity%
\begin{equation*}
K_{\sigma ,\delta }:=\sup_{\substack{ 0\leq u_{1}<u_{2}<v_{1}<v_{2}\leq 1: 
\\ v_{1}-u_{2}\geq \sigma ,  \\ \left\vert u_{2}-u_{1}\right\vert
,\left\vert v_{2}-v_{1}\right\vert \leq \delta }}\frac{\left\vert \mathbb{E}%
^{a,x}\left( A_{u_{1},u_{2}}\cdot A_{v_{1},v_{2}}\right) \right\vert }{%
\left( u_{2}-u_{1}\right) \left( v_{2}-v_{1}\right) }\text{ }
\end{equation*}%
where $\delta ,\sigma \in \left( 0,1\right) $. As above $\cdot $ denotes the
scalar product in $so\left( d\right) $.

\begin{proposition}
\label{productAreaTo0}For fixed $\sigma \in (0,1),$ $k,l\in \left\{
1,..,d\right\} $ we have $\lim_{\delta \rightarrow 0}K_{\sigma ,\delta }=0$.
\end{proposition}

\begin{proof}
By the Markov property,%
\begin{eqnarray*}
\frac{\left\vert \mathbb{E}^{a,x}\left( A_{u_{1},u_{2}}\cdot
A_{v_{1},v_{2}}\right) \right\vert }{\left( u_{2}-u_{1}\right) \left(
v_{2}-v_{1}\right) } &=&\frac{\left\vert \mathbb{E}^{a,x}\left(
A_{u_{1},u_{2}}\cdot \mathbb{E}^{a,\mathbf{X}_{u_{2}}}\left(
A_{v_{1}-u_{2},v_{2}-u_{2}}\right) \right) \right\vert }{\left(
u_{2}-u_{1}\right) \left( v_{2}-v_{1}\right) } \\
&\leq &\frac{\left\vert \mathbb{E}^{a,x}\left( A_{u_{1},u_{2}}\cdot \mathbb{E%
}^{a,\mathbf{X}_{u_{2}}}\left( A_{v_{1}-u_{2},v_{2}-u_{2}};\left\Vert 
\mathbf{X}_{u_{2}}\right\Vert \leq R\right) \right) \right\vert }{\left(
u_{2}-u_{1}\right) \left( v_{2}-v_{1}\right) } \\
&&+\frac{\left\vert \mathbb{E}^{a,x}\left( A_{u_{1},u_{2}}\cdot \mathbb{E}%
^{a,\mathbf{X}_{u_{2}}}\left( A_{v_{1}-u_{2},v_{2}-u_{2}};\left\Vert \mathbf{%
X}_{u_{2}}\right\Vert >R\right) \right) \right\vert }{\left(
u_{2}-u_{1}\right) \left( v_{2}-v_{1}\right) } \\
&\leq &\frac{\mathbb{E}^{a,x}\left( \left\vert A_{u_{1},u_{2}}\right\vert
;\left\Vert \mathbf{X}_{u_{2}}\right\Vert \leq R\right) }{\left(
u_{2}-u_{1}\right) }\sup_{\delta ^{\prime }\leq \delta }\sup_{\substack{ %
\left\Vert y\right\Vert \leq R  \\ u\in \lbrack \sigma ,1]}}\frac{\left\vert 
\mathbb{E}^{a,y}\left( A_{u,u+\delta ^{\prime }}\right) \right\vert }{\delta
^{\prime }} \\
&&+\mathbb{E}^{a,x}\left( \frac{\left\vert A_{u_{1},u_{2}}\right\vert }{%
u_{2}-u_{1}};\left\Vert \mathbf{X}_{u_{2}}\right\Vert >R\right) \sup_{\delta
^{\prime },u,x}\frac{\mathbb{E}^{a,x}\left( \left\vert A_{u,u+\delta
^{\prime }}\right\vert \right) }{\delta ^{\prime }}. \\
&\leq &\frac{\mathbb{E}^{a,x}\left( \left\vert A_{u_{1},u_{2}}\right\vert
\right) }{\left( u_{2}-u_{1}\right) }\sup_{\delta ^{\prime }\leq \delta
}\sup _{\substack{ \left\Vert y\right\Vert \leq R  \\ u\in \lbrack \sigma
,1] }}\frac{\left\vert \mathbb{E}^{a,y}\left( A_{u,u+\delta ^{\prime
}}\right) \right\vert }{\delta ^{\prime }} \\
&&+\sqrt{\mathbb{P}^{a,x}\left( \left\Vert \mathbf{X}_{u_{2}}\right\Vert
>R\right) }\sqrt{\mathbb{E}^{a,x}\left( \left\vert \frac{A_{u_{1},u_{2}}}{%
u_{2}-u_{1}}\right\vert ^{2}\right) }\sup_{\delta ^{\prime },u,x}\frac{%
\mathbb{E}^{a,x}\left( \left\vert A_{u,u+\delta ^{\prime }}\right\vert
\right) }{\delta ^{\prime }} \\
&\leq &C\sup_{\delta ^{\prime }\leq \delta }\sup_{\substack{ \left\vert
y\right\vert \leq R  \\ u\in \lbrack \sigma ,1]}}\frac{\left\vert \mathbb{E}%
^{a,y}\left( A_{u,u+\delta ^{\prime }}\right) \right\vert }{\delta ^{\prime }%
}+C\sqrt{\mathbb{P}^{a,x}\left( \left\Vert \mathbf{X}_{u_{2}}\right\Vert
>R\right) }
\end{eqnarray*}%
for some constant $C=C\left( \left\Vert x\right\Vert ,\sigma ,\Lambda
\right) $ using Lemma \ref{GaussianIntegrabilityMarkov} and Proposition \ref%
{weakConvegenceArea0}, (i). We then fix $\varepsilon >0$ and choose $%
R=R\left( \epsilon \right) $ large enough so that%
\begin{equation*}
C\sup_{u_{2}\in \left[ 0,1\right] }\sqrt{\mathbb{E}^{a,x}\left( \left\vert 
\mathbf{X}_{u_{2}}\right\vert >R\right) }\leq \varepsilon /2.
\end{equation*}%
On the other hand, Theorem \ref{AreaUniformTo0} shows that 
\begin{equation*}
C\sup_{\delta ^{\prime }\leq \delta }\sup_{\substack{ \left\vert
y\right\vert \leq R  \\ u\in \lbrack \sigma ,1]}}\frac{\left\vert \mathbb{E}%
^{a,y}\left( A_{u,u+\delta ^{\prime }}\right) \right\vert }{\delta ^{\prime }%
}\leq \frac{\varepsilon }{2}
\end{equation*}%
for all $\delta $ small enough and the proof is finished.
\end{proof}

\begin{corollary}
\label{preWongZakaiInL4}There exists $C_{\ref{preWongZakaiInL4}}=C_{\ref%
{preWongZakaiInL4}}\left( \Lambda \right) $ such that for all subdivisions $%
D $ of $\left[ 0,1\right] ,$ $s,t\in D,$ for any $\sigma \in \left(
0,1\right) ,$ 
\begin{equation*}
\mathbb{E}^{a,x}\left( \left\vert d\left( S\left( X^{D}\right) _{s,t},%
\mathbf{X}_{s,t}\right) \right\vert ^{4}\right) \leq C_{\ref%
{preWongZakaiInL4}}\left[ \left( t-s\right) ^{2}K_{\sigma ,\left\vert
D\right\vert }+\left( t-s\right) \sigma \right] .
\end{equation*}
\end{corollary}

\begin{proof}
Recalling the discussion around (\ref{SumOfLittleAreas}), equivalence of
homogenous norms leads to%
\begin{equation*}
\mathbb{E}^{a,x}\left( \left\vert d\left( S\left( X^{D}\right) _{s,t},%
\mathbf{X}_{s,t}\right) \right\vert ^{4}\right) \leq C\mathbb{E}%
^{a,x}(|\sum_{i:t_{i}\in D\cap \lbrack s,t)}A_{t_{i},t_{i+1}}|^{2}).
\end{equation*}%
Let us abbreviate $\sum_{i:t_{i}\in D\cap \lbrack s,t)}$ to $\sum_{i}$ in
what follows. Clearly, $\mathbb{E}^{a,x}(|\sum_{i}A_{t_{i},t_{i+1}}|^{2})$
is estimated by $2$ times%
\begin{eqnarray*}
&&\sum_{i\leq j}\mathbb{E}^{a,x}\left( A_{t_{i},t_{i+1}}\cdot
A_{t_{j},t_{j+1}}\right) \\
&\leq &\sum_{\substack{ i\leq j  \\ t_{j}-t_{i+1}\geq \sigma }}\mathbb{E}%
^{a,x}\left( A_{t_{i},t_{i+1}}\cdot A_{t_{j},t_{j+1}}\right) +\sum 
_{\substack{ i\leq j  \\ t_{j}-t_{i+1}<\sigma }}\mathbb{E}^{a,x}\left(
A_{t_{i},t_{i+1}}\cdot A_{t_{j},t_{j+1}}\right) \\
&\leq &K_{\sigma ,\left\vert D\right\vert }\sum_{\substack{ i\leq j  \\ %
t_{j}-t_{i+1}\geq \sigma }}\left( t_{i+1}-t_{i}\right) \left(
t_{j+1}-t_{j}\right) +\sum_{\substack{ i\leq j  \\ t_{j}-t_{i+1}<\sigma }}%
\sqrt{\mathbb{E}^{a,x}\left( \left\vert A_{t_{i},t_{i+1}}\right\vert
^{2}\right) \mathbb{E}^{a,x}\left( \left\vert A_{t_{j},t_{j+1}}\right\vert
^{2}\right) } \\
&\leq &K_{\sigma ,\left\vert D\right\vert }\left( t-s\right) ^{2}+C\sum 
_{\substack{ i,j  \\ t_{j}-t_{i+1}<\sigma }}\left( t_{i+1}-t_{i}\right)
\left( t_{j+1}-t_{j}\right)
\end{eqnarray*}%
and the very last sum is estimated as follows,%
\begin{equation*}
|\sum_{i}\left( t_{i+1}-t_{i}\right) \sum_{\substack{ j  \\ %
t_{j}-t_{i+1}<\sigma }}\left( t_{j+1}-t_{j}\right) |\leq \sigma
\sum_{i}\left( t_{i+1}-t_{i}\right) =\sigma \left( t-s\right) .
\end{equation*}%
The proof is finished.
\end{proof}

\paragraph{Putting Things Together}

\begin{theorem}
Let $D$ be a dissection of $\left[ 0,1\right] $ with mesh $\left\vert
D\right\vert .$Then, for all $1\leq q<\infty $ and $0\leq \alpha <1/2,$%
\begin{equation*}
d_{\alpha \text{-H\"{o}lder}}\left( S\left( X^{D}\right) ,\mathbf{X}\right)
\rightarrow 0\text{ in }L^{q}\left( \mathbb{P}^{a,x}\right) \text{ as }%
\left\vert D\right\vert \rightarrow 0.
\end{equation*}
\end{theorem}

\begin{proof}
We first show pointwise convergence. We fix $\varepsilon >0$ and apply
Corollary \ref{preWongZakaiInL4} with $\sigma =\varepsilon /2C.$ Then,%
\begin{equation*}
\sup_{s,t\in D:s<t}\mathbb{E}^{a,x}\left( \left\vert d\left( S\left(
X^{D}\right) _{s,t},\mathbf{X}_{s,t}\right) \right\vert ^{4}\right) \leq
CK_{\sigma ,\left\vert D\right\vert }+\frac{\varepsilon }{2}
\end{equation*}%
By Proposition \ref{productAreaTo0} it then follows that, for $\left\vert
D\right\vert $ small enough, 
\begin{equation*}
\sup_{s,t\in D:s<t}\left\Vert d\left( S\left( X^{D}\right) _{s,t},\mathbf{X}%
_{s,t}\right) \right\Vert _{L^{4}\left( \mathbb{P}^{a,x}\right) }^{4}\leq
\varepsilon .
\end{equation*}%
By Theorem \ref{UniformGaussTailOfApproximations} we have for all $q\in
\lbrack 1,\infty ),$%
\begin{equation}
\sup_{D}\left\Vert \left\Vert S\left( X^{D}\right) \right\Vert _{\alpha 
\text{-H\"{o}lder}}\right\Vert _{L^{q}\left( \mathbb{P}^{a,x}\right)
}+\left\Vert \left\Vert \mathbf{X}\right\Vert _{\alpha \text{-H\"{o}lder}%
}\right\Vert _{L^{q}\left( \mathbb{P}^{a,x}\right) }<\infty  \label{UGToA}
\end{equation}%
and both results combined yield%
\begin{equation*}
\lim_{\left\vert D\right\vert \rightarrow 0}\sup_{0\leq s<t\leq 1}\left\Vert
d\left( S\left( X^{D}\right) _{s,t},\mathbf{X}_{s,t}\right) \right\Vert
_{L^{4}\left( \mathbb{P}^{a,x}\right) }=0
\end{equation*}%
and by H\"{o}lder's inequality the last statement remains valid even when we
replace $L^{4}$ by $L^{q}$ for any $q\in \lbrack 1,\infty )$. Now, for every 
$m>0,$%
\begin{eqnarray*}
\mathbb{E}^{a,x}\left( d_{\infty }\left( S\left( X^{D}\right) ,\mathbf{X}%
\right) ^{q}\right) &\leq &c_{q}\mathbb{E}^{a,x}\left( \sup_{1\leq i\leq
m}d\left( S\left( X^{D}\right) _{\frac{i}{m}},\mathbf{X}_{\frac{i}{m}%
}\right) ^{q}\right) \\
&&+c_{q}\mathbb{E}^{a,x}\left( \sup_{\left\vert t-s\right\vert <\frac{1}{m}%
}\left( \left\Vert S\left( X^{D}\right) _{s,t}\right\Vert ^{q}+\left\Vert 
\mathbf{X}_{s,t}\right\Vert ^{q}\right) \right) \\
&\leq &c_{q}m\sup_{0\leq t\leq 1}\left\Vert d\left( S\left( X^{D}\right)
_{t},\mathbf{X}_{t}\right) \right\Vert _{L^{q}\left( \mathbb{P}^{a,x}\right)
}^{q} \\
&&+c_{q}\left( \frac{1}{m}\right) ^{\alpha q}\mathbb{E}^{a,x}\left( \left(
\left\Vert S\left( X^{D}\right) \right\Vert _{\alpha \text{-H\"{o}lder}%
}^{q}+\left\Vert \mathbf{X}\right\Vert _{\alpha \text{-H\"{o}lder}%
}^{q}\right) \right) \\
&\leq &c_{q}m\sup_{0\leq t\leq 1}\left\Vert d\left( S\left( X^{D}\right)
_{t},\mathbf{X}_{t}\right) \right\Vert _{L^{q}\left( \mathbb{P}^{a,x}\right)
}^{q}+C\left( \frac{1}{m}\right) ^{\alpha q}.
\end{eqnarray*}%
By choosing first $m$ large enough and then $D$ with $\left\vert
D\right\vert $ small enough we see that $d_{\infty }\left( S\left(
X^{D}\right) ,\mathbf{X}\right) \rightarrow 0$ in $L^{q}$ as $\left\vert
D\right\vert \rightarrow 0$, for all $q<\infty $. An easy application of the
Campell-Hausdorff formula gives a $d_{0}/d_{\infty }$-estimate, 
\begin{equation*}
\forall \mathbf{x,y}\in C\left( \left[ 0,1\right] ,g^{2}\left( \mathbb{R}%
^{d}\right) \right) :d_{0}\left( \mathbf{x,y}\right) \leq d_{\infty }\left( 
\mathbf{x,y}\right) +C\sqrt{\left\Vert \mathbf{y}\right\Vert _{\infty
}d_{\infty }\left( \mathbf{x,y}\right) }\text{.}
\end{equation*}%
With Cauchy-Schwarz and a standard H\"{o}lder interpolation argument, using (%
\ref{UGToA}) with $\alpha ^{\prime }\in (\alpha ,1/2)$, we then see that%
\begin{equation*}
d_{\alpha \text{-H\"{o}lder}}\left( S\left( X^{D}\right) ,\mathbf{X}\right)
\rightarrow 0\text{ in }L^{q}\left( \mathbb{P}^{a,x}\right) \text{ as }%
\left\vert D\right\vert \rightarrow 0.
\end{equation*}
\end{proof}

\begin{remark}
This convergence result implies that $\sigma \left( A_{s,t}:u\leq s\leq
t\leq v\right) \subset \mathcal{F}_{s,t}=\sigma \left( X_{s,r}:s\leq r\leq
t\right) $ where $X=\pi _{1}\left( \mathbf{X}\right) $ and $%
X_{s,r}=X_{r}-X_{s}\in \mathbb{R}^{d}$.
\end{remark}

\begin{corollary}
Let $Y=\pi \left( 0,y_{0};\mathbf{X}\right) \equiv \pi \left( \mathbf{X}%
\right) $ denote the $\mathbb{R}^{e}$-valued (random) RDE solution driven by 
$\mathbf{X}^{a,x}$ along fixed $\mathrm{Lip}^{\gamma }$ vector fields $%
V_{1},...,V_{d}$ on $\mathbb{R}^{e}$, with $\gamma >2$, and started at time $%
0$ from $y_{0}$ fixed. Let $Y^{D}=\pi \left( 0,y_{0},X^{D}\right) $ be the
piecewise smooth solution to corresponding control ODE%
\begin{equation*}
dY^{D}=\sum_{i=1}^{d}V_{i}\left( Y^{D}\right) dX^{D;i}.
\end{equation*}%
Then for any $\alpha \in \lbrack 0,1/2)$ we have $\left\vert
Y-Y^{D}\right\vert _{\alpha \text{-H\"{o}l;}\left[ 0,1\right] }\rightarrow 0$
in $L^{q}\left( \mathbb{P}^{a,x}\right) ,$ for all $q<\infty $.
\end{corollary}

\begin{proof}
The universal limit theorem \cite{lyons-98, lyons-qian-02} shows immediately
that%
\begin{equation*}
\left\vert Y-Y^{D}\right\vert _{\alpha \text{-H\"{o}l;}\left[ 0,1\right]
}\rightarrow 0\text{ in probability.}
\end{equation*}%
It then suffices to remark that the estimates on the It\^{o}-Lyons map in 
\cite{friz-victoir-06-Euler} combined with Theorem \ref%
{UniformGaussTailOfApproximations} show that for all $q<\infty ,$%
\begin{equation*}
\sup_{D}\mathbb{E}\left\vert Y^{D}\right\vert _{\alpha \text{-H\"{o}l;}\left[
0,1\right] }^{q},\,\ \left\vert Y\right\vert _{\alpha \text{-H\"{o}l;}\left[
0,1\right] }^{q}<\infty .
\end{equation*}%
\ 
\end{proof}

\section{RDE Solutions as Markov Processes}

The following is an immediate consequence of the stochastic Taylor formula
for random RDEs \cite{friz-victoir-06-Euler}.

\begin{lemma}
\label{EulerEstimates}Let $\alpha \in \left( 1/3,1/2\right) $, $N=2$. Assume
the random rough path $\mathbf{X}$ is such that $\left\Vert \mathbf{X}%
\right\Vert _{\alpha \text{-H\"{o}l;}\left[ 0,1\right] }<\infty $ has a
Gauss tail and let $Z$ denote the random RDE solution driven by $\mathbf{X}$
along fixed $\mathrm{Lip}^{\gamma }$ vector fields $V_{1},...,V_{d}$, with $%
\gamma >2$, and started from $z$. Then for all $f\in C_{b}^{\infty }$ we have%
\begin{eqnarray*}
\mathbb{E}\left[ f\left( Z_{t}\right) \right] &=&f\left( z\right)
+\sum_{i=1}^{d}V_{i}f\left( z\right) \mathbb{E}\left[ \pi _{1,i}\left( 
\mathbf{X}_{0,t}\right) \right] \\
&&+\frac{1}{2}\sum_{i,j=1}^{d}V_{i}V_{j}f\left( z\right) \mathbb{E}\left[
\pi _{1,i}\left( \mathbf{X}_{0,t}\right) \pi _{1,j}\left( \mathbf{X}%
_{0,t}\right) \right] \\
&&+\frac{1}{2}\sum_{i,j=1}^{d}\left[ V_{i},V_{j}\right] f\left( z\right) 
\mathbb{E}\left[ \pi _{2,i,j}\left( \mathbf{X}_{0,t}\right) \right] +\mathbb{%
E}\left[ R_{2}\left( t,f\right) \right] .
\end{eqnarray*}%
with remainder term,%
\begin{equation*}
\mathbb{E}\left[ \left\vert R_{2}\left( t,f\right) \right\vert \right]
=o\left( t\right) \text{ as }t\rightarrow 0.
\end{equation*}%
As earlier, $\left\langle \cdot ,\cdot \right\rangle $ denotes the inner
product on $L^{2}\left( g^{2}\left( \mathbb{R}^{d}\right) ,m\right) $.
\end{lemma}

\begin{lemma}
\label{weakconvergencemoments}Let $g$ be a compactly supported smooth
function. Then, for all $k,l\in \left\{ 1,...,d\right\} $,%
\begin{eqnarray*}
\lim_{t\rightarrow 0}\left\langle g,\frac{\mathbb{E}^{a,.}\left[ \pi
_{1,k}\left( \mathbf{X}_{0,t}\right) \right] }{t}\right\rangle
&=&\sum_{j=1}^{d}\int_{g^{2}\left( \mathbb{R}^{d}\right) }a^{kj}\left(
y\right) U_{j}g\left( y\right) dy, \\
\lim_{t\rightarrow 0}\left\langle g,\frac{\mathbb{E}^{a,.}\left[ \pi
_{1,k}\left( \mathbf{X}_{0,t}\right) \pi _{1,l}\left( \mathbf{X}%
_{0,t}\right) \right] }{t}\right\rangle &=&-2\int_{g^{2}\left( \mathbb{R}%
^{d}\right) }a^{kl}\left( y\right) g\left( y\right) dy, \\
\lim_{t\rightarrow 0}\left\langle g,\frac{\mathbb{E}^{a,.}\left[ \pi
_{2,i,j}\left( \mathbf{X}_{0,t}\right) \right] }{t}\right\rangle &=&0.
\end{eqnarray*}
\end{lemma}

\begin{proof}
Third equality was shown in Proposition \ref{weakConvegenceArea0}. For the
first statement, almost by definition of $\mathcal{E}^{a},$%
\begin{eqnarray*}
\lim_{t\rightarrow 0}\left\langle g,\frac{\mathbb{E}^{a,.}\left[ \pi
_{1,k}\left( \mathbf{X}_{0,t}\right) \right] }{t}\right\rangle &=&\mathcal{E}%
^{a}\left( \pi _{1,k},g\right) \\
&=&\sum_{j=1}^{d}\int_{g^{2}\left( \mathbb{R}^{d}\right) }a^{kj}\left(
y\right) U_{j}g\left( y\right) dy.
\end{eqnarray*}%
Let us now consider the second equality. First rewrite $\pi _{1,k}\left( 
\mathbf{X}_{0,t}\right) \pi _{1,l}\left( \mathbf{X}_{0,t}\right) $ as 
\begin{equation*}
\pi _{1,k}\left( \mathbf{X}_{t}\right) \pi _{1,l}\left( \mathbf{X}%
_{t}\right) -\pi _{1,k}\left( \mathbf{X}_{0}\right) \pi _{1,l}\left( \mathbf{%
X}_{0}\right) -\pi _{1,k}\left( \mathbf{X}_{0}\right) \pi _{1,l}\left( 
\mathbf{X}_{0,t}\right) -\pi _{1,l}\left( \mathbf{X}_{0}\right) \pi
_{1,k}\left( \mathbf{X}_{0,t}\right) .
\end{equation*}%
Then, by a similar argument as above,%
\begin{eqnarray*}
\lim_{t\rightarrow 0}\left\langle g,\frac{\mathbb{E}^{a,.}\left[ \pi
_{1,k}\left( \mathbf{X}_{0,t}\right) \pi _{1,l}\left( \mathbf{X}%
_{0,t}\right) \right] }{t}\right\rangle &=&\mathcal{E}^{a}\left( \pi
_{1,k}\pi _{1,l},g\right) \\
&&-\mathcal{E}^{a}\left( \pi _{1,k},\pi _{1,l}g\right) -\mathcal{E}%
^{a}\left( \pi _{1,l},\pi _{1,k}g\right) .
\end{eqnarray*}%
By the Leibniz formula, recalling that $d\Gamma ^{a}\left( f,f^{\prime
}\right) \equiv \left( \sum_{i,j}a^{ij}U_{i}fU_{j}f^{\prime }\right) dm$ is
the energy measure of $\mathcal{E}^{a}$, we have%
\begin{eqnarray*}
\mathcal{E}^{a}\left( \pi _{1,k}\pi _{1,l},g\right) &=&\int \pi
_{1,l}d\Gamma ^{a}\left( \pi _{1,k},g\right) +\int \pi _{1,k}d\Gamma
^{a}\left( \pi _{1,l},g\right) \\
\mathcal{E}^{a}\left( \pi _{1,k},\pi _{1,l}g\right) &=&\int \pi
_{1,l}d\Gamma ^{a}\left( \pi _{1,k},g\right) +\int gd\Gamma ^{a}\left( \pi
_{1,k},\pi _{1,l}\right) \\
\mathcal{E}^{a}\left( \pi _{1,l},\pi _{1,k}g\right) &=&\int \pi
_{1,k}d\Gamma ^{a}\left( \pi _{1,l},g\right) +\int gd\Gamma ^{a}\left( \pi
_{1,l},\pi _{1,k}\right) .
\end{eqnarray*}%
and using the symmetry of $a$ we see that%
\begin{multline*}
\lim_{t\rightarrow 0}\left\langle g,\frac{\mathbb{E}^{a,.}\left[ \pi
_{1,k}\left( \mathbf{X}_{0,t}\right) \pi _{1,l}\left( \mathbf{X}%
_{0,t}\right) \right] }{t}\right\rangle \\
\left. 
\begin{array}{l}
=-2\sum_{i,j=1}^{d}\int_{g^{2}\left( \mathbb{R}^{d}\right) }a^{ij}\left(
y\right) g\left( y\right) U_{i}\pi _{1,k}\left( y\right) U_{j}\pi
_{1,l}\left( y\right) dy \\ 
=-2\int_{g^{2}\left( \mathbb{R}^{d}\right) }a^{kl}\left( y\right) g\left(
y\right) dy.%
\end{array}%
\right.
\end{multline*}
\end{proof}

Let us fix a collection $V=\left( V_{1},...,V_{d}\right) $ of $\mathrm{Lip}%
^{3}$ vector fields on $\mathbb{R}^{e}$ on let us consider the RDE\footnote{%
Regularity of the vector fields could be improved to $\mathrm{Lip}%
^{2+\epsilon }$. Also, one can easily add a drift term $V_{0}\left( Y\right)
dt$ by considering the canonical space-time rough path $\left( \mathbf{X}%
^{a,x},t\right) $.} 
\begin{equation*}
\left\{ 
\begin{array}{l}
dY=V\left( Y\right) d\mathbf{X}^{a,x} \\ 
Y_{0}=y.%
\end{array}%
\right.
\end{equation*}%
where $Y$ is the $\mathbb{R}^{e}$-valued solution path\footnote{%
We could construct the solution as (random) geometric rough path with values
in $g^{2}\left( \mathbb{R}^{e}\right) $ and the arguments which follow
extend to this case.}.\ In general, $Y$ is not Markov, but it is easy to see
that $Z^{z}=\left( \mathbf{X}^{a,x}\mathbf{,}Y\right) \in g^{2}\left( 
\mathbb{R}^{d}\right) \oplus \mathbb{R}^{e}$ started at $z=\left( x,y\right) 
$ is Markov and (unique) solution of the RDE 
\begin{equation*}
\left\{ 
\begin{array}{l}
dZ^{z}=W\left( Z^{z}\right) d\mathbf{X}^{a,x}, \\ 
Z_{0}^{z}=\left( x,y\right) .%
\end{array}%
\right.
\end{equation*}%
where $W=\left( W_{1},...,W_{d}\right) $ with vector fields $W_{i}$ on $%
g^{2}\left( \mathbb{R}^{d}\right) \oplus \mathbb{R}^{e}$ given by%
\begin{equation*}
W_{i}\left( x,y\right) =\left( U_{i}\left( x\right) ,V_{i}\left( y\right)
\right) ,\,\,\,\left( x,y\right) \in g^{2}\left( \mathbb{R}^{d}\right)
\oplus \mathbb{R}^{e}\text{.}
\end{equation*}%
Recall that $U_{i}:g^{2}\left( \mathbb{R}^{d}\right) \rightarrow g^{2}\left( 
\mathbb{R}^{d}\right) $ are the vector fields defined in (\ref{DefOf_U_on_g2}%
) and, by the usual identification with first order differential operators,
the $U_{i}$ extend canonically to first order differential operators (and
hence vector fields) on $g^{2}\left( \mathbb{R}^{d}\right) \oplus \mathbb{R}%
^{e}$ which we denote for clarity with $\widetilde{U_{i}}$. We now describe
the infinitesimal behaviour of the associated semigroup $t\mapsto \mathbb{E}%
^{a,x}\left( f\left( Z_{t}^{.}\right) \right) $.

\begin{proposition}
\label{EZonCc}Let $f,g\in C_{c}^{\infty }\left( g^{2}\left( \mathbb{R}%
^{d}\right) \oplus \mathbb{R}^{e}\right) $. Then%
\begin{multline*}
\lim_{t\rightarrow 0}\left\langle \frac{\mathbb{E}^{a,.}\left( f\left(
Z_{t}^{.}\right) -f\left( .\right) \right) }{t},g\right\rangle _{L^{2}\left(
g^{2}\left( \mathbb{R}^{d}\right) \oplus \mathbb{R}^{e}\right) } \\
\left. =-\sum_{i,j=1}^{d}\int_{g^{2}\left( \mathbb{R}^{d}\right) \oplus 
\mathbb{R}^{e}}a^{ij}\left( x\right) W_{i}f\left( x,y\right) W_{j}^{\ast
}g\left( x,y\right) \,dxdy=:\mathcal{E}^{Z}\left( f,g\right) \right.
\end{multline*}%
where $W^{\ast }$ is the adjoint of $W$ with respect to Lebesue measure on $%
g^{2}\left( \mathbb{R}^{d}\right) \oplus \mathbb{R}^{e}$.
\end{proposition}

\begin{proof}
Let us fix $f,g$ $C_{c}^{\infty }\left( g^{2}\left( \mathbb{R}^{d}\right)
\oplus \mathbb{R}^{e}\right) .$ We want to apply Lemma \ref{EulerEstimates}
with unbounded vector fields $W$ (the unboundedness comes from the $U_{i}$)
and we need to localize our problem. Let $R>0$ such that $f$ and $g$ are $0$
outside $B\left( 0,R\right) ,$ and define compactly supported smooth vector
fields $U_{i}^{R}$ such that $U_{i}^{R}$ and $U_{i}$ agree on $B\left(
0,2R\right) $. Let $Z^{R}$ denote the solution of the RDE driven by $X$
along the vector fields $W_{i}^{R}=\left( U_{i}^{R},V_{i}\right) .$ Observe
first that $W_{i}^{R}f=W_{i}f$ by construction. Applying Lemma \ref%
{EulerEstimates}, we obtain 
\begin{multline*}
\lim_{t\rightarrow 0}\left\langle \frac{\mathbb{E}\left( f\left(
Z_{t}^{R,.}\right) -f\left( .\right) \right) }{t},g\right\rangle
_{L^{2}\left( g^{2}\left( \mathbb{R}^{d}\right) \oplus \mathbb{R}^{e}\right)
} \\
\left. 
\begin{array}{l}
=\sum_{i=1}^{d}\lim_{t\rightarrow 0}\left\langle W_{i}f\left( .\right) \frac{%
\mathbb{E}^{a,.}\left[ \pi _{1,i}\left( \mathbf{X}_{0,t}\right) \right] }{t}%
,g\right\rangle _{L^{2}\left( g^{2}\left( \mathbb{R}^{d}\right) \oplus 
\mathbb{R}^{e}\right) } \\ 
+\frac{1}{2}\sum_{i,j=1}^{d}\lim_{t\rightarrow 0}\left\langle
W_{i}W_{j}f\left( .\right) \frac{\mathbb{E}^{a,.}\left[ \pi _{1,i}\left( 
\mathbf{X}_{0,t}\right) \pi _{1,j}\left( \mathbf{X}_{0,t}\right) \right] }{t}%
,g\right\rangle _{L^{2}\left( g^{2}\left( \mathbb{R}^{d}\right) \oplus 
\mathbb{R}^{e}\right) } \\ 
+\frac{1}{2}\sum_{i,j=1}^{d}\lim_{t\rightarrow 0}\left\langle \left[
W_{i},W_{j}\right] f\left( .\right) \frac{\mathbb{E}^{a,.}\left[ \pi
_{2,i,j}\left( \mathbf{X}_{0,t}\right) \right] }{t},g\right\rangle
_{L^{2}\left( g^{2}\left( \mathbb{R}^{d}\right) \oplus \mathbb{R}^{e}\right)
}.%
\end{array}%
\right.
\end{multline*}%
As $Z_{t}^{R,.}$ and $Z$ differ only through the area of $\mathbf{X}^{a,.},$
using that uniformly over $x\in B\left( 0,R\right) $, the probability of $%
\mathbf{X}^{a,x}$ going outside $B\left( 0,2R\right) $ is bounded above by $%
C\exp \left( -CR^{2}\right) ,$ we easily see that 
\begin{equation*}
\lim_{R\rightarrow \infty }\lim_{t\rightarrow 0}\left\langle \frac{\mathbb{E}%
\left( f\left( Z_{t}^{R,.}\right) -f\left( .\right) \right) }{t}%
,g\right\rangle _{L^{2}}=\lim_{t\rightarrow 0}\left\langle \frac{\mathbb{E}%
\left( f\left( Z_{t}^{.}\right) -f\left( .\right) \right) }{t}%
,g\right\rangle _{L^{2}}.
\end{equation*}%
We then use lemma \ref{weakconvergencemoments} to obtain%
\begin{multline*}
\lim_{t\rightarrow 0}\left\langle \frac{\mathbb{E}\left( f\left(
Z_{t}^{.}\right) -f\left( .\right) \right) }{t},g\right\rangle _{L^{2}\left(
g^{2}\left( \mathbb{R}^{d}\right) \oplus \mathbb{R}^{e}\right) } \\
\left. 
\begin{array}{l}
=\sum_{i,j=1}^{d}\int_{\left( x,y\right) \in g^{2}\left( \mathbb{R}%
^{d}\right) \times \mathbb{R}^{e}}a^{ij}\left( x\right) \widetilde{U_{i}}%
\left[ gW_{j}f\right] \left( x,y\right) dxdy \\ 
-\sum_{i,j=1}^{d}\int_{\left( x,y\right) \in g^{2}\left( \mathbb{R}%
^{d}\right) \times \mathbb{R}^{e}}a^{ij}\left( x\right) g\left( x,y\right)
W_{i}W_{j}f\left( x,y\right) dxdy.%
\end{array}%
\right.
\end{multline*}%
The proof is finished if we can show 
\begin{multline*}
\int_{\left( x,y\right) \in g^{2}\left( \mathbb{R}^{d}\right) \times \mathbb{%
R}^{e}}a^{ij}\left( x\right) \left( \widetilde{U_{i}}\left[ gW_{j}f\right]
\left( x,y\right) -g\left( x,y\right) W_{i}W_{j}f\left( x,y\right) \right)
dxdy \\
=\int_{g^{2}\left( \mathbb{R}^{d}\right) \oplus \mathbb{R}^{e}}a^{ij}\left(
x\right) W_{i}f\left( x,y\right) W_{j}^{\ast }g\left( x,y\right) \,dxdy
\end{multline*}%
and to see this we may assume, by a simple limit argument, that $a\in \Xi
\left( \Lambda \right) $ is smooth. We have $a^{ij}\left( x\right) 
\widetilde{U_{i}}\left[ gW_{j}f\right] \left( x,y\right) $ equal to 
\begin{equation*}
\widetilde{U_{i}}\left[ a^{ij}\left( \pi _{g^{2}\left( \mathbb{R}^{d}\right)
}\left( .\right) \right) gW_{j}f\right] \left( x,y\right) -\widetilde{U_{i}}%
\left[ a^{ij}\left( \pi _{g^{2}\left( \mathbb{R}^{d}\right) }\left( .\right)
\right) \right] \left( gW_{j}f\right) \left( x,y\right) ,
\end{equation*}%
and $a^{ij}\left( x\right) g\left( x,y\right) W_{i}W_{j}f\left( x,y\right) $
equal to 
\begin{equation*}
g\left( x,y\right) W_{i}\left[ a^{ij}\left( \pi _{g^{2}\left( \mathbb{R}%
^{d}\right) }\left( .\right) \right) W_{j}f\right] \left( x,y\right) -W_{i}%
\left[ a^{ij}\left( \pi _{g^{2}\left( \mathbb{R}^{d}\right) }\left( .\right)
\right) \right] \left( gW_{j}f\right) \left( x,y\right) .
\end{equation*}%
But by construction of $W_{i}$ we have $W_{i}\left[ a^{ij}\left( \pi
_{g^{2}\left( \mathbb{R}^{d}\right) }\left( .\right) \right) \right] =%
\widetilde{U_{i}}\left[ a^{ij}\left( \pi _{g^{2}\left( \mathbb{R}^{d}\right)
}\left( .\right) \right) \right] $. Moreover, by integration by parts,%
\begin{equation*}
\int_{\left( x,y\right) \in g^{2}\left( \mathbb{R}^{d}\right) \times \mathbb{%
R}^{e}}\widetilde{U_{i}}\left[ a^{ij}\left( \pi _{g^{2}\left( \mathbb{R}%
^{d}\right) }\left( .\right) \right) gW_{j}f\right] \left( x,y\right) dxdy=0,
\end{equation*}%
and we see that 
\begin{multline*}
\int_{\left( x,y\right) \in g^{2}\left( \mathbb{R}^{d}\right) \times \mathbb{%
R}^{e}}a^{ij}\left( x\right) \left( \widetilde{U_{i}}\left[ gW_{j}f\right]
\left( x,y\right) -g\left( x,y\right) W_{i}W_{j}f\left( x,y\right) \right)
dxdy \\
\left. 
\begin{array}{l}
=-\int_{\left( x,y\right) \in g^{2}\left( \mathbb{R}^{d}\right) \times 
\mathbb{R}^{e}}g\left( x,y\right) W_{i}\left[ a^{ij}\left( \pi _{g^{2}\left( 
\mathbb{R}^{d}\right) }\left( .\right) \right) W_{j}f\right] \left(
x,y\right) dxdy \\ 
=-\int_{\left( x,y\right) \in g^{2}\left( \mathbb{R}^{d}\right) \times 
\mathbb{R}^{e}}a^{ij}\left( x\right) W_{i}^{\ast }g\left( x,y\right)
W_{j}f\left( x,y\right) dxdy,\text{ }%
\end{array}%
\right.
\end{multline*}%
by definition of $W_{i}^{\ast }.$
\end{proof}

\begin{remark}
The reader might want to check that when $a\left( x\right) $ is smooth and
depends only on the projection of $x$ onto $\mathbb{R}^{d},$ an application
of It\^{o}'s lemma leads to the same result. In particular, when $a=I$ the
process $Z$ solves a Stratonovich equation along vector fields $%
W=(W_{1},...,W_{d})$ with generator in H\"{o}rmander form%
\begin{equation*}
L^{Z}=\sum_{i=1}^{d}W_{i}^{2}
\end{equation*}%
and the associated form $\left( f,g\right) \mapsto $ $-\left\langle
L^{Z}f,g\right\rangle =-\sum_{i=1}^{d}\int W_{i}fW_{i}^{\ast }g$ agrees with
Proposition \ref{EZonCc}.
\end{remark}

\section{Large Deviations}

We fix $a\in \Xi \left( \Lambda \right) $. The law of $t\mapsto \mathbf{X}%
^{a;x}\left( \varepsilon t\right) $ where $\mathbf{X}^{a;x}$ is the $%
g^{2}\left( \mathbb{R}^{d}\right) $-valued process associated to the
Dirichlet form $\mathcal{E}^{a}$, started at $x$, can be viewed as Borel
measure on $C_{x}\left( \left[ 0,1\right] ,g^{2}\left( \mathbb{R}^{d}\right)
\right) \subset C\left( \left[ 0,1\right] ,g^{2}\left( \mathbb{R}^{d}\right)
\right) $, i.e. the space of continuous paths started at $x$, and is denoted
by $\mathbb{P}_{\varepsilon }^{a;x}$. As usual, we write $\mathbf{X=X}^{a;x}$
when no confusion is possible and in particular under $\mathbb{P}^{a;x}$
where $\mathbf{X}_{t}\left( \omega \right) =\omega \left( t\right) \equiv
\omega _{t}$. We shall see that a sample path large deviation principle
holds w.r.t. to uniform (and then homogenous H\"{o}lder!) topology on $%
C_{x}\left( \left[ 0,T\right] ,g^{2}\left( \mathbb{R}^{d}\right) \right) $.
Having properties (i)-(iii) of the of following propostion, the proof
follows essentially Varadhan \cite{VaII67}, see also \cite{BaKu00}, and we
outline the key steps for the reader's convenience.

\begin{proposition}
\label{Bass22}(i) $\left( g^{2}\left( \mathbb{R}^{d}\right) ,d^{a}\right) $
is a geodesic space.\newline
(ii) The Varadhan-Ram\'{\i}rez short time formula holds,%
\begin{equation}
\lim_{\varepsilon \rightarrow 0}4\varepsilon \log p^{a}\left( \varepsilon
,x,y\right) =-d^{a}\left( x,y\right) ^{2}.  \label{VRrecalled}
\end{equation}%
(iii) For $\alpha \in (0,1/2)$ there exist a constant $C_{\ref{Bass22}}=$ $%
C_{\ref{Bass22}}\left( \alpha ,\Lambda \right) $such that%
\begin{equation*}
\sup_{x\in g^{2}\left( \mathbb{R}^{d}\right) }\mathbb{P}^{a;x}\left(
\sup_{0\leq s<t\leq 1}\frac{d^{a}\left( \mathbf{X}_{s},\mathbf{X}_{t}\right) 
}{\left\vert t-s\right\vert ^{\alpha }}>r\right) \leq C_{\ref{Bass22}}\exp
\left( -\frac{r^{2}}{C_{\ref{Bass22}}}\right)
\end{equation*}%
and the same estimate holds with $d$ instead of $d^{a}$.
\end{proposition}

\begin{proof}
(i) was shown in Proposition \ref{da_is_geodesic}, (ii) was discussed in the
section on short time asymptotics and (iii) follows from Theorem \ref%
{FerniqueEstimates}.
\end{proof}

On $C\left( \left[ 0,1\right] ,g^{2}\left( \mathbb{R}^{d}\right) \right) $,
equipped with uniform topology, we define the \textit{energy or action
functional}%
\begin{equation}
I^{a}\left( \omega \right) =\lim \sup_{\left\vert D\right\vert \rightarrow
0}\sum_{t_{i}\in D}\frac{d^{a}\left( \omega _{t_{i}},\omega
_{t_{i-1}}\right) ^{2}}{t_{i}-t_{i-1}}\in \left[ 0,\infty \right] .
\label{GoodRateFunction}
\end{equation}%
We shall see shortly that $I^{a}$ is a good rate function in the sense that $%
\phi \mapsto I^{a}\left( \phi \right) $ is lower semicontinuous with compact
level sets.

\subsection{Upper Bound}

We first recall that $d^{a}$ is a geodesic distance, i.e. that for all $%
x,y\in g^{2}\left( \mathbb{R}^{d}\right) ,$ there exists a continuous path
joining $x$ to $y,$ of length $d^{a}\left( x,y\right) .$

\begin{proposition}
\label{Bass24} (i) On $C\left( \left[ 0,1\right] ,g^{2}\left( \mathbb{R}%
^{d}\right) \right) $ we have%
\begin{equation*}
\inf_{\omega :\,\,\omega \left( s\right) =y,\omega \left( t\right)
=z}I^{a}\left( \omega \right) =\frac{d^{a}\left( y,z\right) ^{2}}{t-s}
\end{equation*}%
and the infimum is attained by a $d^{a}$-geodesic path.\newline
(ii) More generally,%
\begin{equation*}
\inf_{\substack{ \omega \left( t_{i}\right) =x_{i}  \\ i=1,...,m}}%
I^{a}\left( \omega \right) =I^{a}\left( \omega ^{D}\right) =\sum_{i=1}^{m}%
\frac{d^{a}\left( x_{i},x_{i-1}\right) ^{2}}{t_{i}-t_{i-1}}
\end{equation*}%
where $\omega ^{D}$ is a piecewise $d^{a}$-geodesic path with $\omega
^{D}\left( t_{i}\right) =x_{i}$ for all $i=1,...,m.$\newline
(iii) In particular, 
\begin{equation}
d^{a}\left( \omega _{s},\omega _{t}\right) \leq I^{a}\left( \phi \right)
^{1/2}\left( t-s\right) ^{1/2}.  \label{Bass211}
\end{equation}
\end{proposition}

\begin{proof}
Straight-forward, see \cite{VaII67, BaKu00} for instance.
\end{proof}

\begin{lemma}
\label{Bass25}(i) The functional $I^{a}$ is a good rate-function.\newline
(ii) If $C$ is closed and $C_{\delta }\supset C$ denotes the $\delta $%
-neighbourhood of $C$ (indifferently defined via $d$ or $d^{a}$) then 
\begin{equation*}
\lim_{\delta \rightarrow 0}\inf_{\omega \in C_{\delta }}I^{a}\left( \omega
\right) =\inf_{\omega \in C}I^{a}\left( \omega \right) .
\end{equation*}
\end{lemma}

\begin{proof}
Using (\ref{Bass211}) and Arzela-Ascoli this is proved as in \cite{Va80}.
\end{proof}

\begin{lemma}
\label{Bass27} Let $D$ be a dissection of $[0,1]$ with $\#D$ points and
define the (continuous) evaluation map%
\begin{equation*}
\Pi _{D}\left( \omega \right) :=\left( \omega _{t}\right) _{t\in D}\in \left[
g^{2}\left( \mathbb{R}^{d}\right) \right] ^{\#D}.
\end{equation*}%
Let $C$ be a closed "cylindrical" set of form $\Pi _{D}^{-1}A$ with $A\in %
\left[ g^{2}\left( \mathbb{R}^{d}\right) \right] ^{\#D}$ closed. Then%
\begin{equation*}
\lim \sup_{\varepsilon \rightarrow 0}4\varepsilon \log \mathbb{P}%
_{\varepsilon }^{a,x}\left( C\right) \leq -\inf_{\omega \in C}I^{a}\left(
\omega \right) .
\end{equation*}
\end{lemma}

\begin{proof}
Using the short time formula (\ref{VRrecalled}) and Lemma \ref{Bass24} this
is proved in the same way as \cite[Lemma 3.1]{VaII67}.
\end{proof}

\begin{lemma}
\label{Bass28}For every $\delta >0$,%
\begin{equation*}
\underset{m\rightarrow \infty }{\lim \sup }\lim \sup_{\epsilon \rightarrow
0}\varepsilon \log \sup_{x\in g^{2}\left( \mathbb{R}^{d}\right) }\mathbb{P}%
_{\varepsilon }^{a;x}\left( \sup_{0\leq t\leq 1}d^{a}\left( \mathbf{X}_{t},%
\mathbf{X}_{t}^{D_{m}}\right) >\delta \right) =-\infty
\end{equation*}%
where $\mathbf{X}^{D_{m}}$ is the $d$-geodesic approximation connecting the
points $\left\{ \mathbf{X}_{t}:t\in D^{m}\right\} $ with $D^{m}=\left\{
j/m:j=0,...,m\right\} $.
\end{lemma}

\begin{proof}
For a fixed $t$ and $D=D_{m}$ let $t_{D}$ be the closest point in $D$ to the
left of $t.$ Noting that $\mathbf{X}_{t_{D}}=\mathbf{X}_{t_{D}}^{D}$ and
using Lipschitz equivalence of $d$ and $d^{a}$ we have 
\begin{eqnarray*}
d^{a}\left( \mathbf{X}_{t},\mathbf{X}_{t}^{D}\right) &\leq &d^{a}\left( 
\mathbf{X}_{t},\mathbf{X}_{t_{D}}\right) +d^{a}\left( \mathbf{X}_{t_{D}},%
\mathbf{X}_{t_{D}}^{D}\right) +d^{a}\left( \mathbf{X}_{t}^{D},\mathbf{X}%
_{t_{D}}^{D}\right) \\
&\leq &C_{\ref{Bass28}}^{1}\left( d\left( \mathbf{X}_{t},\mathbf{X}%
_{t_{D}}\right) +d\left( \mathbf{X}_{t}^{D},\mathbf{X}_{t_{D}}^{D}\right)
\right)
\end{eqnarray*}%
We know from the earlier section on strong geodesic approximation that%
\begin{equation*}
\sup_{D}\left\Vert \mathbf{X}^{D}\right\Vert _{\alpha \text{-H\"{o}l;}\left[
0,1\right] }\leq 3\left\Vert \mathbf{X}\right\Vert _{\alpha \text{-H\"{o}l;}%
\left[ 0,1\right] }
\end{equation*}%
and it follows that%
\begin{equation*}
\sup_{0\leq t\leq 1}d^{a}\left( \mathbf{X}_{t},\mathbf{X}_{t}^{D}\right)
\leq 4\left\Vert \mathbf{X}\right\Vert _{\alpha \text{-H\"{o}l;}\left[ 0,1%
\right] }\times \left\vert D\right\vert ^{\alpha }
\end{equation*}%
where $\left\vert D\right\vert $ denotes the mesh of $D$ as usual. By a
simple scaling argument (section \ref{scalingDivForm}) and Proposition \ref%
{Bass22}, (iii) we see that%
\begin{equation*}
\mathbb{P}_{\varepsilon }^{a;x}\left( 4\left\Vert \mathbf{X}\right\Vert
_{\alpha \text{-H\"{o}l;}\left[ 0,1\right] }>\delta m^{\alpha }\right) \leq
C_{\ref{Bass22}}\exp \left( -\frac{1}{C_{\ref{Bass22}}}\frac{\delta
^{2}m^{2\alpha }}{\varepsilon }\right)
\end{equation*}%
and, noting that $C_{\ref{Bass22}}$ does not depend on $x$, 
\begin{eqnarray*}
\sup_{x}\mathbb{P}_{\varepsilon }^{a;x}\left( \sup_{0\leq t\leq
1}d^{a}\left( \mathbf{X}_{t},\mathbf{X}_{t}^{D_{m}}\right) >\delta \right)
&\leq &\sup_{x}\mathbb{P}_{\varepsilon }^{a;x}\left( 4\left\Vert \mathbf{X}%
\right\Vert _{\alpha \text{-H\"{o}l;}\left[ 0,1\right] }>\delta m^{\alpha
}\right) \\
&\leq &C_{\ref{Bass22}}\exp \left( -\frac{1}{C_{\ref{Bass22}}}\frac{\delta
^{2}m^{2\alpha }}{\varepsilon }\right) .
\end{eqnarray*}%
It readily follows that%
\begin{equation*}
\lim \sup_{m\rightarrow \infty }\lim \sup_{\varepsilon \rightarrow
0}\varepsilon \log \sup_{x}\mathbb{P}_{\varepsilon }^{a;x}\left( \sup_{0\leq
t\leq 1}d^{a}\left( \mathbf{X}_{t},\mathbf{X}_{t}^{D_{m}}\right) >\delta
\right) =-\infty
\end{equation*}%
as claimed.
\end{proof}

\begin{theorem}
For any measurable $A\subset C_{x}\left( \left[ 0,1\right] ,g^{2}\left( 
\mathbb{R}^{d}\right) \right) $%
\begin{equation*}
\lim \sup_{\varepsilon \rightarrow 0}4\varepsilon \log \mathbb{P}%
_{\varepsilon }^{a,x}\left( A\right) \leq -\inf_{\omega \in \bar{A}%
}I^{a}\left( \omega \right) .
\end{equation*}%
where $\bar{A}$ is the closure of $A$ w.r.t. to the uniform topology on path
space.
\end{theorem}

\begin{proof}
It suffices to consider $A$ closed. We write $A_{\delta }\supset A$ for the $%
\delta $-neighbourhood of $A$ (indifferently defined via $d$ or $d^{a}$) and
set%
\begin{equation*}
I^{\delta ,a}\left( \omega \right) :=\inf_{\tilde{\omega}:\sup_{t\in \left[
0,1\right] }d^{a}\left( \omega _{t},\tilde{\omega}_{t}\right) <\delta
}I^{a}\left( \tilde{\omega}\right) \text{ \ and \ }T_{\delta }:=\inf_{\omega
\in A_{\delta }}I^{a}\left( \omega \right) .\text{ }
\end{equation*}%
If $\omega \in A$ then $I^{\delta ,a}\left( \omega \right) \geq T_{\delta }$
and therefore, $D^{m}$ being defined as above,%
\begin{eqnarray*}
\mathbb{P}_{\varepsilon }^{a,x}\left( A\right) &\leq &\mathbb{P}%
_{\varepsilon }^{a,x}\left( \omega :I^{\delta ,a}\left( \mathbf{\omega }%
\right) \geq T_{\delta }\right) \\
&\leq &\mathbb{P}_{\varepsilon }^{a,x}\left[ \sup_{t}d^{a}\left( \omega
_{t},\omega _{t}^{D_{m}}\right) \geq \delta \right] +\mathbb{P}_{\varepsilon
}^{a,x}\left[ I^{a}\left( \omega ^{D_{m}}\right) \geq T_{\delta }\right] .
\end{eqnarray*}%
Noting that lemma \ref{Bass28} states precisely that%
\begin{equation*}
\lim \sup_{m\rightarrow \infty }\lim \sup_{\varepsilon \rightarrow
0}\varepsilon \log \sup_{x}\mathbb{P}_{\varepsilon }^{a,x}\left( \sup_{0\leq
t\leq 1}d^{a}\left( \omega _{t},\omega _{t}^{D_{m}}\right) >\delta \right)
=-\infty .
\end{equation*}%
and that, by Proposition \ref{Bass24}, (ii), the set $\left\{ \omega
:I^{a}\left( \omega ^{D_{m}}\right) \geq T_{\delta }\right\} $ is equal to%
\begin{equation*}
C^{m}:=\left\{ \omega :\sum_{i=1}^{m}\frac{d^{a}\left( \omega
_{t_{i}},\omega _{t_{i-1}}\right) ^{2}}{t_{i}-t_{i-1}}\geq T_{\delta
}\right\}
\end{equation*}%
we see from Lemma \ref{Bass27} that for any $m$,%
\begin{equation*}
\lim \sup_{\varepsilon \rightarrow 0}4\varepsilon \log \mathbb{P}%
_{\varepsilon }^{a,x}\left[ C^{m}\right] \leq -\inf_{\mathbf{\omega }\in
C^{m}}I^{a}\left( \omega \right) \leq -T_{\delta }.
\end{equation*}%
By Lemma \ref{Bass25}, $\lim_{\delta \rightarrow 0}T_{\delta }=\inf_{\mathbf{%
\omega }\in A}I^{a}\left( \omega \right) $ and combining all these results
yield the upper LDP\ bound.
\end{proof}

\subsection{Lower bound}

\begin{lemma}
For every $\omega \in C_{x}\left( \left[ 0,1\right] ,g^{2}\left( \mathbb{R}%
^{d}\right) \right) $ and every $\delta >0$,%
\begin{equation*}
\lim \inf_{\varepsilon \rightarrow 0}4\varepsilon \log \mathbb{P}%
_{\varepsilon }^{a,x}\left( B_{\delta }\left( \omega \right) \right) \geq
-I^{a}\left( \omega \right)
\end{equation*}%
where%
\begin{equation*}
B_{\delta }\left( \omega \right) =\left\{ \tilde{\omega}\in C_{x}\left( %
\left[ 0,1\right] ,g^{2}\left( \mathbb{R}^{d}\right) \right) :\sup_{t\in %
\left[ 0,1\right] }d^{a}\left( \omega _{t},\tilde{\omega}_{t}\right) <\delta
\right\} \text{.}
\end{equation*}
\end{lemma}

\begin{proof}
Using the short time formula (\ref{VRrecalled}), Lemma \ref{Bass24} and the
upper LDP\ this is proved as \cite[Lemma 3.4]{VaII67}.
\end{proof}

\begin{corollary}
For any measurable $A\subset C_{x}\left( \left[ 0,1\right] ,g^{2}\left( 
\mathbb{R}^{d}\right) \right) $%
\begin{equation*}
-\inf_{\omega \in A^{\circ }}I^{a}\left( \omega \right) \leq \lim
\inf_{\varepsilon \rightarrow 0}4\varepsilon \log \mathbb{P}_{\varepsilon
}^{a,x}\left( A\right)
\end{equation*}%
where $A^{\circ }$ is the interior of $A$ w.r.t. to the uniform topology on
path space.
\end{corollary}

\begin{proof}
W.l.o.g. assume that $A$ is open. Take any $\omega \in A$ and $\delta >0$
small enough such that $V=B_{\delta }\left( \omega \right) \subset A$. From
the last lemma it then follows that%
\begin{equation*}
\lim \inf_{\varepsilon \rightarrow 0}4\varepsilon \log \mathbb{P}%
_{\varepsilon }^{a,x}\left( A\right) \geq \lim \inf_{\varepsilon \rightarrow
0}4\varepsilon \log \mathbb{P}_{\varepsilon }^{a,x}\left( V\right) \geq
-I^{a}\left( \omega \right) .
\end{equation*}%
As this is true for all $f\in A$ we have the result.
\end{proof}

\subsection{LDP\ in H\"{o}lder topology \&\ Freidlin Wentzell}

The above estimates are summarized in

\begin{theorem}
Let $\mathbb{P}_{\varepsilon }^{a;x}$ be the law of $t\mapsto \mathbf{X}%
^{a;x}\left( \varepsilon t\right) $ where $\mathbf{X}^{a;x}$ is the $%
g^{2}\left( \mathbb{R}^{d}\right) $-valued process associated to the
Dirichlet form $\mathcal{E}^{a}$. Then $\left( \mathbb{P}_{\varepsilon
}^{a;x}\right) _{\varepsilon >0}$ satisfies a large deviation principle in
uniform topology on $C_{x}\left( \left[ 0,T\right] ,g^{2}\left( \mathbb{R}%
^{d}\right) \right) $ with good rate function $I^{a}$ defined in equation (%
\ref{GoodRateFunction}).
\end{theorem}

It would be easy to deduce from this result a functional form of Strassen's
Law of Iterated Logarithm holds, see \cite{DeuSt89}, but we shall not pursue
this here.

\begin{corollary}
Fix $\alpha \in \lbrack 0,1/2)$. Then $\left( \mathbb{P}_{\varepsilon
}^{a;x}\right) _{\varepsilon >0}$ satisfies a large deviation principle in $%
\alpha $-H\"{o}lder topology on $C^{\alpha \text{-H\"{o}lder}}\left( \left[
0,T\right] ,g^{2}\left( \mathbb{R}^{d}\right) \right) $ with good rate
function $I^{a}$.
\end{corollary}

\begin{proof}
The random variable $\left\Vert \mathbf{X}^{a;x}\right\Vert _{\alpha \text{-H%
\"{o}lder}}$ has a Gaussian tail for all $\alpha <1/2$. By the inverse
contraction principle \cite{DeZe93} we see that the large deviation
principle in uniform topology can be strenghtened to $\alpha $-H\"{o}lder
topology.
\end{proof}

From the contraction principle and Lyons' universal limit theorem \cite%
{lyons-98} we obtain

\begin{corollary}[Freidlin-Wentzell]
Let $Y_{\varepsilon }=\pi \left( 0,y_{0};\mathbf{X}_{\varepsilon
}^{a,x}\right) $ denote the $\mathbb{R}^{e}$-valued (random) RDE solution
driven by $\mathbf{X}_{\varepsilon }^{a,x}=\mathbf{X}^{a,x}\left(
\varepsilon \cdot \right) $ along fixed $\mathrm{Lip}^{2+\varepsilon }$
vector fields $V_{1},...,V_{d}$ on $\mathbb{R}^{e}$ and started at time $0$
from $y_{0}$ fixed (i.e. $\pi $ is the It\^{o} map). Let $\mathbb{Q}%
_{\varepsilon }$ denote the law of $Y_{\varepsilon }$. Then $\left( \mathbb{Q%
}_{\varepsilon }:\varepsilon >0\right) $ satisfies a large deviation
principle in $\alpha $-H\"{o}lder topology, $\alpha \in \lbrack 0,1/2)$,
with good rate function 
\begin{equation*}
J^{a}(y\mathbf{)}=\inf \left\{ I^{a}(\omega ):\omega \in C_{x}\left( \left[
0,1\right] ,g^{2}\left( \mathbb{R}^{d}\right) \right) \text{ and }y=\pi
\left( 0,y_{0};\omega \right) \right\} \text{.}
\end{equation*}
\end{corollary}

\section{Support Theorems}

To prove an extension of the Stroock-Varadhan support theorem \cite{StVa72,
ikeda-watanabe-89} (\cite{BeGrLe94} and \cite{MiSaSo94} for H\"{o}lder
topology) to RDEs driven by the "Markovian" rough paths $\mathbf{X}^{a;x}$,
it would be enough to show to for fixed $a\in \Xi \left( \Lambda \right) $
and some $\alpha \in \left( 1/3,1/2\right) $,%
\begin{equation*}
\text{supp}\left( \mathbb{P}^{a;x}\right) =x\ast C_{0}^{0,\alpha }\left( %
\left[ 0,1\right] ,g^{2}\left( \mathbb{R}^{d}\right) \right) .
\end{equation*}%
The $\subset $ direction is obvious (from section \ref{SectionWongZakai} )
but equality remains an open (and challenging) problem. Nonetheless, we are
able to prove the desired extension of the Stroock-Varadhan support theorem.
First, by shifting the argument of $a$ we can and will assume $x=0$. If we
can show that for fixed $a\in \Xi \left( \Lambda \right) $, some $n\geq 2$
and $\alpha \in \left( \frac{1}{n+1},\frac{1}{n}\right) ,$ 
\begin{equation*}
\text{supp}\left( \left( S_{n}\right) _{\ast }\mathbb{P}^{a;0}\right)
=C_{0}^{0,\alpha }\left( \left[ 0,1\right] ,g^{n}\left( \mathbb{R}%
^{d}\right) \right)
\end{equation*}%
where $S_{n}:C_{0}^{0,\gamma }\left( \left[ 0,1\right] ,g^{2}\left( \mathbb{R%
}^{d}\right) \right) \rightarrow C_{0}^{0,\gamma }\left( \left[ 0,1\right]
,g^{n}\left( \mathbb{R}^{d}\right) \right) $ is the continuous Young-Lyons
lift, $\gamma \in \left( 1/3,1/2\right) $, the extended Stroock-Varadhan
support theorem (in H\"{o}lder topology of exponent less than $1/n$ and
hence in uniform topology) is a consequence of basic consistency properties
of RDE solutions and the fundamental continuity result of rough path theory.
Validity of the Stroock-Varadhan support theorem for differential equations
driven by $\mathbf{X}^{a,x}$ in the rough paths sense was conjectured, via
conditional statements, by T. Lyons in \cite{lyons-04}. \ 

\subsection{Support in Uniform Topology\label{SecSuppUniTop}}

Let $h,$ $x\in C_{0}^{1}\left( \left[ 0,1\right] ,\mathbb{R}^{d}\right) $.
Every such $x$ can be lifted to $S\left( x\right) \in C_{0}^{1\text{-var}%
}\left( \left[ 0,1\right] ,g^{2}\left( \mathbb{R}^{d}\right) \right) $ via
iterated integration. Similarly, one can lift $x+h$, the translation of $x$
in direction $h$, to a $g^{2}\left( \mathbb{R}^{d}\right) $-valued path $%
S\left( x+h\right) $. Provided $\alpha \in (1/3,1/2],$this operation extends
to a continuous translation operator $T_{h}$,%
\begin{equation*}
\mathbf{x}\in C^{\alpha \text{-H\"{o}l}}\left( \left[ 0,1\right]
,g^{2}\left( \mathbb{R}^{d}\right) \right) \mapsto T_{h}\mathbf{x}\in
C^{\alpha \text{-H\"{o}l}}\left( \left[ 0,1\right] ,g^{2}\left( \mathbb{R}%
^{d}\right) \right) .
\end{equation*}%
We refer to \cite{lyons-qian-02} for details. We note that for $h\in
C_{0}^{1}\left( \left[ 0,1\right] ,\mathbb{R}^{d}\right) $ fixed and a
sequence $\left( \mathbf{x}_{k}\right) $, 
\begin{equation*}
\left\Vert T_{h}\mathbf{x}_{k}\right\Vert _{\alpha \text{-H\"{o}l}%
}\rightarrow 0\text{ as }k\rightarrow \infty \text{ iff }d_{\alpha \text{-H%
\"{o}l}}\left( S\left( h\right) ,\mathbf{x}_{k}\right) \rightarrow 0\text{
as }k\rightarrow \infty \text{.}
\end{equation*}%
Assuming that $a\left( x\right) $ only depends on $\pi _{1}\left( x\right) $%
, with abuse of notation $a=a\left( \pi _{1}\left( \cdot \right) \right) $,
we have that $X+h$ is Markov with (formal) generator%
\begin{equation*}
\sum_{i,j}\partial _{i}\left( a^{ij}\left( \cdot -h_{t}\right) \partial
_{j}\right) +\sum_{k}\dot{h}_{t}^{k}\partial _{k}
\end{equation*}%
and $T_{h}\mathbf{X}$ is a Markov with (formal) generator%
\begin{equation*}
\sum_{i,j}U_{i}\left( a^{ij}\left( \pi _{1}\left( \cdot \right)
-h_{t}\right) U_{j}\right) +\sum_{k}\dot{h}_{t}^{k}U_{k}
\end{equation*}%
where $U_{1},...,U_{d}$ are the generating left-invariant vector fields on $%
g^{2}\left( \mathbb{R}^{d}\right) $.

\begin{proposition}
\label{subLower}Let $h\in C_{0}^{1}\left( \left[ 0,1\right] ,\mathbb{R}%
^{d}\right) $. There exists a constant $C_{\ref{subLower}}$ depending only
on $\Lambda $ and $\left\vert \dot{h}\right\vert _{\infty ;\left[ 0,1\right]
}$ such that for all $\varepsilon \in (0,1],$%
\begin{equation*}
\mathbb{P}^{a;0}\left( \left\Vert T_{h}\left( \mathbf{X}\right) \right\Vert
_{\infty ,\left[ 0,1\right] }<\varepsilon \right) \geq \frac{1}{C_{\ref%
{subLower}}}\exp \left( -\frac{C_{\ref{subLower}}}{\varepsilon ^{2}}\right) .
\end{equation*}%
As a consequence, the support of $\mathbb{P}^{a,0}$ equals the closure of $%
S\left( C_{0}^{1}\left( \left[ 0,1\right] ,\mathbb{R}^{d}\right) \right) $,
with respect to uniform topology on $C_{0}\left( \left[ 0,1\right]
,g^{2}\left( \mathbb{R}^{d}\right) \right) .$
\end{proposition}

\begin{proof}
We first consider $h=0$. Let $n$ be the smallest integer such that $%
n^{-1/2}\leq \varepsilon /2.$ Set $y_{0}=0\in g^{2}\left( \mathbb{R}%
^{d}\right) $. Clearly, $\mathbb{P}^{a;y_{0}}\left( \left\Vert \mathbf{X}%
\right\Vert _{\infty ,\left[ 0,1\right] }<\varepsilon \right) $ is greater
or equal than%
\begin{equation*}
q_{\varepsilon ,n}:=\mathbb{P}^{a;y_{0}}\left( \forall i\in \left\{ 1,\cdots
,n\right\} :\text{ }\left\Vert X_{\frac{i}{n}}\right\Vert <n^{-1/2}\text{
and }X_{t}\in B\left( 0,\varepsilon \right) \text{ }\forall t\in \left[ 
\frac{i}{n},\frac{i+1}{n}\right] \right) .
\end{equation*}%
Hence, letting $p^{B\left( 0,\varepsilon \right) }$denote the Dirichlet heat
kernel for $\mathbf{X=X}^{a;0}$, the Markov property implies%
\begin{equation*}
q_{\varepsilon ,n}=\int_{B\left( 0,n^{-1/2}\right) }\cdots \int_{B\left(
0,n^{-1/2}\text{ }\right) }p_{B\left( 0,\varepsilon \right) }^{a}\left(
1/n,y_{0},y_{1}\right) \cdots p_{B\left( 0,\varepsilon \right) }^{a}\left(
1/n,y_{n-1},y_{n}\right) dy_{1}\cdots dy_{n}.
\end{equation*}%
We join the points $y_{i}$ and $y_{i+1}$ by the curve $\gamma _{i}$, which
is the concatenation of a geodesic curve joining $y_{i}$ to $0$ and a
geodesic curve between $0$ to $y_{i+1}.$ In particular, the length of $%
\gamma _{i}$ is bounded by $2n^{-1/2},$ and $\gamma _{i}$ remains in the
ball $B\left( 0,n^{-1/2}\right) \subset B\left( 0,\varepsilon /2\right) .$
Hence 
\begin{equation*}
R_{i}\equiv d^{a}\left( \gamma _{i},g^{2}\left( \mathbb{R}^{d}\right)
/B\left( 0,\varepsilon \right) \right) \geq \varepsilon /2\geq n^{-1/2}
\end{equation*}%
and we can apply the lower heat kernel bounds for the killed process with~$%
t=1/n$ and $\delta =\min \left( R_{i}^{2},t\right) =1/n$ to obtain 
\begin{eqnarray*}
p_{B\left( 0,\varepsilon \right) }^{a}\left( 1/n,y_{i},y_{i+1}\right) &\geq &%
\frac{n^{d^{2}/2}}{C_{\ref{subLower}}^{1}}\exp \left( -C_{\ref{subLower}%
}^{1}nd\left( y_{i},y_{i+1}\right) ^{2}\right) \exp \left( -\frac{C_{\ref%
{subLower}}^{1}}{nR_{i}^{2}}\right) \\
&\geq &\frac{n^{d^{2}/2}}{C_{\ref{subLower}}^{1}}\exp \left( -5C_{\ref%
{subLower}}^{1}\right)
\end{eqnarray*}%
where we used $d\left( y_{i},y_{i+1}\right) \leq 2n^{-1/2}$ and $R_{i}\geq
n^{-1/2}$. Since $m\left( B_{r}\left( 0\right) \right) \simeq r^{N}$ with
doubling constant $N=d^{2}$ we find 
\begin{eqnarray*}
q_{\varepsilon ,n} &\geq &\prod_{i=1}^{n}\left\{ \frac{1}{C_{\ref{subLower}%
}^{1}}\exp \left( -5C_{\ref{subLower}}^{1}\right) \times \frac{m\left(
B\left( 0,n^{-1/2}\right) \right) }{\left( n^{-1/2}\right) ^{d^{2}}}\right\}
= \\
&\geq &\left\{ \exp \left( -C_{\ref{subLower}}^{2}\right) \right\} ^{n}\geq
\exp \left( -\frac{C_{\ref{subLower}}^{3}}{\varepsilon ^{2}}\right) .
\end{eqnarray*}%
For $h\neq 0$ we note that the process $T_{h}\left( \mathbf{X}\right) $ is
described by a non-symmetric, time dependent Dirichlet form as in \cite%
{SturmII}, for instance. More precisely, the $\mathbb{R}^{d}$-valued process 
$t\mapsto h\left( t\right) +\mathbf{X}^{1}\left( t\right) $ is desribed by
the form%
\begin{equation*}
\left( f,g\right) \mapsto \int_{\mathbb{R}^{d}}\left[ a^{ij}\left( \cdot
,t\right) \partial _{i}f\partial _{j}g+gb^{i}\left( \cdot ,t\right) \partial
_{i}f\right] dx
\end{equation*}%
and the bilinear form for $T_{h}\left( \mathbf{X}\right) $ its the natural
lift obtained by replacing $\partial _{i}$ by $U_{i}$ for $i=1,...,d$,%
\begin{equation*}
\left( f,g\right) \mapsto \int_{g^{2}\left( \mathbb{R}^{d}\right) }\left[
a^{ij}\left( \pi _{1}\left( \cdot \right) ,t\right)
U_{i}fU_{j}g+gb^{i}\left( \cdot ,t\right) U_{i}f\right] dm
\end{equation*}%
Such lower order perturbations and time-dependence have been discussed in 
\cite{St88, SaSt91, ElRo00}. In particular, there are lower heat kernel
bounds for the killed process which allow the above proof to go through.
\end{proof}

\subsection{Support in H\"{o}lder Topology:\ A Conditional Result}

Motivated by \cite{friz-lyons-stroock-06} we first study the probability
that $\mathbf{X}^{a;x}$ stays in bounded open domain $D\subset g^{2}\left( 
\mathbb{R}^{d}\right) $ for long times.

\begin{proposition}
\label{SpectralProp}Let $D$ be an open domain in $g^{2}\left( \mathbb{R}%
^{d}\right) $ with finite volume, no regularity assumptions are made about $%
\partial D$. Let $a\in \Xi \left( \Lambda \right) $ and $\mathbf{X}^{a}$ be
the process associated to $\mathcal{E}^{a}$ started at $x\in g^{2}\left( 
\mathbb{R}^{d}\right) $ and assume $x\in D$. Then there exist positive
constants $K_{1}=K_{1}\left( x,D,\Lambda \right) $ and $K_{2}=K_{2}\left(
D,\Lambda \right) $ so that for all $t\geq 0$%
\begin{equation*}
K_{1}e^{-\lambda t}\leq \mathbb{P}\left[ \mathbf{X}_{s}^{a,x}\in D\,\forall
s:0\leq s\leq t\right] \leq K_{2}e^{-\lambda t}
\end{equation*}%
where $\lambda \equiv \lambda _{1}^{a}>0$ is the simple and first Dirichlet
eigenvalue of $-L^{a}$ on the domain $D.$ Moreover,%
\begin{equation*}
\forall a\in \Xi \left( \Lambda \right) :0<\lambda _{\min }\leq \lambda
_{1}^{a}\leq \,\,\,\lambda _{\max }<\infty
\end{equation*}%
where $\lambda _{\min },\,\lambda _{\max }$ depend only on $\Lambda $ and $D$%
.
\end{proposition}

\begin{remark}
The proof will show that $K_{1}\sim \psi _{1}^{a}\left( x\right) $. Noting
that $\psi _{1}^{a}\left( x\right) e^{-\lambda _{1}^{a}t}$ solves the same
PDE\ as $u^{a}\left( t,x\right) $, the above can be regarded as a "partial"
parabolic boundary Harnack statement.
\end{remark}

\begin{proof}
If $p_{D}^{a}$ denotes the Dirichlet heat kernel for $D$ we can write%
\begin{equation*}
u^{a}\left( t,x\right) :=\mathbb{P}^{x}\left[ \mathbf{X}_{s}^{a}\in
D\,\forall s:0\leq s\leq t\right] =\int_{D}p_{D}^{a}\left( t,x,y\right) dy.
\end{equation*}%
Recall \cite{Fu94} that $p_{D}^{a}$ is the kernel for a semigroup $%
P_{D}^{a}:L^{2}\left( D\right) \rightarrow L^{2}\left( D\right) $ which
corresponds to the Dirichlet form $\left( \mathcal{E}^{a},\mathcal{F}%
_{D}\right) $ whose domain $\mathcal{F}_{D}$ consists of all $f\in \mathcal{F%
}\equiv D\left( \mathcal{E}^{a}\right) $ with quasi continuous modifications
equal to $0$ q.e. on $D^{c}$. The infinitesimal generator of $P_{D}^{a}$,
denoted by $L_{D}^{a}$, is a self-adjoint, densely defined operator with
spectrum $\sigma \left( -L_{D}^{a}\right) \subset \lbrack 0,\infty )$. We
now use an ultracontractivity argument to show that $\sigma \left(
-L_{D}^{a}\right) $ is discrete. To this end, we note that the upper bound
on $p^{a}$ plainly implies $\left\vert p_{D}^{a}\left( t,\cdot ,\cdot
\right) \right\vert _{\infty }=O(t^{-d^{2}/2})$. Since $\left\vert
D\right\vert <\infty $ if follows that $\left\Vert P_{D}^{a}\left( t\right)
\right\Vert _{L^{1}\rightarrow L^{\infty }}<\infty $ which is, by
definition, ultracontractivity of the semigroup $P_{D}^{a}$. It is now a
standard consequence \cite[Thm 2.1.4]{Dav89} that $\sigma \left(
-L_{D}^{a}\right) =\left\{ \lambda _{1}^{a},\lambda _{2}^{a},...\right\}
\subset \lbrack 0,\infty )$, listed in non-decreasing order. Moreover, it is
clear that $\lambda _{1}^{a}\neq 0$; indeed the kernel estimates are plenty
to see that $\left\Vert P_{D}^{a}\left( t\right) \right\Vert
_{L^{2}\rightarrow L^{2}}\rightarrow 0$ as $t\rightarrow \infty $ which
contradicts the the existence of non-zero $f\in L^{2}\left( D\right) $ so
that $P_{D}^{a}\left( t\right) f=f$ for all $t\geq 0$. Let us note that%
\begin{eqnarray*}
\lambda _{1}^{a} &=&\inf \sigma \left( H\right) \\
&=&\inf \left\{ \mathcal{E}^{a}\left( f,f\right) :f\in \mathcal{F}_{D}\text{
with }\left\vert f\right\vert _{L^{2}\left( D\right) }=1\right\} \text{ \ \
\ (by Rayleigh-Ritz)} \\
&=&\inf \left\{ \int_{D}\Gamma ^{a}\left( f,f\right) dm:f\in \mathcal{F}_{D}%
\text{ with }\left\vert f\right\vert _{L^{2}\left( D\right) }=1\right\}
\end{eqnarray*}%
and since $\Gamma ^{a}\left( f,f\right) /\Gamma ^{I}\left( f,f\right) \in %
\left[ \Lambda ^{-1},\Lambda \right] $ for $f$ $\neq 0$ it follows that $%
\lambda _{1}^{a}\in \left[ \lambda _{\min },\lambda _{\max }\right] $ for
all $a\in \Xi \left( \Lambda \right) $ where we set%
\begin{equation}
\lambda _{\min }=\Lambda ^{-1}\lambda _{1}^{I},\,\,\,\lambda _{\max
}=\Lambda \lambda _{1}^{I}\text{.}  \label{lambdaMinMax}
\end{equation}%
From \cite[Thm 1.4.3]{Dav89} the lower heat kernel estimates for the killed
process imply irreducibility of the semigroup $P_{D}^{a}$, hence simplicity
of the first eigenvalue $\lambda $, and there is an a.s. strictly positive
eigenfunction to $\lambda \equiv \lambda _{1}^{a}$, say $\psi \equiv \psi
_{1}^{a}$, and by De Giorgi-Moser-Nash regularity we may assume that $\psi $
is H\"{o}lder continuous and strictly positive away from the boundary\ (this
follows also from Harnack's inequality). We also can (and will) assume that $%
\left\Vert \psi \right\Vert _{L^{2}\left( D\right) }=1$.\newline
\underline{Lower bound:} Noting that $v\left( t,x\right) =e^{-\lambda t}\psi
\left( x\right) $ is a weak solution of $\partial _{t}v=L_{D}^{a}v$ with $%
v\left( 0,\cdot \right) =\psi $ we have 
\begin{equation*}
v\left( t,x\right) =\int_{D}p_{D}^{a}\left( t,x,y\right) \psi \left(
y\right) dy,
\end{equation*}%
at first for a.e. $x$ but by using a H\"{o}lder regular version of $%
p_{D}^{a} $ the above holds for all $x\in D$. It follows that%
\begin{eqnarray*}
0 &<&\psi \left( x\right) \\
&=&e^{\lambda t}\int_{D}p_{D}^{a}\left( t,x,y\right) \psi \left( y\right) dy
\\
&\leq &e^{\lambda \left( t+1\right) }\int_{D}p_{D}^{a}\left( t,x,y\right)
\int_{D}p_{D}^{a}\left( 1,y,z\right) \psi \left( z\right) dzdy \\
&\leq &e^{\lambda \left( t+1\right) }\int_{D}\left( p_{D}^{a}\left(
t,x,y\right) \sqrt{\int_{D}\left[ p_{D}^{a}\left( 1,y,z\right) \right] ^{2}dz%
}\sqrt{\int \psi ^{2}\left( z\right) dz}\right) dy \\
&\leq &C\left( \Lambda ,D\right) e^{\lambda \left( t+1\right) }u^{a}\left(
t,x\right) \\
&=&\left[ C\left( \Lambda ,D\right) e^{\lambda _{\max }}\right] \times
e^{\lambda t}u^{a}\left( t,x\right)
\end{eqnarray*}%
and this gives the lower bound with $K_{1}=\psi \left( x\right) /\left[
C\left( \Lambda ,D\right) e^{\lambda _{\max }}\right] $. Clearly $\psi =\psi
_{1}^{a}$ depends on $a$ and a piori so does $K_{1}$. We now show that $\psi 
$ (and hence $K_{1}$) depends on $a$ only through $\Lambda $. From%
\begin{equation*}
p_{D}^{a}\left( t,y,y\right) =\sum_{i=1}^{\infty }e^{-\lambda
_{i}^{a}t}\left\vert \psi _{i}^{a}\left( y\right) \right\vert ^{2}
\end{equation*}%
evaluated at $t=1$ say we see that%
\begin{equation*}
\left\vert \psi \left( y\right) \right\vert ^{2}\leq e^{\lambda
}p_{D}^{a}\left( 1,y,y\right) \leq e^{\lambda _{\max }}p_{D}^{a}\left(
1,y,y\right) \leq e^{\lambda _{\max }}p^{a}\left( 1,y,y\right)
\end{equation*}%
and by using our upper heat kernel estimates for $p^{a}$ we see that there
is a constant $M=M\left( \Lambda ,D\right) $ such that $\left\vert \psi
\right\vert _{\infty }\leq M$. Given $x$ and $M$ we can find a compact set $%
\mathfrak{K}$ $\subset $ $D$ so that $m\left( D\backslash \mathfrak{K}%
\right) \leq 1/(4M^{2})$ and $x\in \mathfrak{K}$ (recall that $m$ is Haar
measure on $g^{2}\left( \mathbb{R}^{d}\right) $). By\ Harnack's inequality%
\begin{equation*}
\sup_{\mathfrak{K}}\psi \leq C\psi \left( x\right) .
\end{equation*}%
for $C=C\left( \mathfrak{K},\Lambda \right) =C\left( \,x,D,\Lambda \right) .$%
We then have%
\begin{equation*}
1=\left\vert \psi \right\vert _{L^{2}}\leq M\sqrt{m\left( D\backslash 
\mathfrak{K}\right) }+C\psi \left( x\right) \sqrt{m\left( \mathfrak{K}%
\right) }\leq 1/2+C\psi \left( x\right) \sqrt{m\left( D\right) }
\end{equation*}%
which gives the required lower bound on $\psi \left( x\right) \equiv \psi
_{1}^{a}\left( x\right) $ which only depends on $\,x,D$ and $\Lambda $. 
\newline
\underline{Upper bound:} Recall that $-\lambda \equiv -\lambda _{1}^{a}$
denotes the first eigenvalue of $L_{D}^{a}$ with associated semigroup $%
P_{D}^{a}$. It follows that 
\begin{equation*}
\left\vert P_{D}^{a}\left( t\right) f\right\vert _{L^{2}}\leq e^{-\lambda
t}\left\vert f\right\vert _{L^{2}}
\end{equation*}%
which may be rewritten as%
\begin{equation*}
\left\vert \int_{D}p_{D}^{a}\left( t,\cdot ,z\right) f\left( z\right)
dz\right\vert _{L^{2}}\leq e^{-\lambda t}\left\vert f\right\vert _{L^{2}}.
\end{equation*}%
Let $t>1$. Using Chapman-Kolmogorov and symmetry of the kernel,%
\begin{eqnarray*}
u\left( t,x\right) &=&\int_{D}p_{D}^{a}\left( t,x,z\right)
dz=\int_{D}\int_{D}p_{D}^{a}\left( 1,x,y\right) p_{D}^{a}\left(
t-1,z,y\right) dydz \\
&=&\sqrt{m\left( D\right) }\left( \int_{D}\left( \int_{D}p_{D}^{a}\left(
t-1,z,y\right) p_{D}^{a}\left( 1,x,y\right) dy\right) ^{2}dz\right) ^{1/2} \\
&=&\sqrt{m\left( D\right) }\left\vert \left( \int_{D}p_{D}^{a}\left(
t-1,\cdot ,y\right) p_{D}^{a}\left( 1,x,y\right) dy\right) \right\vert
_{L^{2}\left( D\right) } \\
&=&\sqrt{m\left( D\right) }\left\vert P_{D}^{a}\left( t-1\right)
\,p_{D}^{a}\left( 1,x,\cdot \right) \right\vert _{L^{2}\left( D\right) } \\
&\leq &\sqrt{m\left( D\right) }e^{-\lambda \left( t-1\right) }\,\left\vert
p_{D}^{a}\left( 1,x,\cdot \right) \right\vert _{L^{2}\left( D\right) } \\
&\leq &\sqrt{m\left( D\right) }e^{\lambda _{\max }}e^{-\lambda t}\sqrt{%
p_{D}^{a}\left( 2,x,x\right) } \\
&\leq &K_{2}e^{-\lambda t}.
\end{eqnarray*}%
where we used upper heat kernel estimates in the last step to obtain $%
K_{2}=K_{2}\left( D,\Lambda \right) .$
\end{proof}

\begin{corollary}
\label{SpectralCor}Fix $a\in \Xi \left( \Lambda \right) $. There exists $%
K=K\left( \Lambda \right) $ and for all $\varepsilon >0$ there exist $%
\lambda =\lambda ^{\left( \varepsilon \right) }$ such that%
\begin{eqnarray}
K^{-1}e^{-\lambda t\varepsilon ^{-2}} &\leq &\mathbb{P}^{a,0}\left[
\left\vert \left\vert \mathbf{X}\right\vert \right\vert _{0,\left[ 0,t\right]
}<\varepsilon \right]  \label{LemmaStr1} \\
\forall x &:&\mathbb{P}^{a,x}\left[ \left\vert \left\vert \mathbf{X}%
\right\vert \right\vert _{0,\left[ 0,t\right] }<\varepsilon \right] \leq
Ke^{-\lambda t\varepsilon ^{-2}}.  \label{LemmaStr2}
\end{eqnarray}
\end{corollary}

\begin{proof}
A straight-forward consequence of scaling and Proposition \ref{SpectralProp}
applied to%
\begin{equation*}
D=B\left( 0,1\right) =\left\{ y:\left\Vert y\right\Vert <1\right\}
\end{equation*}%
where $\left\Vert \cdot \right\Vert $ is the standard CC norm on $%
g^{2}\left( \mathbb{R}^{d}\right) $. Then $\lambda $ is the first eigenvalue
corresponding to $a$ scaled by factor $\varepsilon $.
\end{proof}

\begin{proposition}
\label{FLSestimate}\label{SmallTimeEstimate}\bigskip Let $\alpha \in \lbrack
0,1/2)$. There exists a constant $C_{\ref{SmallTimeEstimate}}$ such that for
all $\varepsilon \in (0,1]$ and $R>0$%
\begin{equation*}
\mathbb{P}^{a,0}\left( \left. \sup_{\left\vert t-s\right\vert <\varepsilon
^{2}}\frac{\left\Vert \mathbf{X}_{s,t}\right\Vert }{\left\vert
t-s\right\vert ^{\alpha }}>R\right\vert \left\Vert \mathbf{X}\right\Vert _{0;%
\left[ 0,1\right] }<\varepsilon \right) \leq C_{\ref{SmallTimeEstimate}}\exp
\left( -\frac{1}{C_{\ref{SmallTimeEstimate}}}\frac{R^{2}}{\varepsilon
^{2\left( 1-2\alpha \right) }}\right) .\text{ }
\end{equation*}
\end{proposition}

\begin{proof}
There will be no confusion to write $\mathbb{P}^{x}\equiv \mathbb{P}^{a,x}$
and $\mathbb{P}_{\varepsilon }^{x}$ $\equiv \mathbb{P}^{x}\left( \left.
\cdot \right\vert \left\Vert \mathbf{X}\right\Vert _{0;\left[ 0,1\right]
}<\varepsilon \right) $. Suppose there exists a pair of times $s,t\in \left[
0,1\right] $ such that%
\begin{equation*}
s<t,\,|t-s|<\varepsilon ^{2}\text{ and }\frac{\left\Vert \mathbf{X}%
_{s,t}\right\Vert }{\left\vert t-s\right\vert ^{\alpha }}>R.
\end{equation*}%
Then there exists a $k\in \{1,...,\left\lceil 1/\varepsilon ^{2}\right\rceil
\}$ so that $\left[ s,t\right] \subset \left[ \left( k-1\right) \varepsilon
^{2},\left( k+1\right) \varepsilon ^{2}\right] $. In particular, the
probability that such a pair of times exists is at most 
\begin{equation*}
\sum_{k=1}^{\left\lceil 1/\varepsilon ^{2}\right\rceil }\mathbb{P}%
_{\varepsilon }^{0}\left( \left\Vert \mathbf{X}\right\Vert _{\alpha ,\left[
\left( k-1\right) \varepsilon ^{2},\left( k+1\right) \varepsilon ^{2}\right]
}>R\right) .
\end{equation*}%
Set $\left[ \left( k-1\right) \varepsilon ^{2},\left( k+1\right) \varepsilon
^{2}\right] =:\left[ T_{1},T_{2}\right] $. The rest of the proof is
concerned with the existence of $C$ such that%
\begin{equation*}
\mathbb{P}_{\varepsilon }^{0}\left( \left\vert \left\vert \mathbf{X}%
\right\vert \right\vert _{\alpha ,\left[ T_{1},T_{2}\right] }>R\right) \leq
C\exp \left( -C^{-1}\frac{R^{2}}{^{\varepsilon ^{2\left( 1-2\alpha \right) }}%
}\right)
\end{equation*}%
since the factor $\left\lceil 1/\varepsilon ^{2}\right\rceil $ can be
absorbed in the exponential factor be making $C$ bigger. We estimate%
\begin{eqnarray*}
&&\mathbb{P}^{0}\left( \left\vert \left\vert \mathbf{X}\right\vert
\right\vert _{\alpha ,\left[ T_{1},T_{2}\right] }>R\left\vert \left\vert
\left\vert \mathbf{X}\right\vert \right\vert _{0,\left[ 0,1\right]
}<\varepsilon \right. \right) \\
&\leq &\frac{\mathbb{P}^{0}\left( \left\vert \left\vert \mathbf{X}%
\right\vert \right\vert _{\alpha ,\left[ T_{1},T_{2}\right] }>R;\left\vert
\left\vert \mathbf{X}\right\vert \right\vert _{0,\left[ 0,T_{1}\right]
}<\varepsilon ;\left\vert \left\vert \mathbf{X}\right\vert \right\vert _{0,%
\left[ T_{2},1\right] }<\varepsilon \right) }{\mathbb{P}^{0}\left[
\left\vert \left\vert \mathbf{X}\right\vert \right\vert _{0,\left[ 0,1\right]
}<\varepsilon \right] }.
\end{eqnarray*}%
By using the Markov-property and the above lemma, writing $\lambda ^{\left(
\varepsilon \right) }=\lambda ^{a;\varepsilon }$, this equals 
\begin{eqnarray*}
&&\frac{\mathbb{E}^{0}\left[ \mathbb{P}^{\mathbf{X}_{T_{2}}}\left(
\left\vert \left\vert \mathbf{X}\right\vert \right\vert _{0,\left[ 0,1-T_{2}%
\right] }<\varepsilon \right) ;\left\vert \left\vert \mathbf{X}\right\vert
\right\vert _{\alpha ,\left[ T_{1},T_{2}\right] }>R;\left\vert \left\vert 
\mathbf{X}\right\vert \right\vert _{0,\left[ 0,T_{1}\right] }<\varepsilon %
\right] }{\mathbb{P}^{0}\left[ \left\vert \left\vert \mathbf{X}\right\vert
\right\vert _{0,\left[ 0,1\right] }<\varepsilon \right] } \\
&\leq &Ce^{\lambda ^{\left( \varepsilon \right) }\varepsilon ^{-2}}\mathbb{E}%
^{0}\left[ e^{-\lambda ^{\left( \varepsilon \right) }\left( 1-T_{2}\right)
\varepsilon ^{-2}};\left\vert \left\vert \mathbf{X}\right\vert \right\vert
_{\alpha ,\left[ T_{1},T_{2}\right] }>R;\left\vert \left\vert \mathbf{X}%
\right\vert \right\vert _{0,\left[ 0,T_{1}\right] }<\varepsilon \right] \\
&=&Ce^{\lambda ^{\left( \varepsilon \right) }T_{2}\varepsilon ^{-2}}\mathbb{P%
}^{0}\left[ \left\vert \left\vert \mathbf{X}\right\vert \right\vert _{\alpha
,\left[ T_{1},T_{2}\right] }>R;\left\vert \left\vert \mathbf{X}\right\vert
\right\vert _{0,\left[ 0,T_{1}\right] }<\varepsilon \right]
\end{eqnarray*}%
where constants were allowed to change in insignificant ways. If $\mathbf{X}$
had indepedent increments in the group (such as is the case for Enhanced
Brownian motion $\mathbf{B}$) $\mathbb{P}^{0}\left[ ...\right] $ would split
up immediately. This is not the case here but the Markov property serves as
a substitute; using the Dirichlet heat kernel $p_{B\left( 0,\varepsilon
\right) }^{a}$we can write%
\begin{equation*}
\mathbb{P}^{0}\left[ \left\vert \left\vert \mathbf{X}\right\vert \right\vert
_{\alpha ,\left[ T_{1},T_{2}\right] }>R;\left\vert \left\vert \mathbf{X}%
\right\vert \right\vert _{0,\left[ 0,T_{1}\right] }<\varepsilon \right]
=\int_{B\left( 0,\varepsilon \right) }dx\,\,p_{B\left( 0,\varepsilon \right)
}^{a}\left( T_{1},0,x\right) \mathbb{P}^{x}\left[ \left\vert \left\vert 
\mathbf{X}\right\vert \right\vert _{\alpha ,\left[ 0,T_{2}-T_{1}\right] }>R%
\right] .
\end{equation*}%
Then, scaling and the usual Fernique-type estimates for the H\"{o}lder norm
of $\mathbf{X}$ gives 
\begin{equation*}
\sup_{x}\mathbb{P}^{x}\left[ \left\vert \left\vert \mathbf{X}\right\vert
\right\vert _{\alpha ,\left[ 0,T_{2}-T_{1}\right] }>R\right] \leq C\exp
\left( -\frac{1}{C}\left( \frac{R}{\varepsilon ^{1-2\alpha }}\right)
^{2}\right) ,\text{ \ }
\end{equation*}%
where we used $T_{2}-T_{1}=2\varepsilon ^{2}$, and we obtain 
\begin{eqnarray*}
&&\mathbb{P}^{0}\left[ \left\vert \left\vert \mathbf{X}\right\vert
\right\vert _{\alpha ,\left[ T_{1},T_{2}\right] }>R;\left\vert \left\vert 
\mathbf{X}\right\vert \right\vert _{0,\left[ 0,T_{1}\right] }<\varepsilon %
\right] \\
&\leq &C\exp \left( -\frac{1}{C}\left( \frac{R}{\varepsilon ^{1-2\alpha }}%
\right) ^{2}\right) \mathbb{P}^{0}\left[ \left\vert \left\vert \mathbf{X}%
\right\vert \right\vert _{0,\left[ 0,T_{1}\right] }<\varepsilon \right] \\
&\leq &C\exp \left( -\frac{1}{C}\left( \frac{R}{\varepsilon ^{1-2\alpha }}%
\right) ^{2}\right) e^{-\lambda ^{\left( \varepsilon \right)
}T_{1}\varepsilon ^{-2}}.
\end{eqnarray*}%
Putting things together we have%
\begin{eqnarray*}
\mathbb{P}^{0}\left( \left\vert \left\vert \mathbf{X}\right\vert \right\vert
_{\alpha ,\left[ T_{1},T_{2}\right] }>R\left\vert \left\vert \left\vert 
\mathbf{X}\right\vert \right\vert _{0,\left[ 0,1\right] }<\varepsilon
\right. \right) &\leq &Ce^{\lambda ^{\left( \varepsilon \right) }\left(
T_{2}-T_{1}\right) \varepsilon ^{-2}}\exp \left( -\frac{1}{C}\left( \frac{R}{%
\varepsilon ^{1-2\alpha }}\right) ^{2}\right) \\
&\leq &Ce^{2\lambda _{\max }}\exp \left( -\frac{1}{C}\left( \frac{R}{%
\varepsilon ^{1-2\alpha }}\right) ^{2}\right)
\end{eqnarray*}%
and the proof is finished.
\end{proof}

\begin{corollary}
\bigskip Let $\alpha \in \lbrack 0,1/2)$. For all $R>0$ the ball $\left\{ 
\mathbf{x}:\left\Vert \mathbf{x}\right\Vert _{\alpha \text{-H\"{o}l;}\left[
0,1\right] }<R\right\} $ has positive $\mathbb{P}^{a,0}$-measure and 
\begin{equation}
\lim_{\epsilon \rightarrow 0}\mathbb{P}^{a,0}\left( \left. \left\Vert 
\mathbf{X}\right\Vert _{\alpha \text{-H\"{o}l;}\left[ 0,1\right]
}<R\right\vert \left\Vert \mathbf{X}\right\Vert _{0;\left[ 0,1\right]
}<\varepsilon \right) \rightarrow 1.\text{ }  \label{CorSmallBall}
\end{equation}
\end{corollary}

\begin{proof}
We first observe that the uniform conditioning allows to localise the H\"{o}%
lder norm. More precisely, take $s<t$ in$\left[ 0,1\right] $ with $t-s\geq
\epsilon ^{2}$ and note that from $\left\Vert \mathbf{X}\right\Vert _{0;%
\left[ 0,1\right] }<\varepsilon $ we get $\left\Vert \mathbf{X}%
_{s,t}\right\Vert /\left\vert t-s\right\vert ^{\alpha }\leq \epsilon
^{1-2\alpha }$. It follows that for fixed $R$ and $\epsilon $ small enough,%
\begin{equation*}
\mathbb{P}^{a,0}\left( \left. \left\Vert \mathbf{X}\right\Vert _{\alpha 
\text{-H\"{o}l;}\left[ 0,1\right] }\geq R\right\vert \left\Vert \mathbf{X}%
\right\Vert _{0;\left[ 0,1\right] }<\varepsilon \right) =\mathbb{P}%
^{a,0}\left( \left. \sup_{\left\vert t-s\right\vert <\varepsilon ^{2}}\frac{%
\left\Vert \mathbf{X}_{s,t}\right\Vert }{\left\vert t-s\right\vert ^{\alpha }%
}\geq R\right\vert \left\Vert \mathbf{X}\right\Vert _{0;\left[ 0,1\right]
}<\varepsilon \right)
\end{equation*}%
and the preceding proposition shows convergence to zero with $\epsilon $ and
(\ref{CorSmallBall}) follows. Finally,%
\begin{eqnarray*}
\mathbb{P}^{a,0}\left( \left\Vert \mathbf{X}\right\Vert _{\alpha \text{-H%
\"{o}l;}\left[ 0,1\right] }<R\right) &\geq &\mathbb{P}^{a,0}\left( \left.
\left\Vert \mathbf{X}\right\Vert _{\alpha \text{-H\"{o}l;}\left[ 0,1\right]
}<R\right\vert \left\Vert \mathbf{X}\right\Vert _{0;\left[ 0,1\right]
}<\varepsilon \right) \times \mathbb{P}^{a,0}\left( \left\Vert \mathbf{X}%
\right\Vert _{0;\left[ 0,1\right] }<\varepsilon \right) \\
&\geq &\mathbb{P}^{a,0}\left( \left\Vert \mathbf{X}\right\Vert _{0;\left[ 0,1%
\right] }<\varepsilon \right) /2\text{ \ \ \ \ (for }\epsilon \text{ small
enough)}
\end{eqnarray*}%
and this is positive by either Proposition \ref{SpectralCor} or Proposition %
\ref{subLower}.
\end{proof}

\begin{corollary}
Let $Y=\pi \left( 0,y_{0};\mathbf{X}^{a;0}\right) \equiv \pi \left( \mathbf{X%
}\right) $ denote the $\mathbb{R}^{e}$-valued (random) RDE solution driven
by $\mathbf{X}^{a,x}$ along fixed $\mathrm{Lip}^{2+\varepsilon }$ vector
fields $V_{1},...,V_{d}$ on $\mathbb{R}^{e}$ and started at time $0$ from $%
y_{0}$ fixed. Then, for any $R>0$,%
\begin{equation*}
\mathbb{P}^{a,0}\left( \left. \left\vert Y\right\vert _{\alpha \text{-H\"{o}%
l;}\left[ 0,1\right] }>R\right\vert \left\Vert \mathbf{X}\right\Vert _{0;%
\left[ 0,1\right] }<\varepsilon \right) \rightarrow 0\text{ with }\epsilon
\rightarrow 0.
\end{equation*}
\end{corollary}

\begin{proof}
From, Lyons' limit theorem, $\left\Vert \mathbf{X}\right\Vert _{\alpha \text{%
-H\"{o}l;}\left[ 0,1\right] }\rightarrow 0$ implies, deterministically, $%
\left\vert Y\right\vert _{\alpha \text{-H\"{o}l;}\left[ 0,1\right]
}\rightarrow 0$.
\end{proof}

\subsection{The Stroock-Varadhan support theorem for Markov RDEs}

Let $h\in C_{0}^{1-var}\left( \left[ 0,1\right] ,\mathbb{R}^{d}\right) $.
Give a uniformly elliptic $a:\mathbb{R}^{d}\rightarrow \mathbb{R}^{d}\otimes 
\mathbb{R}^{d}$, so that $a\circ \pi _{1}\in \Xi \left( \Lambda \right) $ we
know that $T_{h}\left( \mathbf{X}^{a}\right) $ is Markov. Furthermore, $%
\mathbf{X}_{0,\cdot }\mathbf{=X}_{0,\cdot }^{a}\in C_{0}^{\gamma \text{-H%
\"{o}l}}\left( \left[ 0,1\right] ,g^{2}\left( \mathbb{R}^{d}\right) \right) $
for $\gamma \in (1/3,1/2)$ and from basic facts of the translation operator
we also have%
\begin{equation*}
T_{h}\left( \mathbf{X}\right) \in C_{0}^{\gamma \text{-H\"{o}l}}\left( \left[
0,1\right] ,g^{2}\left( \mathbb{R}^{d}\right) \right)
\end{equation*}%
This "step-$2$" $\gamma $-H\"{o}lder rough path lifts uniquely and
continuously to any step-$N$ rough path%
\begin{equation*}
S_{N}\left( T_{h}\left( \mathbf{X}\right) \right) \in C_{0}^{\gamma \text{-H%
\"{o}l}}\left( \left[ 0,1\right] ,g^{N}\left( \mathbb{R}^{d}\right) \right) .
\end{equation*}%
Obviously, specializing to~$h=0$ and it is clear that $S_{N}\left( \mathbf{X}%
\right) $ is also $\alpha $-H\"{o}lder for $1/\left( N+1\right) <\alpha <1/N$
and thus a "step-$N$" $\alpha $-rough path in its own right. By basic
consistency properties of rough differential equations, the solutions
corresponding to driving $S_{N}\left( \mathbf{X}\right) $, as step-$N$ rough
path, and $\mathbf{X}$ as step-$2$ rough path, coincide. Hence, it is good
enough to obtain a support description for $S_{N}\left( \mathbf{X}\right) $
in $\alpha $-H\"{o}lder topology and we are able to do this with $N=6$ and
any $\alpha <1/6$.

\begin{lemma}
\label{ControlSNyOverTSsmallerThanDeltaSquare}Let $\mathbf{y}\in
C_{0}^{\gamma \text{-H\"{o}l}}\left( \left[ 0,1\right] ,g^{2}\left( \mathbb{R%
}^{d}\right) \right) ,\,\,\,\gamma \in \left( 1/3,1/2\right) $. For every
integer $N$, there exists $K_{N}$ such that for all $s<t$ in $\left[ 0,1%
\right] $,%
\begin{equation*}
\left\Vert S_{N}\left( \mathbf{y}\right) \right\Vert _{\gamma \text{-H\"{o}l;%
}\left[ s,t\right] }\leq K_{N}\left\Vert \mathbf{y}\right\Vert _{\gamma 
\text{-H\"{o}l;}\left[ s,t\right] }
\end{equation*}%
(Notice that the respective $\gamma $-H\"{o}lder "norms" are with respect to 
$g^{N}\left( \mathbb{R}^{d}\right) $-valued paths on the left-hand-side and
with respect to $g^{2}\left( \mathbb{R}^{d}\right) $-valued paths on the
right-hand-side.)
\end{lemma}

\begin{proof}
See \cite[p242]{lyons-98}.
\end{proof}

\begin{proposition}
\label{ThXvsX}Let $\gamma \in \lbrack 0,1/2)$. Let $h\in C_{0}^{1\text{-var}%
}\left( \left[ 0,1\right] ,\mathbb{R}^{d}\right) $ and $\mathbf{y}\in
C_{0}^{\gamma \text{-H\"{o}l}}\left( \left[ 0,1\right] ,g^{2}\left( \mathbb{R%
}^{d}\right) \right) .$ Then, there exists a constant $C_{\ref{ThXvsX}}>0$
such that for all $0\leq s<t\leq 1$ ,%
\begin{equation}
\left\Vert T_{h}\left( \mathbf{y}\right) _{s,t}\right\Vert \leq C_{\ref%
{ThXvsX}}\left( \left\Vert \mathbf{y}\right\Vert _{\gamma \text{-H\"{o}l;}%
\left[ s,t\right] }\left\vert t-s\right\vert ^{\gamma }+\left\vert
h\right\vert _{1-var,\left[ s,t\right] }\right) .  \label{estimateThCrude}
\end{equation}%
In particular, if $h\in \mathcal{H}$, the usual Cameron-Martin space with $%
\left\vert h\right\vert _{\mathcal{H}}\equiv \left\vert \dot{h}\right\vert
_{L^{2}\left( \left[ 0,1\right] ,\mathbb{R}^{d}\right) },$ then%
\begin{equation}
\left\Vert T_{h}\left( \mathbf{y}\right) \right\Vert _{\gamma \text{-H\"{o}l;%
}\left[ s,t\right] }\leq C_{\ref{ThXvsX}}\left( \left\Vert \mathbf{y}%
\right\Vert _{\gamma \text{-H\"{o}l;}\left[ s,t\right] }+\left\vert
h\right\vert _{\mathcal{H}}\left\vert t-s\right\vert ^{\frac{1}{2}-\gamma
}\right) .  \label{estimteThCameronMartin}
\end{equation}
\end{proposition}

\begin{proof}
It is easy to see that $\left\Vert T_{h}\left( \mathbf{y}\right)
_{s,t}\right\Vert $ is less equal than a constant times%
\begin{equation*}
\left\vert h_{s,t}+\mathbf{y}_{s,t}^{1}\right\vert +\sqrt{\left\vert \pi
_{2}\left( \mathbf{y}_{s,t}\right) \right\vert }+\sqrt{\left\vert
\int_{s}^{t}h_{s,r}\otimes dh_{r}\right\vert }+\sqrt{\left\vert
\int_{s}^{t}h_{s,r}\otimes dy_{r}\right\vert }+\sqrt{\left\vert
\int_{s}^{t}y_{s,r}\otimes dh_{r}\right\vert }.
\end{equation*}
The first three summand are easy to estimate.\ To deal with the last two it
suffices to note%
\begin{eqnarray*}
\sqrt{\left\vert \int_{s}^{t}y_{s,r}\otimes dh_{r}\right\vert } &\leq &\sqrt{%
\sup_{r\in \left[ s,t\right] }\left\vert y_{s,r}\right\vert .\left\vert
h\right\vert _{1\text{-var},\left[ s,t\right] }} \\
&\leq &\frac{1}{2}\sup_{r\in \left[ s,t\right] }\left\vert
y_{s,r}\right\vert +\frac{1}{2}\left\vert h\right\vert _{1\text{-var},\left[
s,t\right] },
\end{eqnarray*}%
then use integration by parts for the last summand. Finally, $\left\vert
h\right\vert _{1\text{-var};\left[ s,t\right] }\leq \left\vert h\right\vert
_{\mathcal{H}}\left\vert t-s\right\vert ^{\frac{1}{2}}$ implies (\ref%
{estimteThCameronMartin}).
\end{proof}

Let us remark that Proposition \ref{subLower} and remains valid with
identical proof in the step-$N$ setting. (The toolbox of Dirchlet forms
applies immediately with $g^{N}\left( \mathbb{R}^{d}\right) $ instead of $%
g^{2}\left( \mathbb{R}^{d}\right) $. Constants may depend on $N$, but $N=6$
will suffice for us.

\begin{theorem}
\label{SupportXforAlphaLToneQuarter}Let $h$ be a Lipschitz path and $\alpha
\in \lbrack 0,1/6)$. Then, for all $\varepsilon >0,$%
\begin{equation*}
\mathbb{P}^{a}\left( \left\Vert S_{6}\left( T_{h}\left( \mathbf{X}\right)
\right) \right\Vert _{\alpha \text{-H\"{o}l};\left[ 0,1\right] }<\varepsilon
\right) >0.
\end{equation*}
\end{theorem}

\begin{proof}
By take $\alpha $ close enough to $1/6$ we may assyme that $N=\left[
1/\alpha \right] =6$. We shall choose a (good) H\"{o}lder exponent $\gamma
=\gamma \left( \alpha \right) \in \left( 1/3,1/2\right) $, to be chosen
below ($\gamma =1/3+\left( 1/6-\alpha \right) /2$ will do). For any $p>0$
(to be choose large later on),%
\begin{eqnarray*}
&&\mathbb{P}\left( \left\Vert S_{N}\left( T_{h}\left( \mathbf{X}\right)
\right) \right\Vert _{\alpha \text{-H\"{o}l};\left[ 0,1\right] }>\varepsilon
\right) \\
&\leq &\mathbb{P}\left( \sup_{\left\vert t-s\right\vert <\varepsilon ^{p}}%
\frac{\left\Vert S_{N}\left( T_{h}\left( \mathbf{X}\right) \right)
_{s,t}\right\Vert }{\left\vert t-s\right\vert ^{\alpha }}>\varepsilon
\right) +\mathbb{P}\left( \text{ }\sup_{\left\vert t-s\right\vert \geq
\varepsilon ^{p}}\frac{\left\Vert S_{N}\left( T_{h}\left( \mathbf{X}\right)
\right) _{s,t}\right\Vert }{\left\vert t-s\right\vert ^{\alpha }}%
>\varepsilon \right) \\
&\leq &\mathbb{P}\left( \sup_{\left\vert t-s\right\vert <\varepsilon ^{p}}%
\frac{\left\Vert S_{N}\left( T_{h}\left( \mathbf{X}\right) \right)
_{s,t}\right\Vert }{\left\vert t-s\right\vert ^{\gamma }}>\varepsilon
\right) +\text{ }\mathbb{P}\left( \text{ }\sup_{\left\vert t-s\right\vert
\geq \varepsilon ^{p}}\frac{\left\Vert S_{N}\left( T_{h}\left( \mathbf{X}%
\right) \right) _{s,t}\right\Vert }{\left\vert t-s\right\vert ^{\alpha }}%
>\varepsilon \right) .
\end{eqnarray*}%
Using $\mathbb{P}\left( A\right) \leq \mathbb{P}\left( B\right) +\mathbb{P}%
\left( C\right) \implies \mathbb{P}\left( A^{c}\right) \geq -\mathbb{P}%
\left( B\right) +\mathbb{P}\left( C^{c}\right) $ we see that%
\begin{eqnarray*}
\mathbb{P}\left( \left\Vert S_{N}\left( T_{h}\left( \mathbf{X}\right)
\right) \right\Vert _{\alpha \text{-H\"{o}l};\left[ 0,1\right] }\leq
\varepsilon \right) &\geq &-\mathbb{P}\left( \sup_{\left\vert t-s\right\vert
<\varepsilon ^{p}}\frac{\left\Vert S_{N}\left( T_{h}\left( \mathbf{X}\right)
\right) _{s,t}\right\Vert }{\left\vert t-s\right\vert ^{\gamma }}%
>\varepsilon \right) \\
&&+\mathbb{P}\left( \text{ }\left\Vert S_{N}\left( T_{h}\left( \mathbf{X}%
\right) \right) \right\Vert _{\infty ;\left[ 0,1\right] }\leq \varepsilon
.\left( \varepsilon ^{p}\right) ^{\alpha }\right) \\
&\equiv &-\left( I\right) +\left( II\right) .
\end{eqnarray*}%
and the proof will be finished if we can find $p=p\left( \alpha \right) $
such that $\left( I\right) /(II)\rightarrow 0$ as $\varepsilon \rightarrow
0. $ It follows from Proposition \ref{subLower} that as $\varepsilon
\rightarrow 0$, 
\begin{equation*}
\left( II\right) \geq \frac{1}{c_{1}}\exp \left( -c_{1}\left( \varepsilon
^{1+p\alpha }\right) ^{-2}\right) =\frac{1}{c_{1}}\exp \left( -c_{1}\left( 
\frac{1}{\varepsilon }\right) ^{2+2p\alpha }\right)
\end{equation*}%
which we express without the irrelevant positive constants as%
\begin{equation}
\log \left( II\right) \gtrsim -\left( \frac{1}{\varepsilon }\right)
^{2+2p\alpha }.  \label{logII_lowerBd}
\end{equation}%
At the same time,%
\begin{eqnarray*}
\left( I\right) &=&\mathbb{P}\left( \sup_{\left\vert t-s\right\vert
<\varepsilon ^{p}}\frac{\left\Vert S_{N}\left( T_{h}\left( \mathbf{X}\right)
\right) _{s,t}\right\Vert }{\left\vert t-s\right\vert ^{\gamma }}%
>\varepsilon \right) \\
&\leq &\sum_{k=1}^{\left\lceil 1/\varepsilon ^{p}\right\rceil }\mathbb{P}%
\left( \left\Vert S_{N}\left( T_{h}\left( \mathbf{X}\right) \right)
\right\Vert _{\gamma \text{-H\"{o}l;}\left[ \left( k-1\right) \varepsilon
^{p},\left( k+1\right) \varepsilon ^{p}\right] }>\varepsilon \right) \\
&\leq &\sum_{k=1}^{\left\lceil 1/\varepsilon ^{p}\right\rceil }\mathbb{P}%
\left( \left\Vert T_{h}\left( \mathbf{X}\right) \right\Vert _{\gamma \text{-H%
\"{o}l;}\left[ \left( k-1\right) \varepsilon ^{p},\left( k+1\right)
\varepsilon ^{p}\right] }>\frac{\varepsilon }{c_{2}}\right) \text{ }
\end{eqnarray*}%
where $c_{2}=K_{N}$ is the constant from Lemma \ref%
{ControlSNyOverTSsmallerThanDeltaSquare}. (Here we used $\gamma >1/3$.) By
Proposition \ref{ThXvsX} this estimate continues with%
\begin{eqnarray*}
&\leq &\sum_{k=1}^{\left\lceil 1/\varepsilon ^{p}\right\rceil }\mathbb{P}%
\left( \left\Vert \mathbf{X}\right\Vert _{\gamma \text{-H\"{o}l;}\left[
\left( k-1\right) \varepsilon ^{p},\left( k+1\right) \varepsilon ^{p}\right]
}+\underbrace{\left\vert h\right\vert _{\mathcal{H}}\left( 2\varepsilon
^{p}\right) ^{\frac{1}{2}-\gamma }}>\frac{\varepsilon }{c_{3}}\right) \\
&\leq &\sum_{k=1}^{\left\lceil 1/\varepsilon ^{p}\right\rceil }\mathbb{P}%
\left( \left\Vert \mathbf{X}\right\Vert _{\gamma \text{-H\"{o}l;}\left[
\left( k-1\right) \varepsilon ^{p},\left( k+1\right) \varepsilon ^{p}\right]
}>\frac{\varepsilon }{c_{4}}\right)
\end{eqnarray*}%
where the term indicated by the curley bracket can indeed by omitted as $%
\varepsilon \rightarrow 0$ provided $p$ is chosen large enough so that $%
p\left( 1/2-\gamma \right) >1$. With scaling and Fernique estimates we see
that%
\begin{equation*}
\sum_{k=1}^{\left\lceil 1/\varepsilon ^{p}\right\rceil }\mathbb{P}\left(
\left\Vert \mathbf{X}\right\Vert _{\gamma \text{-H\"{o}l;}\left[ \left(
k-1\right) \varepsilon ^{p},\left( k+1\right) \varepsilon ^{p}\right] }>%
\frac{\varepsilon }{c_{4}}\right) \leq c_{5}\varepsilon ^{p}\exp \left( -%
\frac{1}{c_{5}}\left( \frac{\varepsilon }{\left( \varepsilon ^{p}\right)
^{1/2-\gamma }}\right) ^{2}\right)
\end{equation*}%
Focusing on the decay rate of $\left( I\right) $ and again ignoring
irrelevant positive constants, we see that%
\begin{equation*}
\log \left( I\right) \lesssim -\left( \frac{\varepsilon }{\left( \varepsilon
^{p}\right) ^{1/2-\gamma }}\right) ^{2}=-\left( \frac{1}{\varepsilon }%
\right) ^{-2+p\left( 1-2\gamma \right) }.
\end{equation*}%
Recalling $\log \left( II\right) \gtrsim -\left( 1/\varepsilon \right)
^{2+p\left( 2\alpha \right) }$ it is clear that, by choosing $p$ large
enough, $\left( I\right) /\left( II\right) \rightarrow 0$ as $\varepsilon
\rightarrow 0$ \textit{provided that }$1-2\gamma >2\alpha $. Our only
constraint is $\gamma >1/3$ and we now see that this is precisely possible
when $\alpha <1/6$ and so the proof is finished.
\end{proof}

\begin{corollary}
The support of the law of $S_{6}\left( \mathbf{X}_{0,\cdot }^{a,0}\right) $
in $\alpha $-H\"{o}lder topology, $\alpha \in \lbrack 0,1/6)$, equals $%
C_{0}^{0,\alpha \text{-H\"{o}lder}}\left( \left[ 0,1\right] ,g^{6}\left( 
\mathbb{R}^{d}\right) \right) $.
\end{corollary}

\begin{proof}
Given $\alpha \in \left( 0,1\right) ,N=\left[ 1/\alpha \right] $, a fixed
Lipschitz $h$ and $\mathbf{x}\in C_{0}^{0,\alpha \text{-H\"{o}l}}\left( %
\left[ 0,1\right] ,g^{N}\left( \mathbb{R}^{d}\right) \right) $ we know \cite[%
p57]{lyons-qian-02} that%
\begin{equation*}
\mathbf{x}\mapsto T_{h}\left( \mathbf{x}\right)
\end{equation*}%
is continuous under $d_{\alpha \text{-H\"{o}l}}$ on the pathspace $%
C_{0}^{0,\alpha \text{-H\"{o}l}}$. It then easily follows that 
\begin{equation*}
T_{-h}\left( \mathbf{x}^{n}\right) \rightarrow 0\text{ \ }%
\Longleftrightarrow \text{ \ }\mathbf{x}^{n}\rightarrow S_{N}\left( h\right)
.
\end{equation*}%
Indeed, $"\Longleftarrow "$ comes from continuity of $\mathbf{x}\mapsto
T_{-h}\left( \mathbf{x}\right) $ and $T_{-h}\left( S_{N}\left( h\right)
\right) =S_{N}\left( h-h\right) =0$ while $"\Longrightarrow "$ follows from%
\begin{equation*}
\underset{=\mathbf{x}^{n}}{\underbrace{T_{h}T_{-h}\left( \mathbf{x}%
^{n}\right) }}\rightarrow \underset{S_{N}\left( H\right) }{\underbrace{%
T_{h}\left( 0\right) }}.
\end{equation*}%
Then use Theorem \ref{SupportXforAlphaLToneQuarter}.
\end{proof}

\begin{corollary}[Stroock-Varadhan]
Let $Y=\pi \left( 0,y_{0};\mathbf{X}^{a,x}\right) \equiv \pi \left( \mathbf{X%
}^{a,x}\right) $ denote the $\mathbb{R}^{e}$-valued (random) RDE solution
driven by $\mathbf{X}^{a,x}$ along fixed $\mathrm{Lip}^{6+\varepsilon }$
vector fields $V_{1},...,V_{d}$ on $\mathbb{R}^{e}$ and started at time $0$
from $y_{0}$ fixed. Let $\mathbb{Q}$ denote the law of $\left( Y_{t}:0\leq
t\leq 1\right) $. Then the support of $\mathbb{Q}$ in uniform topology is
the closure of all control ODE solution,%
\begin{equation*}
\mathcal{S=}\left\{ \pi \left( 0,y_{0},h\right) :h\in C^{1}\left( \left[ 0,1%
\right] ,\mathbb{R}^{d}\right) \right\} .
\end{equation*}%
Here $y\equiv \pi \left( 0,y_{0},h\right) $ denotes the unique solution,
started at time $0$ from $y_{0}$, of the ordinary differential equation%
\begin{equation*}
dy=\sum_{i=1}^{d}V_{i}\left( y\right) dh^{i}.
\end{equation*}
\end{corollary}

\begin{proof}
$Y$ is obtained as RDE solution driven by a $\mathbf{X}^{a,0}$. By a basic
consistency properties of RDE solutions, it is \textit{also }the RDE
solution driven by $S_{6}\left( \mathbf{X}^{a,0}\right) $. \ By continuity
of the It\^{o}-Lyons map, the support description of the later implies the
Stroock-Varadhan support description for $Y.$
\end{proof}

\begin{acknowledgement}
The authors would like to thank T.J. Lyons, S.R.S. Varadhan and J.R. Norris
for helpful discussions. The first author is grateful to T. Coulhon, A.A.
Grigor'yan and E.B. Davies for conversations during the writing of Section
8.2.
\end{acknowledgement}

\bibliographystyle{plain}
\bibliography{roughpaths}

\def\cprime{$'$} \def\cprime{$'$}
\begin{thebibliography}{10}

\bibitem{Ar67}
D.~G. Aronson.
\newblock Bounds for the fundamental solution of a parabolic equation.
\newblock {\em Bull. Amer. Math. Soc.}, 73:890--896, 1967.

\bibitem{BaKu00}
Richard~F. Bass and Takashi Kumagai.
\newblock Laws of the iterated logarithm for some symmetric diffusion
  processes.
\newblock {\em Osaka J. Math.}, 37(3):625--650, 2000.

\bibitem{BeGrLe94}
G{\'e}rard Ben~Arous, Mihai Gr{\u{a}}dinaru, and Michel Ledoux.
\newblock H\"older norms and the support theorem for diffusions.
\newblock {\em Ann. Inst. H. Poincar\'e Probab. Statist.}, 30(3):415--436,
  1994.

\bibitem{CKS87}
E.~A. Carlen, S.~Kusuoka, and D.~W. Stroock.
\newblock Upper bounds for symmetric {M}arkov transition functions.
\newblock {\em Ann. Inst. H. Poincar\'e Probab. Statist.}, 23(2,
  suppl.):245--287, 1987.

\bibitem{Dav89}
E.~B. Davies.
\newblock {\em Heat kernels and spectral theory}, volume~92 of {\em Cambridge
  Tracts in Mathematics}.
\newblock Cambridge University Press, Cambridge, 1989.

\bibitem{DeZe93}
Amir Dembo and Ofer Zeitouni.
\newblock {\em Large deviations techniques and applications}, volume~38 of {\em
  Applications of Mathematics (New York)}.
\newblock Springer-Verlag, New York, second edition, 1998.

\bibitem{DeuSt89}
Jean-Dominique Deuschel and Daniel~W. Stroock.
\newblock {\em Large deviations}, volume 137 of {\em Pure and Applied
  Mathematics}.
\newblock Academic Press Inc., Boston, MA, 1989.

\bibitem{FaSt86}
E.~B. Fabes and D.~W. Stroock.
\newblock A new proof of {M}oser's parabolic {H}arnack inequality using the old
  ideas of {N}ash.
\newblock {\em Arch. Rational Mech. Anal.}, 96(4):327--338, 1986.

\bibitem{friz-lyons-stroock-06}
P.~Friz, T.~Lyons, and D.~Stroock.
\newblock L\'evy's area under conditioning.
\newblock {\em Ann. Inst. H. Poincar\'e Probab. Statist.}, 42(1):89--101, 2006.

\bibitem{friz-victoir-04-Note}
Peter Friz and Nicolas Victoir.
\newblock A note on the notion of geometric rough paths.
\newblock {\em Probab. Theory Related Fields}, 136:395--416, 2006.

\bibitem{friz-victoir-06-Euler}
Peter Friz and Nicolas Victoir.
\newblock Euler estimates for rough differential equations.
\newblock {\em Accepted, Journal of Differential Equations}, 2007.

\bibitem{Fu94}
Masatoshi Fukushima, Y{\=o}ichi {\=O}shima, and Masayoshi Takeda.
\newblock {\em Dirichlet forms and symmetric {M}arkov processes}, volume~19 of
  {\em de Gruyter Studies in Mathematics}.
\newblock Walter de Gruyter \& Co., Berlin, 1994.

\bibitem{ikeda-watanabe-89}
Nobuyuki Ikeda and Shinzo Watanabe.
\newblock {\em Stochastic differential equations and diffusion processes}.
\newblock North-Holland Publishing Co., Amsterdam, second edition, 1989.

\bibitem{Je86}
David Jerison.
\newblock The {P}oincar\'e inequality for vector fields satisfying
  {H}\"ormander's condition.
\newblock {\em Duke Math. J.}, 53(2):503--523, 1986.

\bibitem{LejI06}
Antoine Lejay.
\newblock Stochastic differential equations driven by processes generated by
  divergence form operators {I}: a {W}ong-{Z}akai theorem.
\newblock {\em ESAIM Prob. and Stat.}, 10:356--379, 2006.

\bibitem{LejII06}
Antoine Lejay.
\newblock Stochastic differential equations driven by processes generated by
  divergence form operators {II}: Convergence results.
\newblock {\em Accepted, ESAIM Prob. and Stat.}, 2007.

\bibitem{lyons-98}
Terry Lyons.
\newblock Differential equations driven by rough signals.
\newblock {\em Rev. Mat. Iberoamericana}, 14(2):215--310, 1998.

\bibitem{lyons-04}
Terry Lyons.
\newblock St. {F}lour {L}ectures on {R}ough {P}aths, 2004.
\newblock Handwritten notes available at sag.maths.ox.ac.uk/tlyons/st-flour/.

\bibitem{lyons-qian-02}
Terry Lyons and Zhongmin Qian.
\newblock {\em System {C}ontrol and {R}ough {P}aths}.
\newblock Oxford University Press, 2002.
\newblock Oxford Mathematical Monographs.

\bibitem{LySt99}
Terry Lyons and Lucre{\c{t}}iu Stoica.
\newblock The limits of stochastic integrals of differential forms.
\newblock {\em Ann. Probab.}, 27(1):1--49, 1999.

\bibitem{MiSaSo94}
Annie Millet and Marta Sanz-Sol{\'e}.
\newblock A simple proof of the support theorem for diffusion processes.
\newblock In {\em S\'eminaire de Probabilit\'es, XXVIII}, volume 1583 of {\em
  Lecture Notes in Math.}, pages 36--48. Springer, Berlin, 1994.

\bibitem{montgomery-02}
Richard Montgomery.
\newblock {\em A tour of subriemannian geometries, their geodesics and
  applications}, volume~91 of {\em Mathematical Surveys and Monographs}.
\newblock American Mathematical Society, Providence, RI, 2002.

\bibitem{Mo71}
J.~Moser.
\newblock On a pointwise estimate for parabolic differential equations.
\newblock {\em Comm. Pure Appl. Math.}, 24:727--740, 1971.

\bibitem{Mo64}
J{\"u}rgen Moser.
\newblock A {H}arnack inequality for parabolic differential equations.
\newblock {\em Comm. Pure Appl. Math.}, 17:101--134, 1964.

\bibitem{Ra01}
Jos{\'e}~A. Ram{\'{\i}}rez.
\newblock Short-time asymptotics in {D}irichlet spaces.
\newblock {\em Comm. Pure Appl. Math.}, 54(3):259--293, 2001.

\bibitem{Rob91}
Derek~W. Robinson.
\newblock {\em Elliptic operators and {L}ie groups}.
\newblock Oxford Mathematical Monographs. The Clarendon Press Oxford University
  Press, New York, 1991.
\newblock Oxford Science Publications.

\bibitem{SaSt91}
L.~Saloff-Coste and D.~W. Stroock.
\newblock Op\'erateurs uniform\'ement sous-elliptiques sur les groupes de
  {L}ie.
\newblock {\em J. Funct. Anal.}, 98(1):97--121, 1991.

\bibitem{St88}
Daniel~W. Stroock.
\newblock Diffusion semigroups corresponding to uniformly elliptic divergence
  form operators.
\newblock In {\em S\'eminaire de Probabilit\'es, XXII}, volume 1321 of {\em
  Lecture Notes in Math.}, pages 316--347. Springer, Berlin, 1988.

\bibitem{StVa72}
Daniel~W. Stroock and S.~R.~S. Varadhan.
\newblock On the support of diffusion processes with applications to the strong
  maximum principle.
\newblock In {\em Proceedings of the Sixth Berkeley Symposium on Mathematical
  Statistics and Probability (Univ. California, Berkeley, Calif., 1970/1971),
  Vol. III: Probability theory}, pages 333--359, Berkeley, Calif., 1972. Univ.
  California Press.

\bibitem{SturmIII}
K.~T. Sturm.
\newblock Analysis on local {D}irichlet spaces. {III}. {T}he parabolic
  {H}arnack inequality.
\newblock {\em J. Math. Pures Appl. (9)}, 75(3):273--297, 1996.

\bibitem{SturmII}
Karl-Theodor Sturm.
\newblock Analysis on local {D}irichlet spaces. {II}. {U}pper {G}aussian
  estimates for the fundamental solutions of parabolic equations.
\newblock {\em Osaka J. Math.}, 32(2):275--312, 1995.

\bibitem{Sturm93}
Karl-Theodor Sturm.
\newblock On the geometry defined by {D}irichlet forms.
\newblock In {\em Seminar on Stochastic Analysis, Random Fields and
  Applications (Ascona, 1993)}, volume~36 of {\em Progr. Probab.}, pages
  231--242. Birkh\"auser, Basel, 1995.

\bibitem{ElRo00}
A.~F.~M. ter Elst and Derek~W. Robinson.
\newblock Second-order subelliptic operators on {L}ie groups. {II}. {R}eal
  measurable principal coefficients.
\newblock In {\em Semigroups of operators: theory and applications (Newport
  Beach, CA, 1998)}, volume~42 of {\em Progr. Nonlinear Differential Equations
  Appl.}, pages 103--124. Birkh\"auser, Basel, 2000.

\bibitem{VaII67}
S.~R.~S. Varadhan.
\newblock Diffusion processes in a small time interval.
\newblock {\em Comm. Pure Appl. Math.}, 20:659--685, 1967.

\bibitem{Va80}
S.~R.~S. Varadhan.
\newblock {\em Lectures on diffusion problems and partial differential
  equations}, volume~64 of {\em Tata Institute of Fundamental Research Lectures
  on Mathematics and Physics}.
\newblock Tata Institute of Fundamental Research, Bombay, 1980.
\newblock With notes by Pl. Muthuramalingam and Tara R. Nanda.

\bibitem{Varo90}
N.~Th. Varopoulos.
\newblock Small time {G}aussian estimates of heat diffusion kernels. {II}.
  {T}he theory of large deviations.
\newblock {\em J. Funct. Anal.}, 93(1):1--33, 1990.

\bibitem{CoSaVa92}
N.~Th. Varopoulos, L.~Saloff-Coste, and T.~Coulhon.
\newblock {\em Analysis and geometry on groups}, volume 100 of {\em Cambridge
  Tracts in Mathematics}.
\newblock Cambridge University Press, Cambridge, 1992.

\end{thebibliography}

\bigskip

\end{document}